\documentclass{article}
\usepackage[utf8]{inputenc}

\usepackage{epstopdf}

\makeatletter
\renewcommand{\@biblabel}[1]{#1\hfill \hspace{-0.2cm}}
\makeatother

\usepackage{amsmath}
\usepackage{amssymb}
\newcommand\Peclet{\mathbb{P}\mathrm{e}}  
\usepackage{bm}
\usepackage{graphicx}
\usepackage{caption}
\usepackage{subcaption}
\usepackage{multirow}
\usepackage{todonotes}
\usepackage{placeins}
\usepackage{epstopdf} 
\usepackage[pagewise]{lineno}
\usepackage{float}
\usepackage{url}

\begin{document}

\title{A Machine Learning approach to enhance the SUPG stabilization method for advection-dominated differential problems}
\author{Tommaso Tassi$^1$ \and Alberto Zingaro$^{2,*}$ \and Luca Dede'$^2$ }
\date{%
	$^1$\small{Oliver Wyman Srl, Via Broletto 16, 20121, Milano, Italy}\\[0.3cm]
	$^2$\small{MOX, Dipartimento di Matematica, Politecnico di Milano,\\ Piazza Leonardo da Vinci 32, 20133, Milano, Italy}\\[0.3cm]%
	$^*$\small{Corresponding author: {\tt alberto.zingaro@polimi.it}}\\[0.3cm]
	\today
}

\maketitle

\begin{abstract}
  We propose using machine learning and artificial neural networks (ANNs) to enhance residual-based stabilization methods for advection-dominated differential problems. Specifically, in the context of the finite element method, we consider the streamline upwind Petrov-Galerkin (SUPG) stabilization method and we employ ANNs to optimally choose the stabilization parameter on which the method relies. We generate our dataset by solving optimization problems to find the optimal stabilization parameters that minimize the distances among the numerical and the exact solutions for different data of differential problem and the numerical settings of the finite element method, e.g., mesh size and polynomial degree. The dataset generated is used to train the ANN, and we used the latter ``online'' to predict the optimal stabilization parameter to be used in the SUPG method for any given numerical setting and problem data. We show, by means of 1D and 2D numerical tests for the advection-dominated differential problem, that our ANN approach yields more accurate solution than using the conventional stabilization parameter for the SUPG method.
\end{abstract}

\textbf{Keywords: }{partial differential equations; finite element method; stabilization methods; streamline upwind Petrov-Galerkin; machine learning; artificial neural networks 

\section{Introduction}

The Galerkin-finite element (FE) method applied to partial differential equations (PDEs) with advection terms dominating over diffusion ones may suffer of numerical instability \cite{hughes2012finite,Quarteroni_2017}. These numerical instabilities cause the numerical solution to exhibit oscillations that increase in amplitude with a local increment of the transport dominance over diffusion, i.e., as soon as the local Péclet number becomes larger than one. In order to eliminate (or at least mitigate) numerical instabilities, a general employed strategy is to consider the generalized Galerkin method, i.e., the Galerkin method with additional stabilization terms \cite{Quarteroni_2017}. Examples of stabilization methods to reduce numerical oscillations in the advection dominated regimes are for instance the upwind and streamline-diffusion method, both methods that are not strongly consistent. Instead, examples of strongly consistent methods are the Galerkin-least-squares, the streamline upwind Petrov-Galerkin (SUPG), and the Douglas-Wang methods \cite{Roos_1996,Brooks_1982,Quarteroni_2017}. In particular, for these strongly consistent stabilization methods, the formulation depends on the residual of the PDE (in strong formulation), other than on a parameter -- called \textit{stabilization parameter} -- whose definition is crucial for the success of the stabilization strategies. Determining the values of such stabilization parameter is not straightforward, especially for 2D and 3D problems. In addition, a universal formulation of such parameter is lacking, especially as it is strongly dependent on data of the model and the numerical setting used for the FE approximation of the PDEs, e.g., low and high order FE methods, spectral element methods, isogeometric methods,~etc.~\cite{Schwab, spectralbook, cottrell2009isogeometric}. Some of the formulations for the stabilization parameter have been derived from analytical considerations on the differential problem in 1D, others are suitable only for a specific numerical approximation method, while the most comprehensive ones incorporate some (empirical) dependence on numerical setting, as e.g., the order of the FE \cite{Bochev,Brooks_1982, CODINA200061,CODINA2008264,CODINA20072413,franca1992stabilized,galeao2004finite,REBOLLO2015406,ScovazziRossi,SCOVAZZI2007923,tezduyar2000finite}.

In this work, we propose using machine learning (ML) \cite{mitchell1997machine,Yegnanarayana,goodfellow2016deep} to make computers learn autonomously the values of the stabilization parameter for the FE approximation of advection dominated PDEs. Artificial neural networks (ANNs) \cite{university1988continuous,Kutyniok_2021} are widely popular in ML and Deep Learning for a wide array of applications: they are indeed very versatile tools that are increasingly finding their way in scientific computing \cite{Neittaanmaki202227, Mishra}, especially in the context of numerical approximation of PDEs. For instance, as substitute to standard numerical methods, ANNs can be employed as meshless methods in physics informed neural networks (PINNs) to directly approximate the solution of the PDE as it is trained by minimizing the (strong) residual of the PDE \cite{raissi2017machine,raissi2018hidden,raissi2019physics}. ANNs are also largely employed in a data-driven fashion in the context of model order reduction for parametric PDEs \cite{fresca2021comprehensive, regazzoni2019machine, hesthaven2018non, guo2019data, zancanaro2021hybrid} and to enhance the stability properties of numerical methods for PDEs \cite{discacciati2020controlling}. In fluid dynamics modelling, ANNs are massively adopted for flow features extractions, modelling, optimization and control. For flow features extractions, ANNs are used through clustering and classification to classify wake topologies \cite{annurev_wake}; for modeling fluid dynamics by reconstructing specific flows such as the near wall field in a turbulent flow \cite{annurev_turbulent_wall} and for flow optimization \cite{annurev_gliding}, and control  for aerodynamics applications \cite{annurev_control}. ANNs are also used in a data-driven framework as a manner for providing alternative closure models for stress tensor \cite{annurev_rans} in Reynolds-Average Navier-Stokes (RANS) equations, for sub-grid scale models in large eddy simulation (LES) turbulence models \cite{janssensadvancing, janssenthesis, xie2019artificial, zhou2019subgrid, xie2020modeling}, or for model learning input-output relationships in complex physical processes \cite{REGAZZONI2020113268}. 

In this work, we use ANNs to learn the optimal stabilization parameter in advection dominated PDEs that are discretized by means of the FE method: the goal is to enhance the accuracy of the SUPG FE method by optimally selecting the stabilization parameter under different data of the PDE and numerical settings of the FE method. We found that the proposed ANN-enhanced stabilization method allows to improve accuracy and stabilization properties of the numerical solution compared to those results obtained by analytical expressions of the stabilization parameter. The numerical results obtained shed light also on the possibility to apply the presented strategy to learn closure laws for stabilization and turbulence models of fluid dynamics, as for instance to learn the stabilization parameters in the Variational Multiscale--LES model to model transitional and turbulent flows \cite{forti2015semi, bazilevs2007variational, zingaro2021hemodynamics}.

This work is organized as follows: in Section~\ref{section:advection-diffusion}, we recall the SUPG stabilization method for advection-diffusion equations; in Section~\ref{sec:strategy}, we present our numerical strategy to compute an optimal SUPG stabilization parameter through a feed-forward fully connected ANN. In Section~\ref{sec:validation} we validate our method by comparing the ANN results with those obtained with the 1D advection-diffusion problem from which the expression of the theoretical stabilization parameter has been derived. In Section~\ref{sec:results} we first show the ANN's training and we present our numerical results on the 2D advection-diffusion problem used for training and we finally generalize our findings using the ANN's prediction on a different advection-diffusion problem. Finally, in Section~\ref{sec:conclusions} we draw our conclusions highlighting possible future developments.

\section{The SUPG method for advection-diffusion problems}
\label{section:advection-diffusion}
We briefly recall the advection-diffusion differential problem and the SUPG stabilization method for the advection dominated regime.

Let $\Omega \in \mathbb{R}^d, \; d = 1, 2, 3$ be the physical domain with $\partial \Omega$ being its boundary. We consider the following problem in the unknown function $u$:
\begin{equation}
    \label{eq:transport_diffusion}
    \left\{ \; \begin{aligned}
        - \nabla \cdot (\mu \nabla u) + \boldsymbol{\beta} \cdot \nabla u &= f & \qquad \text{in} \; \Omega, \\
        u &= g & \qquad \text{on} \; \partial \Omega,
    \end{aligned} \right.
\end{equation}
where $\mu$, $\boldsymbol{\beta}$, and $f$ are assigned functions or constants, with $\mu \in L^\infty(\Omega)$, $\bm \beta \in [L^\infty(\Omega)]^d$, with $\nabla \cdot \bm \beta \in L^2(\Omega)$, and $f \in L^2(\Omega)$. 
The Dirichlet datum on the boundary is $g\in H^{1/2}(\partial \Omega)$.
Let $V_g = \{ v \in H^1(\Omega): \, v|_{\partial \Omega} = g \}$ and $V_0 = H_0^1(\Omega)$; the weak formulation of Eq~\eqref{eq:transport_diffusion} reads
\begin{equation}
    \label{eq:weak_bilinear}
    \text{find} \; u \in V_g \ : \ a(u,v) = F(v) \quad \text{for all } v \in V_0.
\end{equation}
with the bilinear form $a(u,v)$ and the linear functional $F(v)$  respectively defined as:
\begin{equation*}
    \label{eq:bilinearForm}
    \begin{aligned}
        a(u,v) &:= \int_\Omega \mu \, \nabla u \cdot \nabla v \, d\Omega + \int_\Omega v \, \boldsymbol{\beta} \cdot \nabla u \, d\Omega, \\
        F(v) &:= \int_\Omega f \, v \, d\Omega.
    \end{aligned}
\end{equation*}
We consider a family of function spaces $V_h \subset V$ (either for $V_g$ and $V_0$) dependent on a parameter $h$ such that $\dim(V_h) = N_h < \infty$.  Let $X_r^h = \{ v^h \in C^0(\overline{ \Omega}):\,  v^h|_K \in \mathbb{P}_r,\,  \, \text{for all } \, K \in \mathcal T_h\}$ be the function space of the FE discretization  with piecewise Lagrange polynomials of degree $r \geq 1$, $\mathcal T_h$ a triangulation of $\Omega$ and $h$ the characteristic size of the mesh, comprised of elements $K \in \mathcal T_h$. By setting $V_h = {X}_h^r \bigcap V_0$, the Galerkin FE method applied to Eq~\eqref{eq:weak_bilinear} reads
\begin{equation}
    \label{eq:galerkin_1}
    \text{find} \; u_h \in V_{g,h} \ : \ a(u_h,v_h) = F(v_h) \quad \text{for all } v_h \in V_{h},
\end{equation}
where $V_{g,h}= {X}_h^r \bigcap V_g$.
The standard Galerkin-FE method, in Eq~\eqref{eq:galerkin_1}, can generate numerical oscillations on $u_h$ if the problem is dominated by the advection term. In particular, these numerical instabilities can arise if the local P\'eclet number is $\Peclet_h > 1$, where
\begin{equation}
    \label{eq:peclet}
    \Peclet_h = \frac{|\boldsymbol{\beta}| h}{2 \mu}.
\end{equation}
The {generalized Galerkin} method (\cite{QV94}) allows for eliminating or mitigating numerical oscillations by adding  stabilization terms to the standard Galerkin formulation as
\begin{equation}
    \label{eq:generalized_galerkin}
    \text{find} \; u_h \in V_{g,h} \ : \ a(u_h,v_h) + b_h(u_h,v_h) = F(v_h) \quad \text{for all } v_h \in V_{h}.
\end{equation}
In the SUPG method, the additional stabilization term reads:
\begin{equation}
    \label{eq:supg}
    b_h(u_h,v_h) = 
    \sum_{K \in \mathcal T _h}
    \int_K R(u_h) \; \tau_K \; \frac{1}{2} \left( \nabla \cdot (\boldsymbol{\beta} v_h) + \boldsymbol{\beta} \cdot \nabla v_h \right) \, d\Omega,  
\end{equation}
where $R(u_h)$ is the residual in strong formulation of Eq~\eqref{eq:transport_diffusion}, which is defined as:
\begin{equation*}
    R(u_h) := - \nabla \cdot (\mu \nabla u_h) + \boldsymbol{\beta} \cdot \nabla u_h - f.
\end{equation*}
The term $\tau_K$, appearing in Eq~\eqref{eq:supg}, is the \textit{stabilization parameter}, which is the focus of this work. The stabilization parameter ${\tau}_K$ is generally defined locally, i.e., mesh element by element. In this paper, we consider a uniform stabilization parameter, thus $\tau_K = \tau$ for all $K \in \mathcal T_h$.

A universal and optimal definition of $\tau$ in terms of the problem data and numerical settings like the mesh size and FE degree is lacking. An extensive review of stabilization parameters for the SUPG method is reported in \cite{john2007spurious}. The most common choices for $\tau$ come from an analytic derivation made for the advection-diffusion problem in 1D with $f=0$ and approximated by means of linear finite elements, which reads:
\begin{equation}
    \label{eq:tau_theory}
    \widetilde{\tau}_1 = \frac{h}{2 |\boldsymbol{\beta}|} \, \xi(\Peclet_h),
\end{equation}
where $\xi(\theta)$ is the upwind function:
\begin{equation*}
    \xi(\theta) = \coth({\theta}) - \frac{1}{\theta},   \quad \theta>0.
\end{equation*}
If a uniform mesh is used, as we consider in this paper, then the value of $h$ is uniform over $\mathcal T_h$; if in addition $\boldsymbol{\beta}$ and $\mu$ are constant, then this implies that $\tau$ is uniform over $\mathcal T_h$ . The choice of $\tau$ made in Eq~\eqref{eq:tau_theory} represents an optimal choice of the stabilization parameter as it yields a nodally exact numerical solution for the 1D advection-diffusion problem if $f = 0$ and the FE polynomial order is $r = 1$ \cite{Quarteroni_2017, galeao2004finite}. 
Thus, the stabilization parameter of Eq~\eqref{eq:tau_theory} may not be fully effective to provide the optimal stabilization for a general advection-diffusion problem, e.g., to guarantee a nodally exact solution in 2D/3D or when using FE degrees larger than $r=1$. 
A commonly used generalization of the formula of Eq~\eqref{eq:tau_theory} to higher FE degrees ($r>1$) is presented in \cite{galeao2004finite} and reads:
\begin{equation}
    \label{eq:tau_theory_r}
    \widetilde{\tau}_r = \frac{h}{2 |\boldsymbol{\beta}| r} \, \xi \left( \frac{\Peclet_h}{r} \right).
\end{equation}
Differently from Eq \eqref{eq:tau_theory}, the stabilization parameter $\widetilde{\tau}_r$ in Eq~\eqref{eq:tau_theory_r} takes into account the contribution of higher polynomial degrees $r$. Still, this formula is not optimal as it does not guarantee a nodally exact numerical solution. 
Our goal is to find a general and optimal expression of the stabilization parameter holding for advection-diffusion problems and FE approximations of degree $r\geq 1$ to be dependent on: the dimension $d$, the FE degree $r$, the mesh size $h$, the forcing term $f$, the diffusion coefficient $\mu$ and the transport coefficient $\beta$.

\section{ANN-based approach for determining the optimal stabilization parameter}
\label{sec:strategy}

We present our approach to determine the optimal stabilization parameter $\tau$ for the SUPG stabilization method by using an ANN.  We consider as {features} (inputs) of our ANN: the FE degree $r$, the mesh size $h$ and the global P\'eclet number $\mathbb Pe_g := | \boldsymbol{\beta}| L / (2\mu)$ of the advection-diffusion problem, where $L$ is the characteristic length of the problem,  $\boldsymbol \beta \in \mathbb R^d$ and $\mu \in \mathbb R$ . As we consider $\boldsymbol{\beta}$ to be uniform and fixed, using $\Peclet_g$ as input corresponds to varying the value of the diffusion coefficient $\mu$. The features (input) of the ANN read:
\begin{equation}
\mathbf x^{(i)} = \left [ r, \, h,\,  \Peclet_g \right ]; 
\label{input}
\end{equation}
the \textit{target} (output) of interest is the optimal SUPG stabilization parameter that we denote with $\tau^*$:
\begin{equation}
\mathbf y^{(i)} = \left [ \tau^* \right ].
\label{output}
\end{equation}
The first step consist in generating the dataset, i.e., pairs of inputs and outputs to be used for training the ANN (\textit{data generation} step). This consists in choosing repeatedly and randomly values of the features $\mathbf x^{(i)}$ in given ranges. These features are used to feed an \textit{optimization problem} that, through an optimizer and a suitable error measure, provide an optimal stabilization parameter $\mathbf y^{(i)}$, i.e., the target of the given feature. Such optimization problem considers as error measure $E(\tau)$: the mismatch between the numerical solution $u_h$ and the exact one $u_\mathrm{ex}$ over the nodes of the FE mesh. $E(\tau)$ reads as:
\begin{equation}
    \label{eq:td_error_definition}
    E(\tau) = \sum_{k=1}^{K_h} | \, e(\bm x_k, \tau) \, |, \qquad e(\bm x_k, \tau) = u_h(\bm x_k;\tau) - u(\bm x_k).
\end{equation}
being $\bm x_k$ the $k$-th node of the FE mesh $\mathcal T_h$ and $K_h$ the number of nodes. $E(\tau)$ is an approximation of the $L_1(\Omega)$ norm.  To find the minimum of $E(\tau)$ for different problem configurations we solve, for each instance of the input parameters $\mathbf x^{(i)}$, the following optimization problem:
\begin{equation}
    \label{eq:td_findtau}
    \text{find} \; \tau^* \ : \quad \min_\tau \, E(\tau).
\end{equation}
Thus, the dataset generated consists of $m$ pairs of inputs-outputs $(\mathbf x^{(i)}, \mathbf y^{(i)})$, with $i=1, \dots, m$. The latter is used for \textit{training} the ANN. The latter will be then used to predict the optimal stabilization parameter $\mathbf y^{(j)}= \tau_\mathrm{ANN}$ to be used for the FE approximation of the advection-diffusion problem in the settings provided by $\mathbf x^{(j)}$.

We report in Figure~\ref{fig:td_scheme} a sketch of the procedure used to build the ANN for the prediction of the optimal stabilization parameter for any new feature $\mathbf x^{(j)}$.

\begin{figure}[H]
	\centering
	\includegraphics[width=0.8\textwidth]{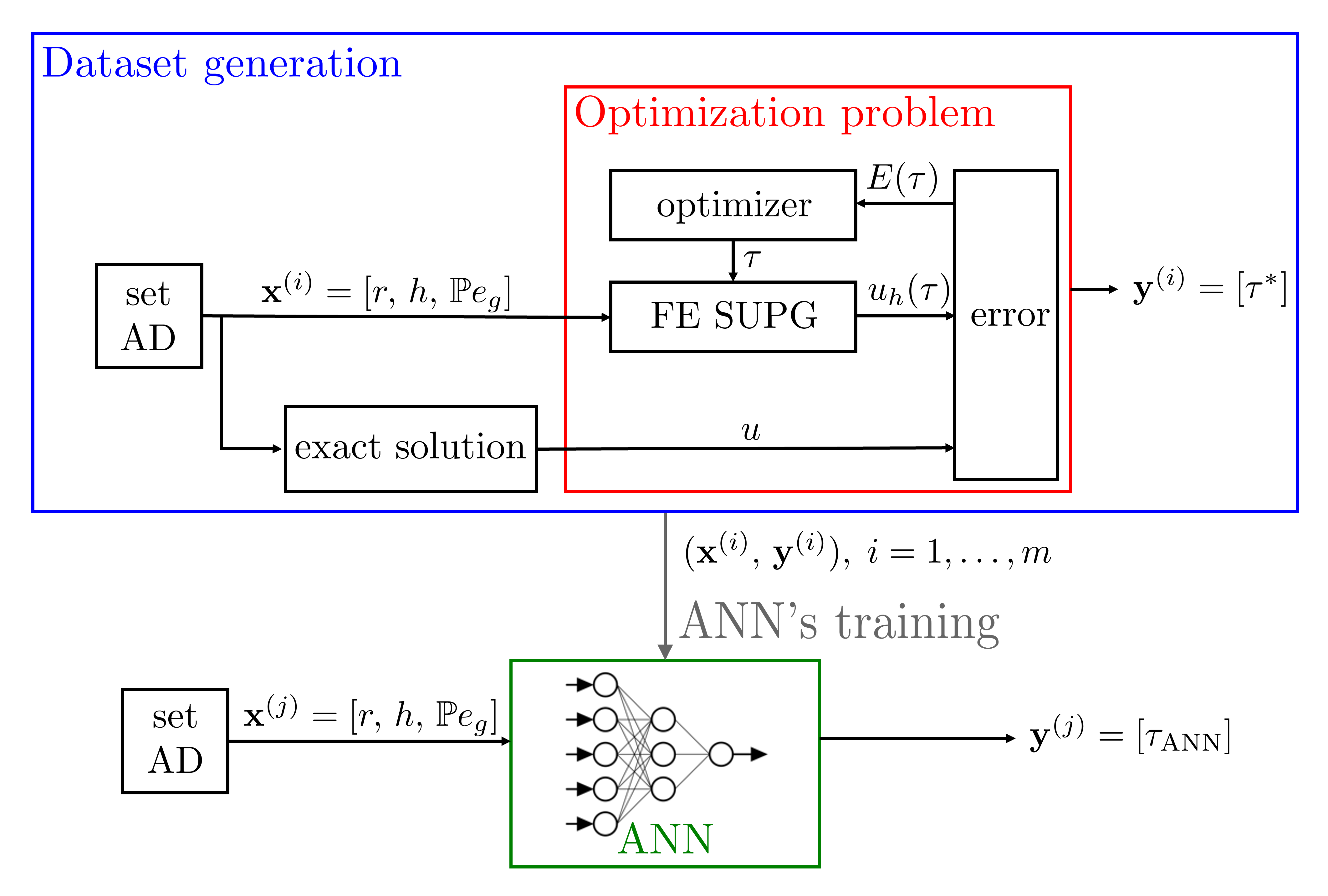}
	\caption{Representation of the strategy used for learning the the relation among $\tau$ and the advection-diffusion parameters and numerical settings.}
	\label{fig:td_scheme}
\end{figure}

\section{Numerical validation}
\label{sec:validation}

In this section, we validate the proposed numerical strategy by means of a 1D advection-diffusion problem from which the expression of $\widetilde{\tau}_1$ in Eq~\eqref{eq:tau_theory} has been derived. In particular, we consider the advection-diffusion problem in Eq~\eqref{eq:transport_diffusion} with $\Omega = (0, 1)$, $f = 0$ and with Dirichlet BCs $ u = 0$ on $x = 0$ and $ u = 1$ on $x = 1$. It admits the following exact solution:
\begin{equation*}
    u(x) = \frac{\exp{\left( \frac{\beta}{\mu} x \right)} - 1}{\exp{\left( \frac{\beta}{\mu} \right)} - 1}.
\end{equation*}

We recall that the stabilization parameter $\widetilde{\tau}_1$ of Eq~\eqref{eq:tau_theory} yields a nodally exact numerical solution if Eq~\eqref{eq:transport_diffusion} is solved by linear FE ($r = 1$).

We plot in Figure~\ref{fig:td_minimumErr} the error measure $E(\tau)$ against the value of the stabilization parameter $\tau$, with $h = 1/20$, $\Peclet_h = 12.5$ and by using $r = 1$ (Figure \ref{fig:td_minimumErr_r1}) and $r = 3$ (Figure \ref{fig:td_minimumErr_r3}). We observe that $E(\tau)$ shows a minimum in $\tau^*$, thus suggesting the possibility to use an optimization algorithm to solve the problem. Specifically, we employ the L-BFGS-B optimization algorithm from {\tt SciPy}, an open source library for {\tt Python} \cite{scipy}.  Instead, the advection-diffusion problem has been solved using the FE open source library {\tt FEniCS} \cite{fenics}, by applying the SUPG method on a uniform mesh.
Particularly, the optimal value $\tau^*$ found for the case $r=1$ corresponds to the one provided by the theory in Eq~\eqref{eq:tau_theory}. For the case $r=3$, we observe a optimal value of $\tau^*$ for which the error is minimized: in this case, it does not corresponds to the stabilization parameter $\widetilde{\tau}_r$ provided by the theory in Eq~\eqref{eq:tau_theory_r}: as a matter of fact, the latter arises from empirical considerations to extend the parameter $\widetilde{\tau}_1$ in Eq~\eqref{eq:tau_theory} to polynomials degrees $r>1$. However, this does not ensure a nodally exact numerical solution. 

\begin{figure}[H]
	\centering
	\begin{subfigure}{0.4\textwidth}
		\includegraphics[width=\textwidth]{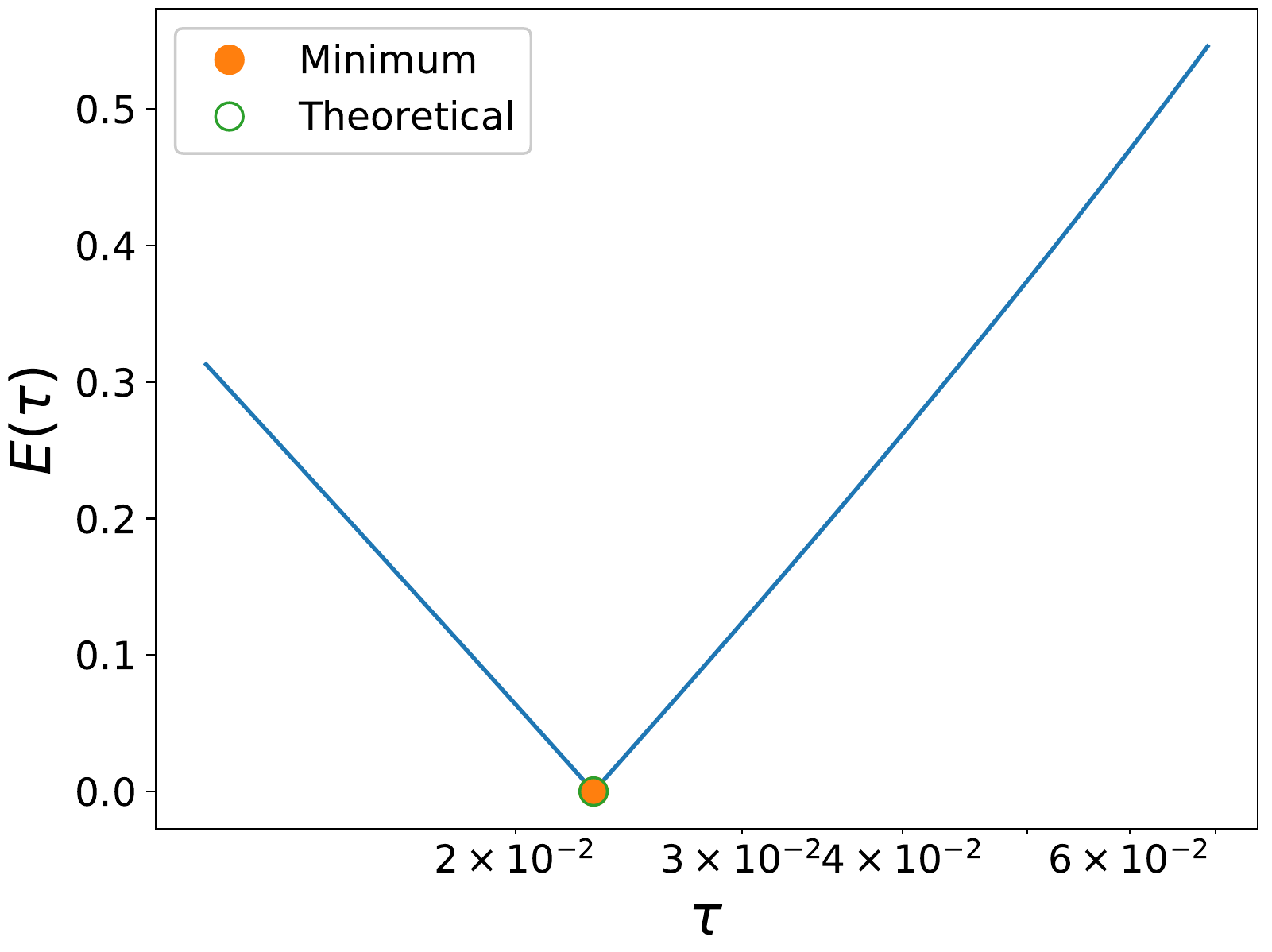}
		\caption{$r = 1$}
		\label{fig:td_minimumErr_r1}
	\end{subfigure}
	\begin{subfigure}{0.4\textwidth}
		\includegraphics[width=\textwidth]{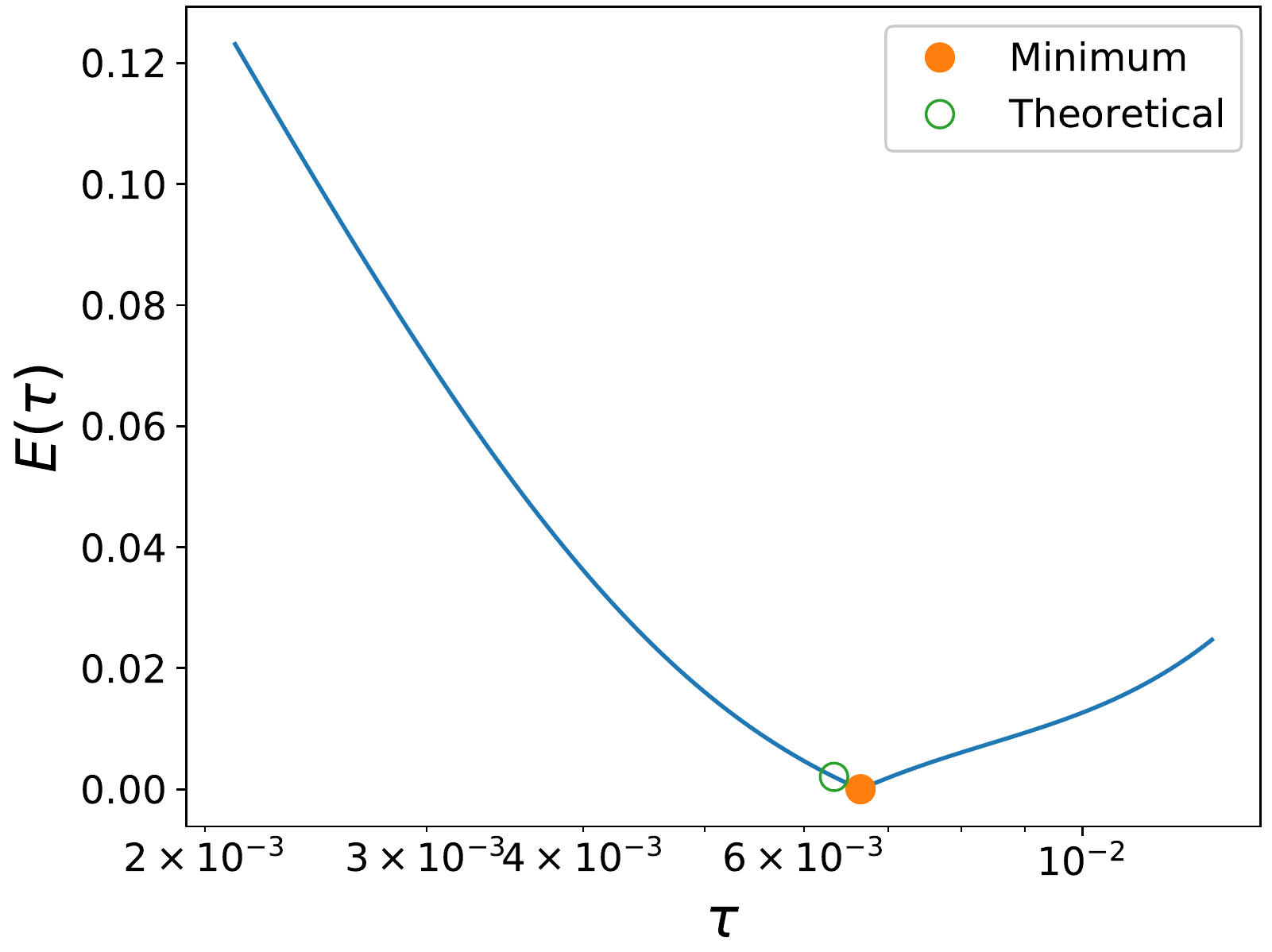}
		\caption{$r = 3$}
		\label{fig:td_minimumErr_r3}
	\end{subfigure}
	\caption{Error (\ref{eq:td_error_definition}) of numerical solutions with SUPG method vs. the stabilization parameter $\tau$ for the 1D problem in Section~\ref{sec:validation} with $h = 1/20$ and $\Peclet_h = 12.5$; comparison between the theoretical value $\widetilde{\tau}_r$ of Eq~\eqref{eq:tau_theory_r} and the optimal one $\tau^*$.}
	\label{fig:td_minimumErr}
\end{figure}

\begin{figure}[H]
	\centering
	\begin{subfigure}{0.42\textwidth}
		\includegraphics[width=\textwidth]{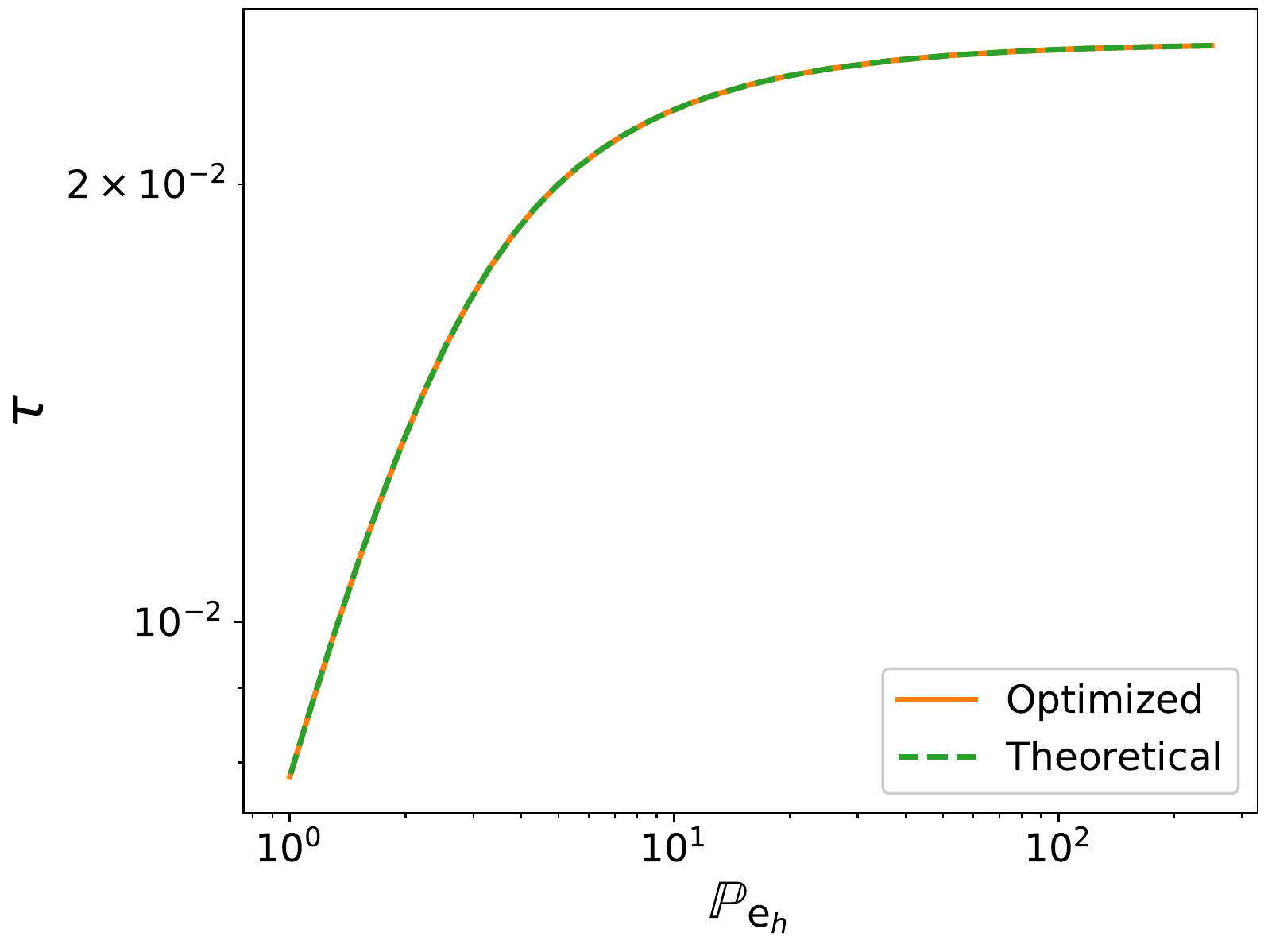}
		\caption{$r = 1$}
		\label{fig:td_1d_comparison_r1}
	\end{subfigure}
	\begin{subfigure}{0.42\textwidth}
		\includegraphics[width=\textwidth]{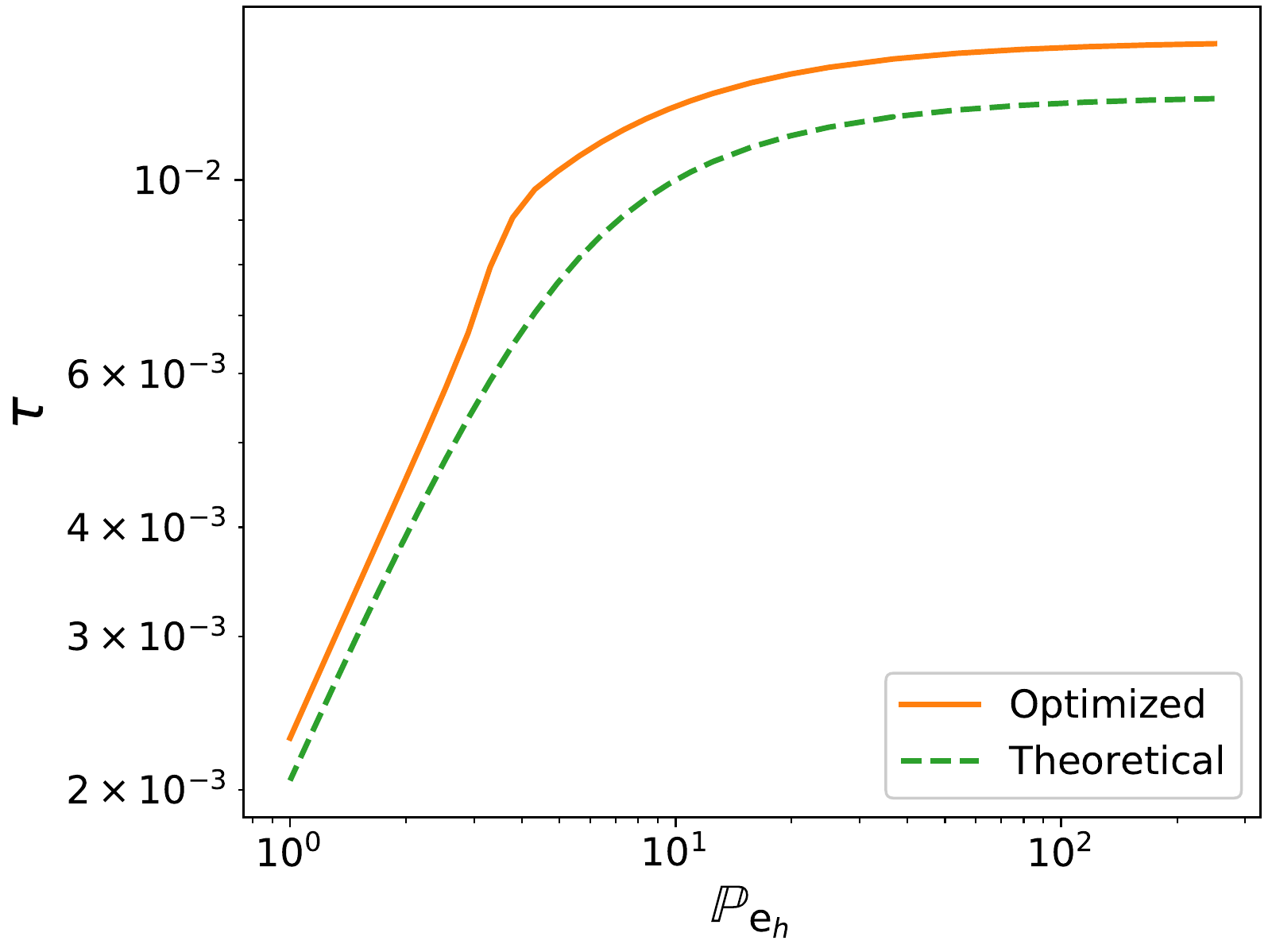}
		\caption{$r = 2$}
		\label{fig:td_1d_comparison_r2}
	\end{subfigure}
	\\
	\begin{subfigure}{0.42\textwidth}
		\includegraphics[width=\textwidth]{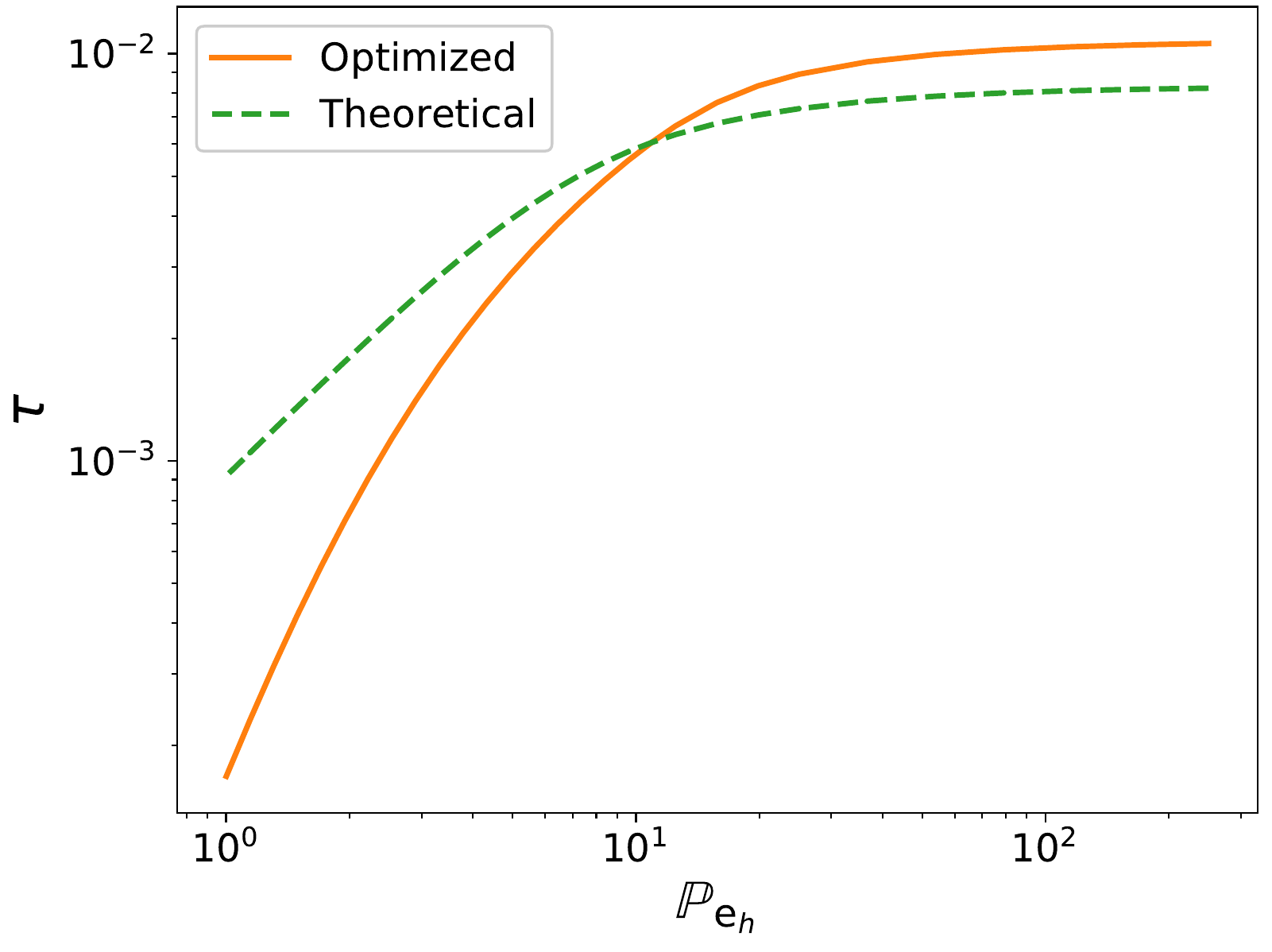}
		\caption{$r = 3$}
		\label{fig:td_1d_comparison_r3}
	\end{subfigure}
	\caption{Comparison of the optimal stabilization parameter $\tau^*$ (full line) and the theoretical one $\widetilde{\tau}_r$ of Eq~(\ref{eq:tau_theory_r}) (dashed line) against $\Peclet_h$, with $h = 1/20$. Results are referred to the 1D advection-diffusion problem described Section~\ref{sec:validation}.}
	\label{fig:td_1d_comparison}
\end{figure}

We solve the optimization problem for different values of the local P\'eclet number $1 \leq \Peclet_h \leq 250$ and we plot in Figure~\ref{fig:td_1d_comparison} the optimal stabilization parameter $\tau^*$ against the local P\'eclet number for $r = 1, 2$, and $3$.

Figure~\ref{fig:td_1d_comparison_r1} shows that, by employing a linear FE space $r=1$, the optimal stabilization parameter $\tau^*$ obtained by minimizing $E(\tau)$ exactly matches the theoretical one $\widetilde{\tau}_1$ for every value of $\Peclet_h$, thus confirming the findings of Figure~\ref{fig:td_minimumErr_r1}. This also serves as validation of the optimization procedure that we proposed. By recalling that $\widetilde{\tau}_1$ of Eq~\eqref{eq:tau_theory} guarantees the numerical solution to be exact at nodes for this particular advection-diffusion problem and $r=1$, we can infer that the optimization procedure is meaningful and it can be therefore exploited further, for example for $r>1$. On the other hand, Figure~\ref{fig:td_1d_comparison_r2} and~\ref{fig:td_1d_comparison_r3} show instead a mismatch between the theoretical $\widetilde{\tau}_r$ of Eq~\eqref{eq:tau_theory_r} and the optimal one $\tau^*$, thus confirming the findings reported in Figure~\ref{fig:td_minimumErr_r3}.

In order to better appreciate the differences in the numerical solutions due to the choice of the stabilization parameters, we report, in Figure~\ref{fig:td_1d_error}, a comparison of the error $E(\tau)$ obtained by means of the optimal $\tau^*$ and theoretical $\widetilde{\tau}_r$ stabilization parameters for $r=2$ and $3$. Moreover, we report in Figure~\ref{fig:td_boundary_layer_r3} a comparison among the exact solution $u$, the SUPG-stabilized numerical solution $u^*_h$ obtained with the optimal stabilization parameter $\tau^*$, and the numerical solution $\widetilde{u}_h$ obtained with theoretical one $\widetilde{\tau}_r$. Specifically, we consider FE of degree $r=3$ and two different values of $\Peclet_h$. In both these cases, the optimal parameter $\tau^*$ leads to a more accurate solution with respect to using the theoretical parameter $\widetilde{\tau}_r$. In particular, when $\Peclet_h$ is ``small", $\widetilde{\tau}_r$ leads to overshooting in the numerical solution, while if $\Peclet_h$ is ``large" the theoretical stabilization parameter leads to undershooting. Conversely, the optimal parameter $\tau^*$ accurately intercepts the exact solution at the nodes.

\begin{figure}[H]
	\centering
	\begin{subfigure}{0.49\textwidth}
		\includegraphics[width=\textwidth]{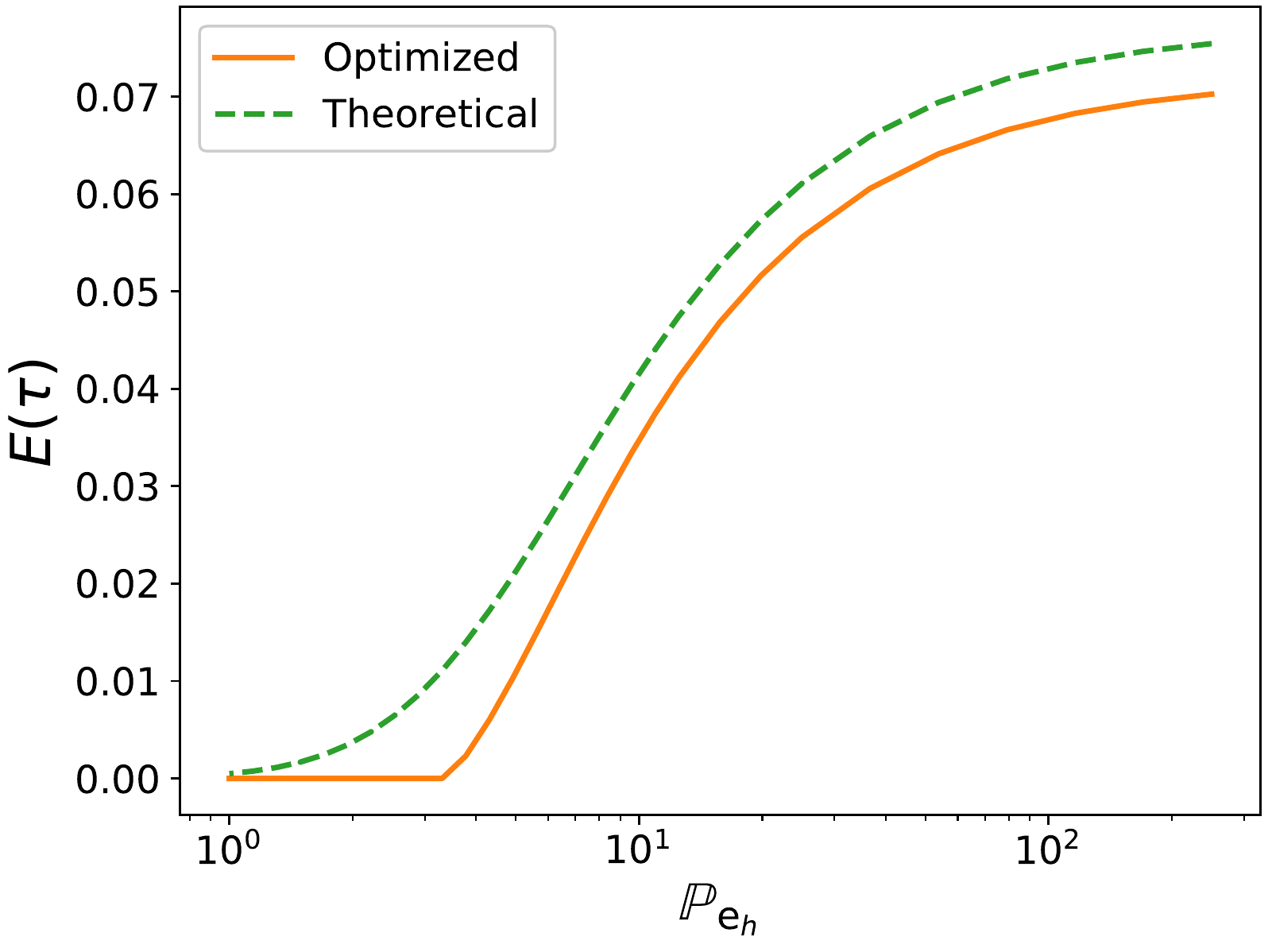}
		\caption{$r = 2$}
	\end{subfigure}
	\begin{subfigure}{0.49\textwidth}
		\includegraphics[width=\textwidth]{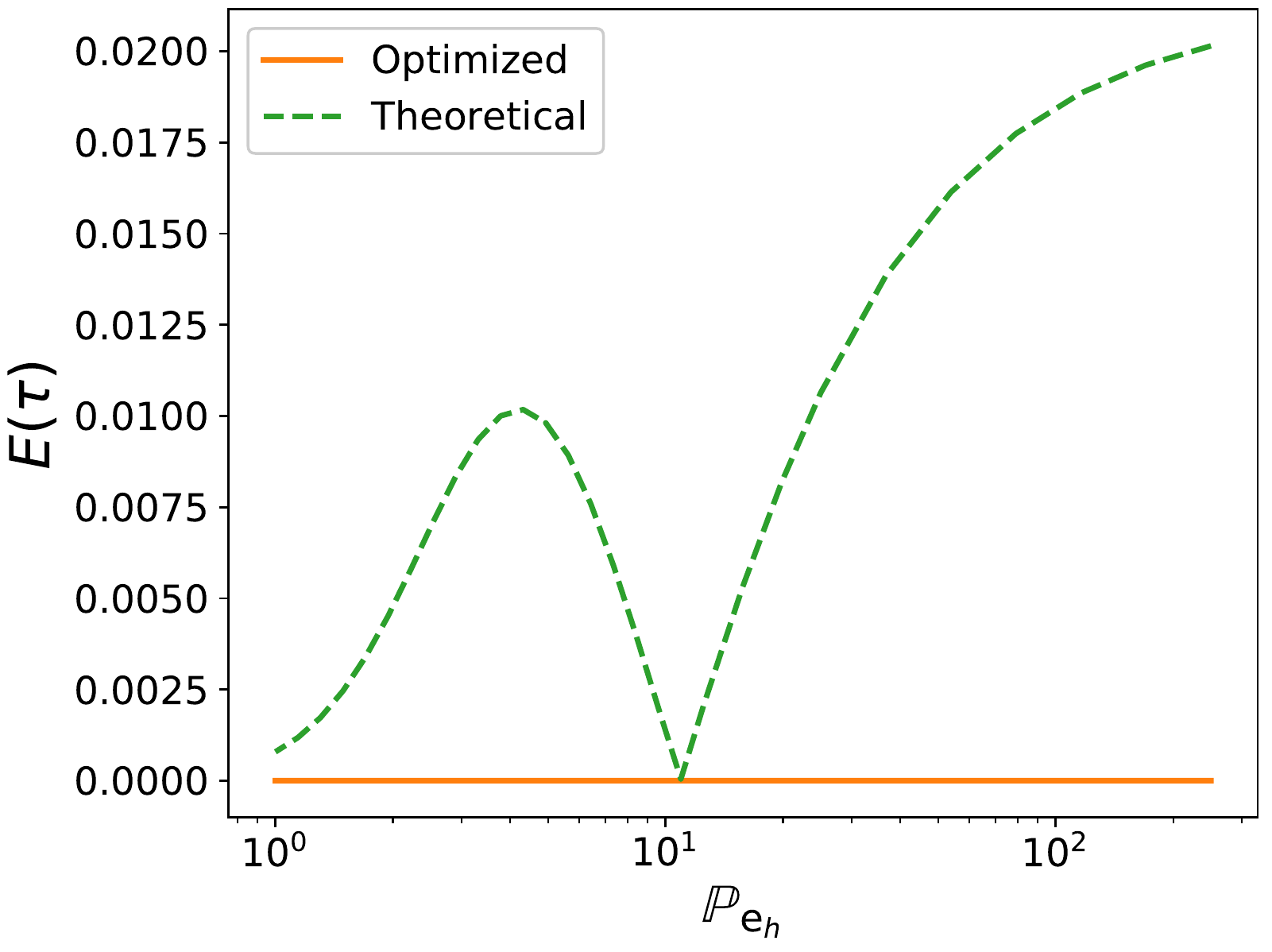}
		\caption{$r = 3$}
	\end{subfigure}
	\caption{Comparison of the error of Eq (\ref{eq:td_error_definition}) obtained with the optimal $\tau^*$ (full line) and theoretical $\widetilde{\tau}_r$ (dashed line) for varying $\Peclet_h$ with $h = 1/20$. Results referred to the 1D advection-diffusion problem of Section~\ref{sec:validation}.}
	\label{fig:td_1d_error}
\end{figure}
\newpage
\begin{figure}[H]
    \centering
    \begin{subfigure}{0.49\textwidth}
         \includegraphics[width=\textwidth]{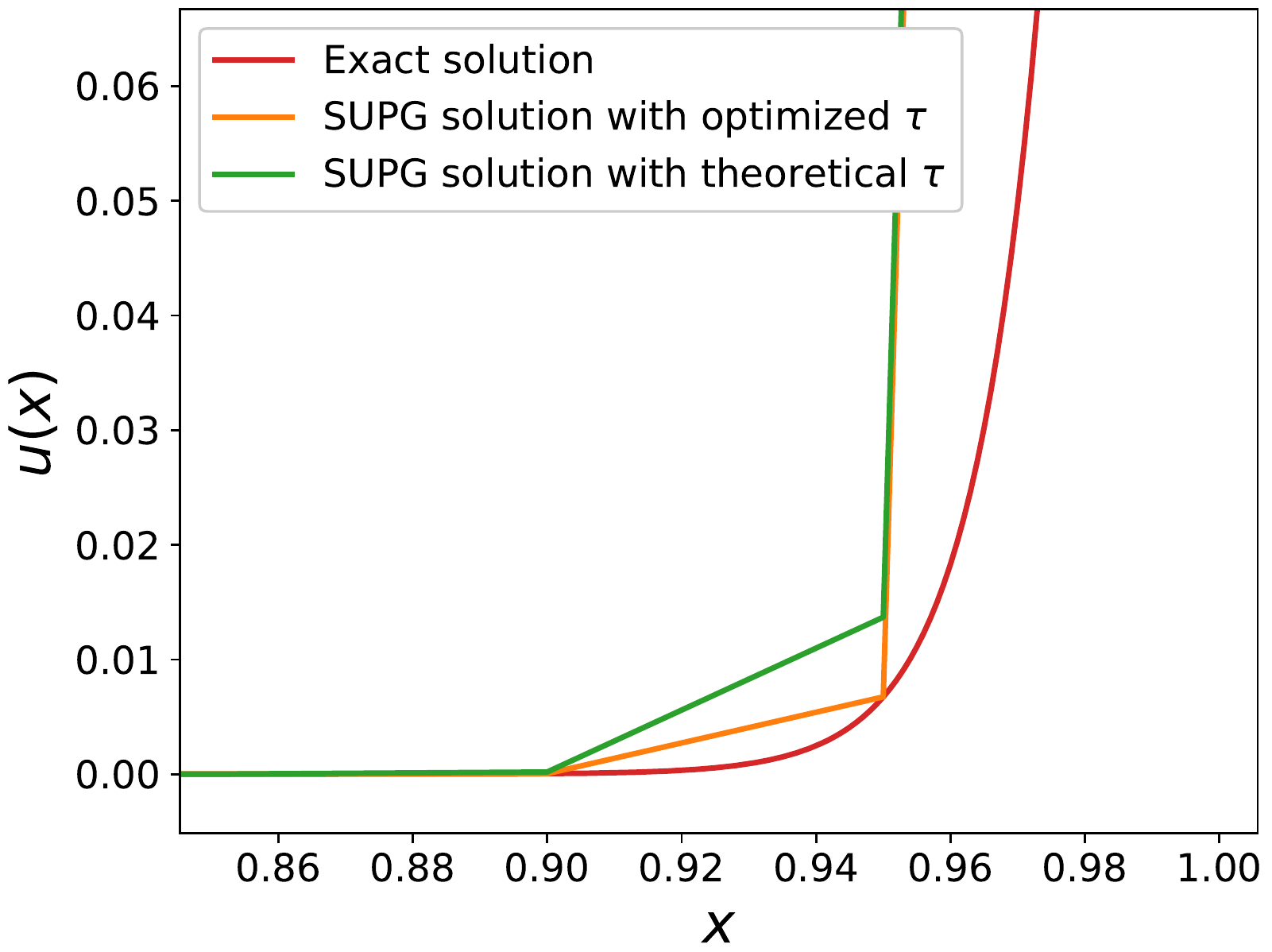}
         \caption{$\Peclet_h = 2.5$}
    \end{subfigure}
    \begin{subfigure}{0.49\textwidth}
         \includegraphics[width=\textwidth]{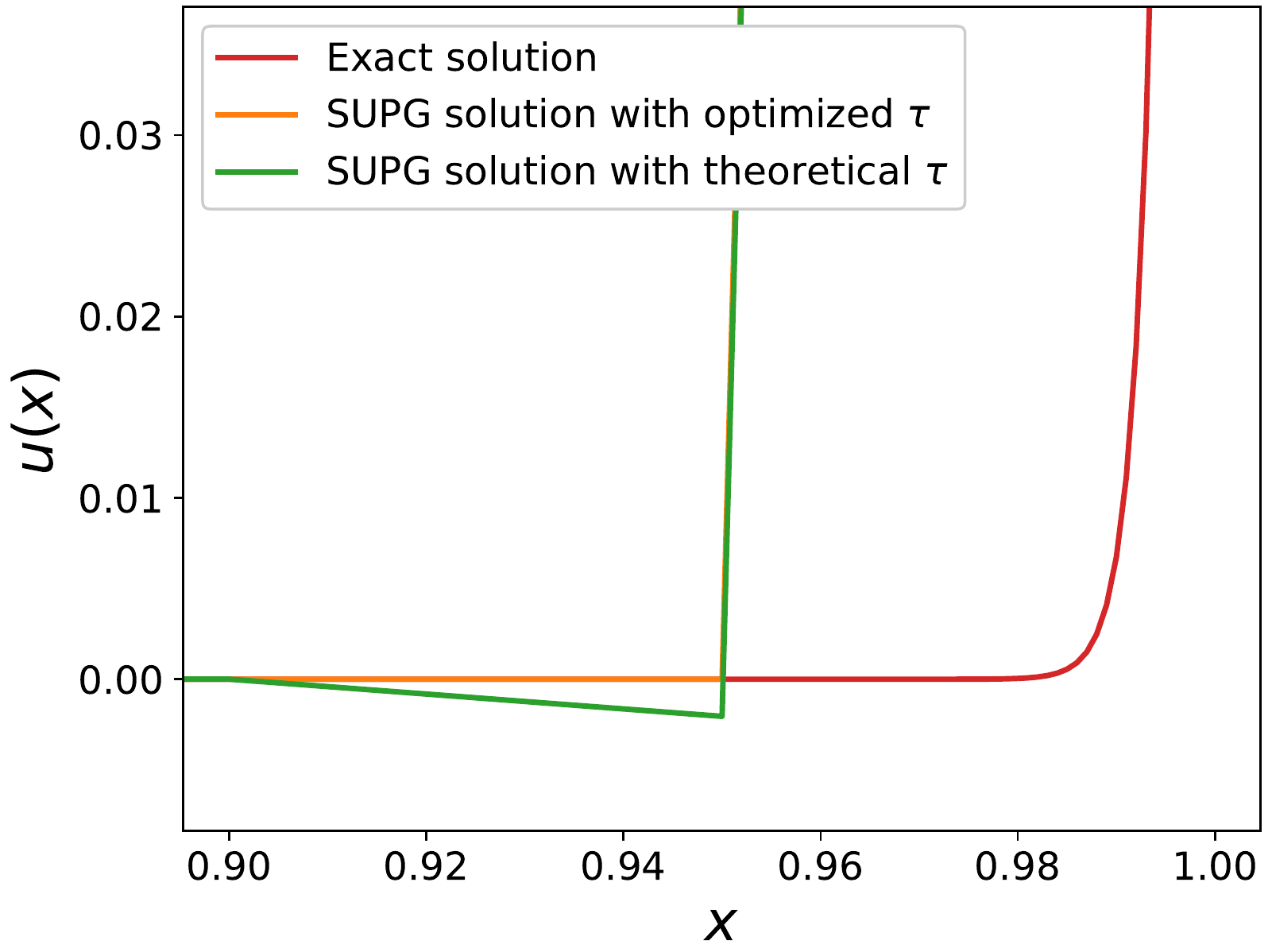}
         \caption{$\Peclet_h = 12.5$}
    \end{subfigure}
    \caption{Boundary layers of the numerical solutions $u^*_h$ and $\widetilde{u}_h$ obtained with the SUPG method with optimal $\tau^*$ and theoretical one $\widetilde{\tau}_r$ for $h = 1/20$ and $r = 3$, respectively; comparison with the exact solution $u$ ($u^*_h$ is nodally exact at the node in $x=0.95$). Results are referred to the 1D advection-diffusion problem of Section~\ref{sec:validation}.}
    \label{fig:td_boundary_layer_r3}
\end{figure}


\section{Numerical results}
\label{sec:results}
We first introduce the training set of our problem and we detail the setup of the ANN that we use in this work. Then, we show the prediction of the stabilization parameter by means of the ANN on different advection-diffusion problems.

\subsection{Training the ANN for a 2D advection-diffusion problem}
\label{sec:training_2D}
We apply the strategy presented so far to a 2D advection-diffusion problem to generate the dataset for the training of the ANN. Specifically, we consider in Eq~\eqref{eq:transport_diffusion}: $\Omega = (0, 1)^2$, $f = 0$, and $\boldsymbol{\beta} = (1,1)$. We prescribe the following exact solution on the whole boundary $\partial \Omega$:
\begin{equation}
    \label{eq:td_2d_exact}
    u(x, y) = \frac{e^{(x / \mu)} - 1}{e^{(1 / \mu)} - 1} + \frac{e^{(y / \mu)} - 1}{e^{(1 / \mu)} - 1}.
\end{equation}

We generate in $\Omega$ a structured mesh $\mathcal T_h$ of triangles with {\tt FEniCS} \cite{fenics}, as shown in Figure~\ref{fig:td_2d_mesh}.
We generate the dataset by repeatedly solving the optimization problem of Eq~(\ref{eq:td_findtau}) for varying set of features as described in Section~\ref{sec:strategy}; specifically, we choose the features as reported in Eq~\eqref{input}. The complete dataset contains $m = \, 900$ examples with $r = \left \{1, 2, 3 \right \}$, $h = \left \{ \frac{\sqrt{2}}{10}, \frac{\sqrt{2}}{20}, \frac{\sqrt{2}}{40} \right \}$ and values of $\Peclet_g$ randomly chosend in the range $[7, 70'710]$ (uniform distribution), which yields $\mu \in [10^{-5}, 10^{-1}]$. We summarize the details for generating the dataset in Table~\ref{tab:dataset}. Figure~\ref{fig:td_ann2_dataset} provides an overview of the dataset used for the training of the ANN, specifically the features $\mathbf x^{(i)}=[r, \, h,\, \mathbb Pe_g]$ and the target $\mathbf y^{(i)}=[\tau^*]$. In this case, the characteristic length $L$ is fixed and set equal to $L=1$.

\begin{figure}[H]
	\centering
	\begin{subfigure}{0.32\textwidth}
		\includegraphics[width=\textwidth]{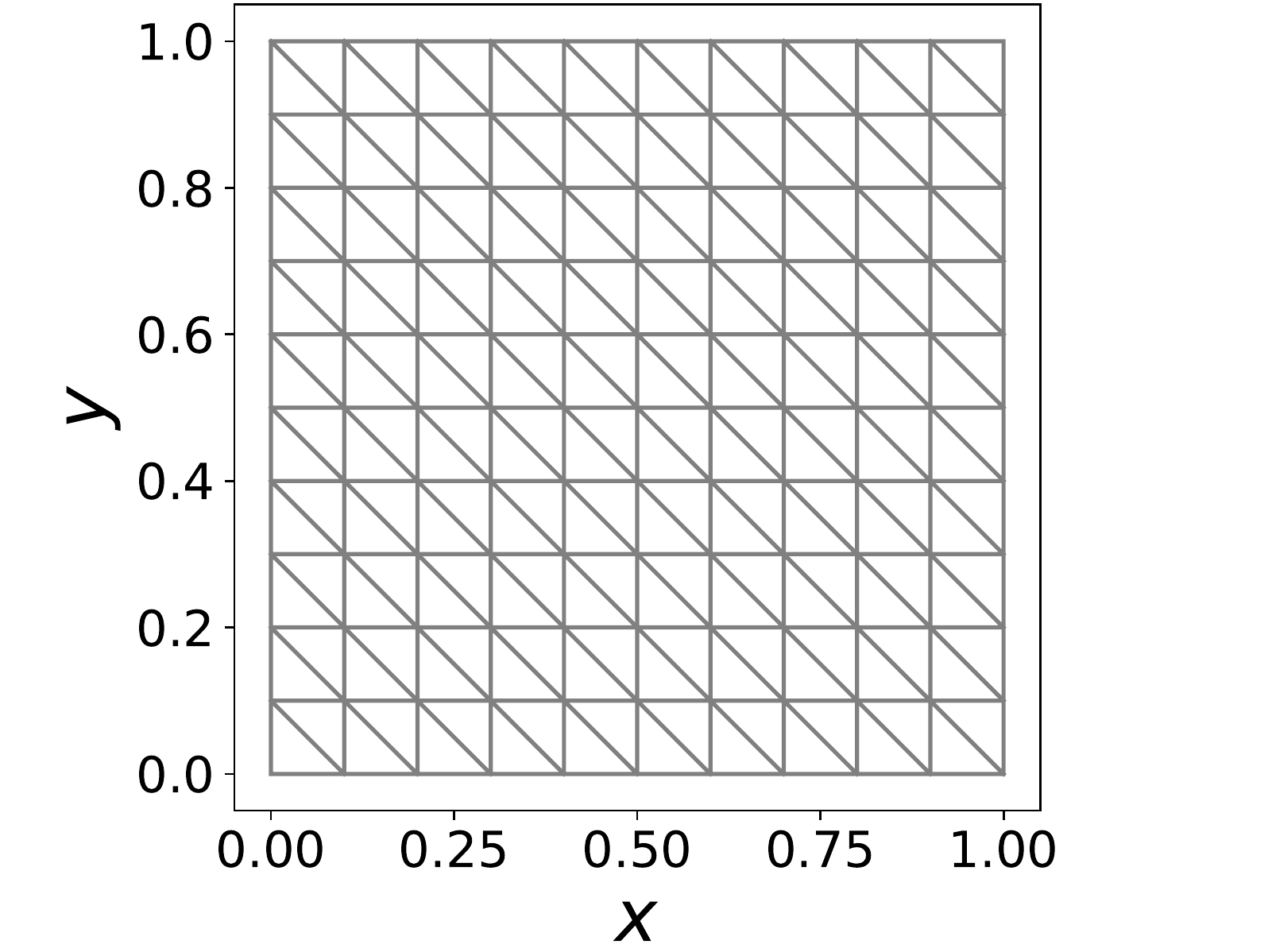}
		\caption{$h = \sqrt{2} / 10$}
	\end{subfigure}
	\hfill
	\begin{subfigure}{0.32\textwidth}
		\includegraphics[width=\textwidth]{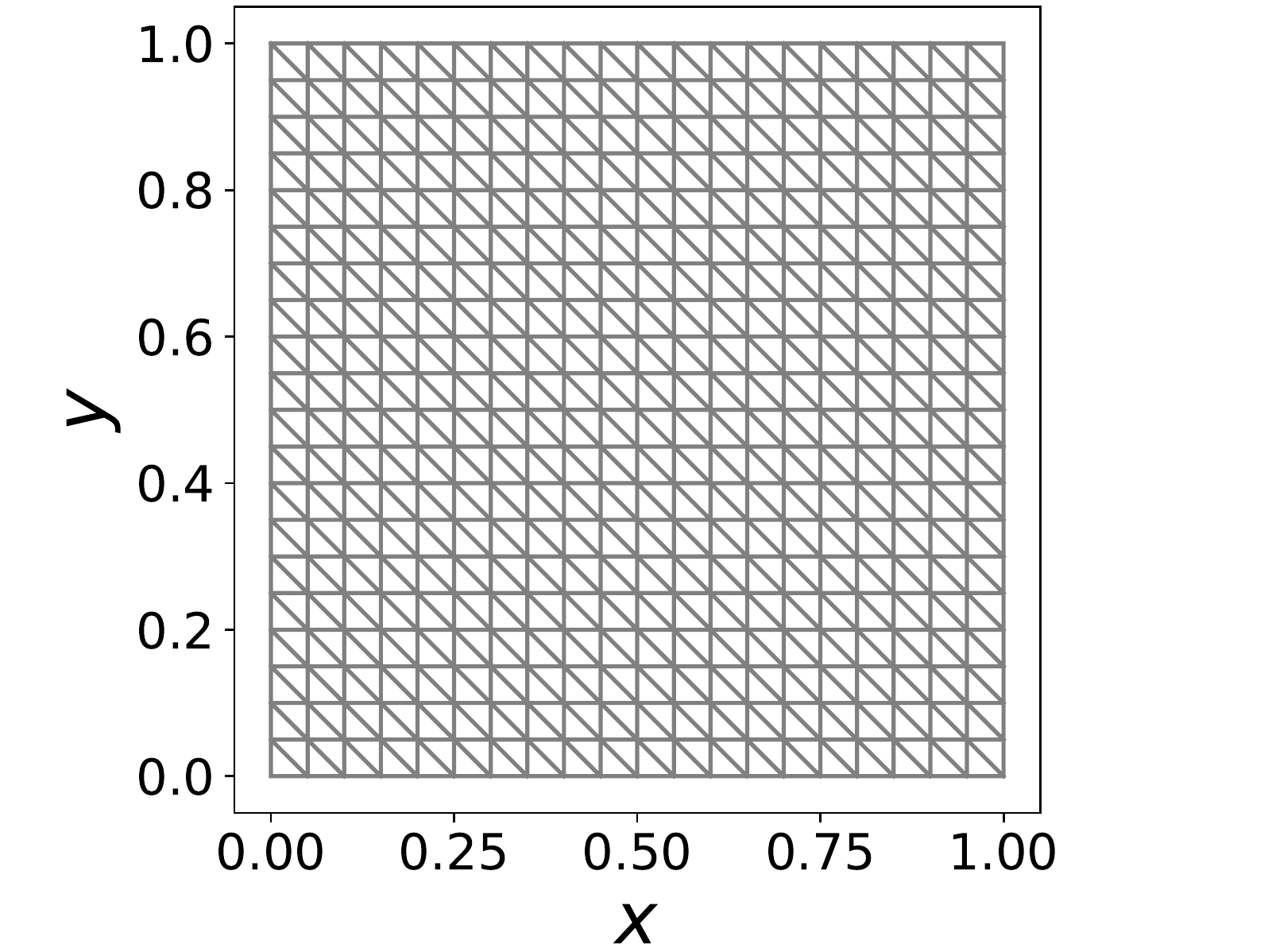}
		\caption{$h = \sqrt{2} / 20$}
	\end{subfigure}
	\hfill
	\begin{subfigure}{0.32\textwidth}
		\includegraphics[width=\textwidth]{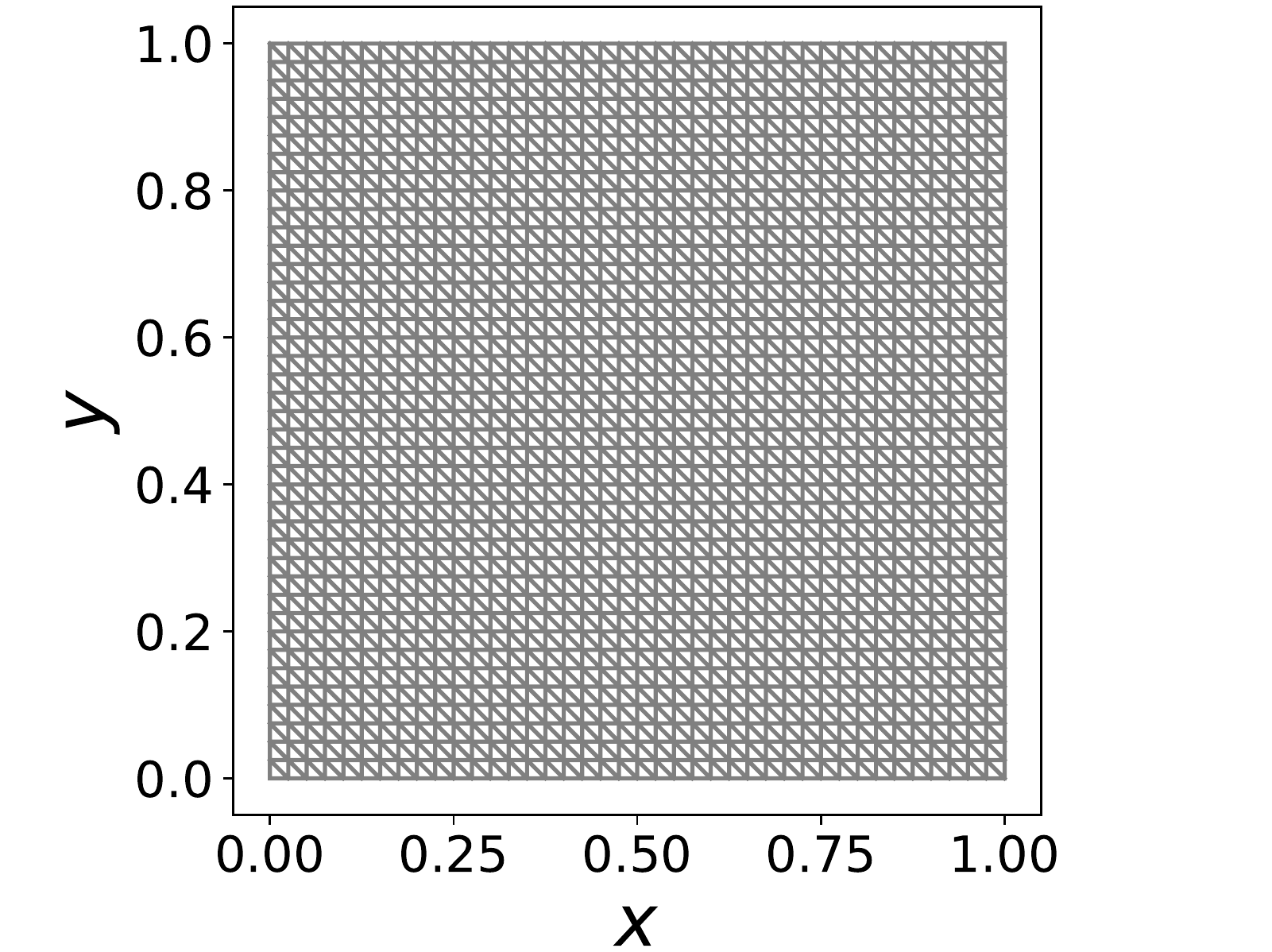}
		\caption{$h = \sqrt{2} / 40$}
	\end{subfigure}
	\caption{Structured meshes used for the FE approximation of the 2D advection-diffusion problem.}
	\label{fig:td_2d_mesh}
\end{figure}

\begin{figure}[H]
	\centering
	\includegraphics[width=0.65\textwidth]{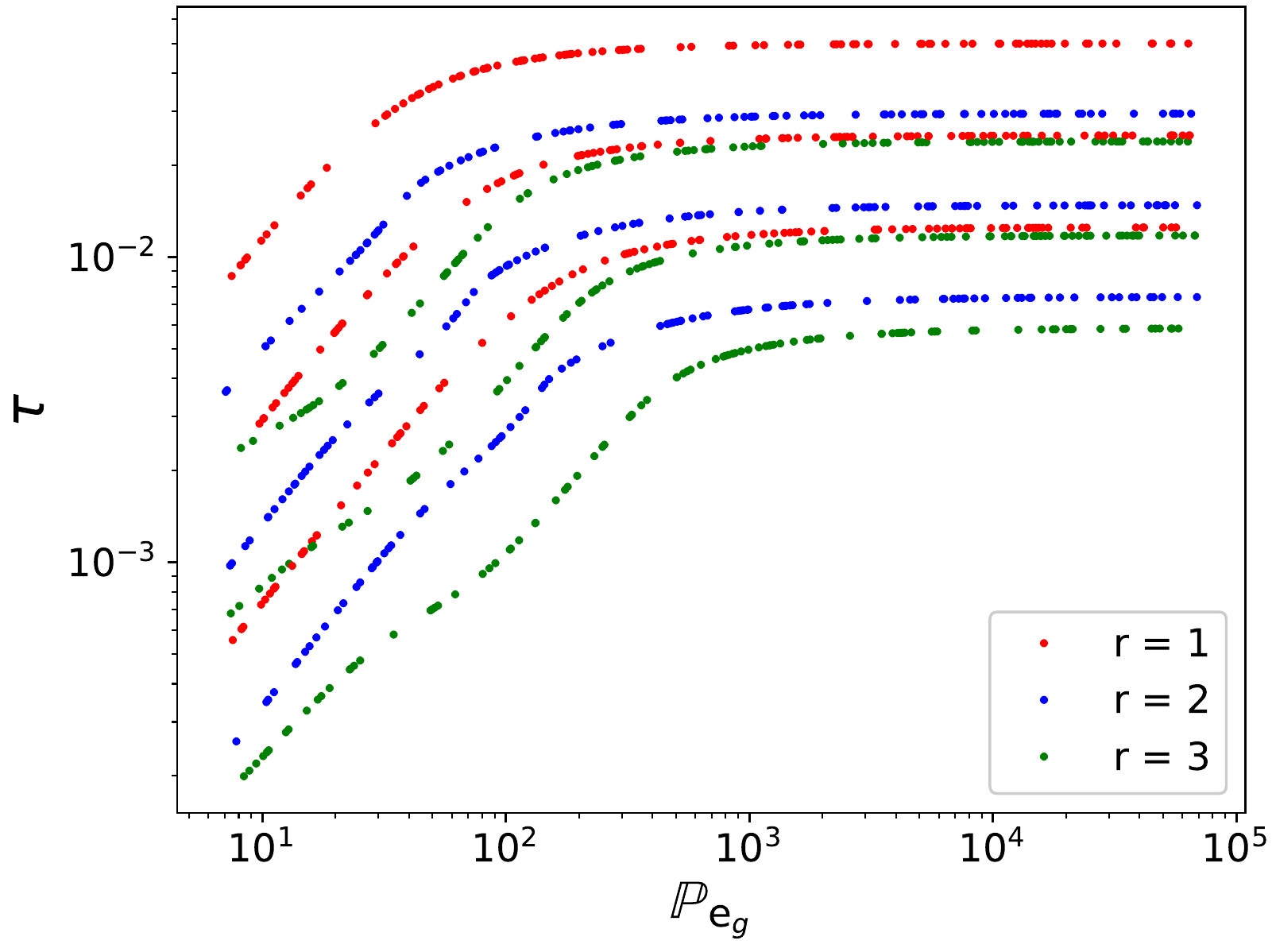}
	\caption{Visualization of the dataset used for the ANN training: target $\tau^*$ against feature $\Peclet_g$ colored by feature $r=1$ (red), $2$ (blue), and $3$ (green) for different values of the feature $h$ (increasing values of $h$ from bottom to top) as listed in Table~\ref{tab:dataset}.}
	\label{fig:td_ann2_dataset}
\end{figure}

\begin{table}[H]
\caption{Details of the dataset used for the ANN training.}
\label{tab:dataset}
\centering
\begin{tabular}{llllll}
\hline
\# Data & \# Training & \# Validation & \multirow{2}{*}{$r$} & \multirow{2}{*}{$h$} & \multirow{2}{*}{$\Peclet_g$} \\
set ($m$) & set & set & & & \\
\hline
\multirow{2}{*}{$900$} & \multirow{2}{*}{$720\, (80 \%)$} & \multirow{2}{*}{$180\, (20 \%)$} & \multirow{2}{*}{$\left \{1, 2, 3 \right \}$} &\multirow{2}{*}{  $\left \{ \frac{\sqrt{2}}{10}, \frac{\sqrt{2}}{20}, \frac{\sqrt{2}}{40} \right \} $}
& randomly \\
& & & & & in $[7, 70'710]$.\\\hline
\end{tabular}
\end{table}

We train a fully-connected feed-forward ANN on the generated dataset by using the open source library {\tt Keras} \cite{keras} built on top of {\tt TensorFlow} \cite{tensorflow}. We divide the dataset into two parts: a training dataset that takes $80\%$ of the examples to be used for the ANN training and a validation dataset that takes the remaining $20\%$. We choose the loss function as the mean squared error that measures, for each training feature, the squared mismatch between the prediction of the ANN $\widehat{y}^{(j)}$ and the actual target $y^{(j)}$. Specifically, the loss function is defined as
\begin{equation}
    \mathcal J = \frac{1}{2m} \sum_{j = 1}^m \left ( \widehat{y}^{(j)} - y^{(j)}\right )^2.
    \label{cost}
\end{equation}
We normalize the features by subtracting their sample mean and dividing by their sample standard deviation in order to help the weights to better adapt to the different scales of the features. Moreover, the targets and the feature $\Peclet_g$ need special care as they are distributed over a wide range of values. Thus, we normalize them using by applying a base $10$ logarithm. The normalized features and targets read:
\begin{equation*}
    \widetilde{\mathbf x} = \left[ 
    \begin{array}{c}
        \frac{r - \overline{r}}{\sigma_r}, \; \frac{h - \overline{h}}{\sigma_h}, \; \frac{\log_{10}(\Peclet_g) - \overline{\log_{10}(\Peclet_g)}}{\sigma_{\log_{10}(\Peclet_g)}}
    \end{array} \right], \qquad
    \widetilde{\mathbf y} = \left[ 
    \begin{array}{c}
        - \log_{10} (\tau^*)
    \end{array} \right],
\end{equation*}
where $\overline{r}$, $\overline{h}$, $\overline{\log_{10}(\Peclet_g)}$ are the sample mean of the training features $r$, $h$ and the logarithm of $\Peclet_g$ respectively, while $\sigma_r$, $\sigma_h$ and $\sigma_{\log_{10}(\Peclet_g)}$ are their sample standard deviations.

\begin{figure}[H]
    \centering
    \includegraphics[width=0.8\textwidth]{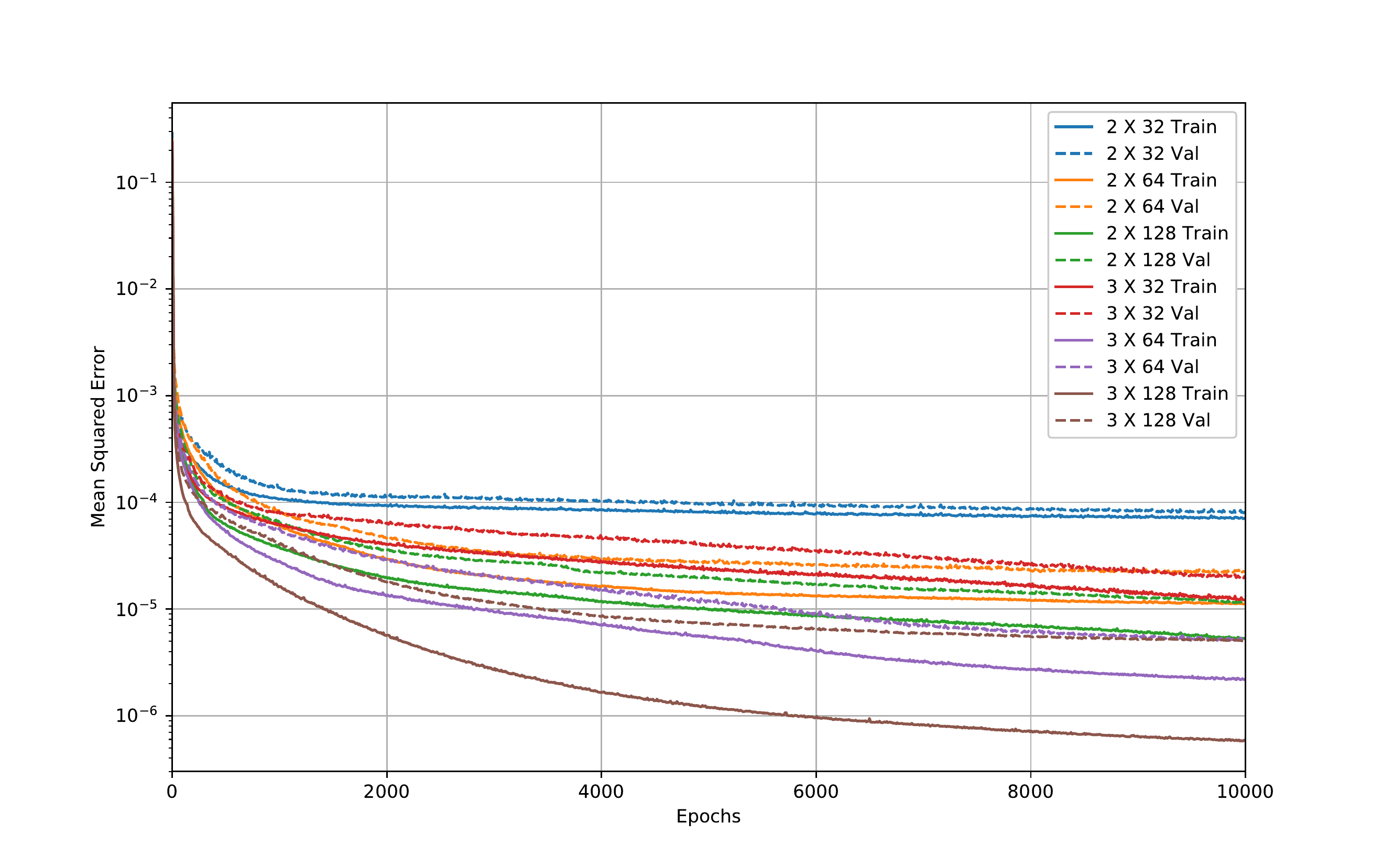}
    \caption{ Comparison of training and validation errors over training epochs with different architectures: 2 and 3 hidden layers and 32, 64 and 128 nodes per layer. }
    \label{fig:training-validation-error}
\end{figure}

In order to find the best performing ANN architecture, we carried out a study by testing different numbers of hidden layers, nodes per layer, optimization algorithm, its learning rate, and batch size.  A comparison of the loss function with different architectures is given in Figure~\ref{fig:training-validation-error}. Using 3 hidden layers is beneficial in terms of validation error drops, while using more than 64 nodes per layer does not bring to considerable advantages.  Finally, we choose an ANN with $3$ hidden layers, $64$ nodes per layer and an output layer with a single node, all using a rectified linear unit (\textit{ReLU}) activation function. We display the ANN architecture in Figure~\ref{fig:ann_architecture}. We trained the ANN with the SGD optimization algorithm, a constant learning rate of $0.01$, and mini-batch size of $32$ samples. Moreover, we employed a momentum of $0.9$ in the optimization algorithm to update the weights. The trained ANN is available in the {\tt GitLab} repository \cite{ann_repo}.

\begin{figure}[H]
	\centering
	\includegraphics[width=0.5\textwidth]{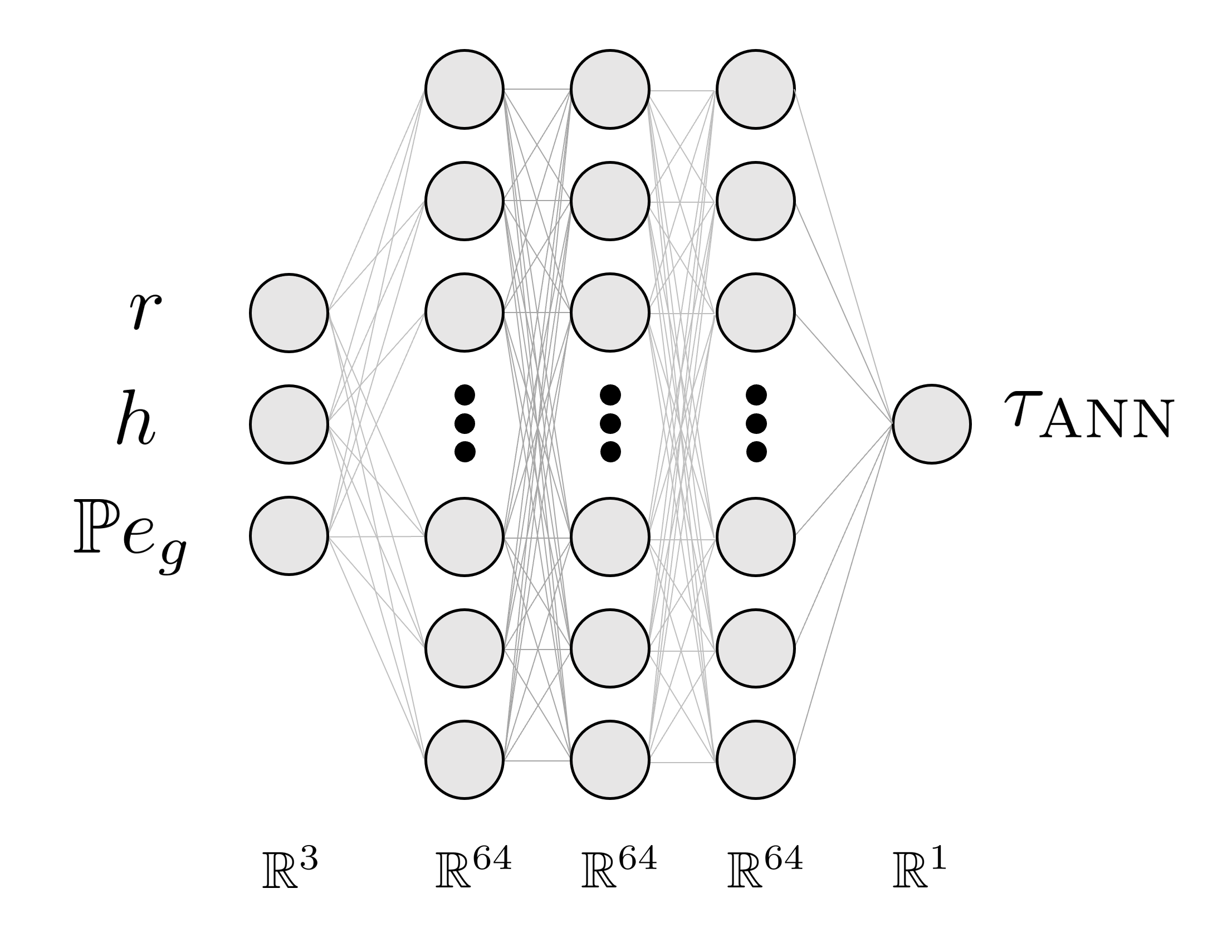}
	\caption{Architecture of the feed-forward fully-connected ANN.}
	\label{fig:ann_architecture}
\end{figure}

Regarding the computational efficiency of the proposed strategy, we stress the clear distinction between the offline phase (dataset generation and ANN's training) and the online phase, where we use the ANN to predict a new stabilization parameter -- alongside the use of the FE solver -- for unseen input parameter values. The most demanding part of the strategy is the dataset generation, requiring approximately 2 h on a standard laptop to repeatedly solve the optimization problem 900 times. Moreover, the ANN's training phase required approximately 10'. The testing (online) phase, which is the one that is performed for application purposes, is considerably inexpensive, requiring only the real-time evaluation of a composition of linear functions, comparable with the evaluation of the empirical relation that brings to the theoretical values ( few milliseconds).

\subsection{Predictions of the stabilization parameter by ANN}
Now, we compare the predictions of the ANN with the theoretical stabilization parameter $\widetilde{\tau}_r$ of Eq~\eqref{eq:tau_theory_r} \cite{Quarteroni_2017,galeao2004finite}. We show in Figure~\ref{fig:td_2d_predictions_comparison} (left) the stabilization parameter $\tau_{\mathrm{ANN}}$ predicted by the ANN by varying mesh size $h$ and global P\'eclet number $\mathbb Pe_g$. For comparison, we report in Figure~\ref{fig:td_2d_predictions_comparison} (right) the corresponding theoretical stabilization parameter $\widetilde{\tau}_r$. We observe that the overall behavior of the trained ANN's predictions are qualitatively similar of the theoretical stabilization parameter $\widetilde{\tau}_r$. Nevertheless, it can be inferred that the values of the stabilization parameters are almost completely matched with linear FE ($r=1$), while they quantitatively differ for $r=2$ and $r=3$. This was expected as $\widetilde{\tau}_r$ for $r>1$ is an empirical extension of the formula for the case $r=1$.

The ANN allows to make predictions with features outside the range of values for which it has been trained, that is for unseen values of such features. With this aim, we report in Figure~\ref{fig:td_2d_predictions_comparison_r4} the comparison of the ANN's predictions $\tau_\mathrm{ANN}$ with $\widetilde{\tau}_r$ for the FE degree $r = 4$. This comparison shows a clear difference between the theoretical and ANN's $\tau$ for $r = 4$ even though the general trend is maintained similar. 
\begin{figure}[H]
    \centering
    \begin{subfigure}{0.49\textwidth}
         \includegraphics[width=\textwidth]{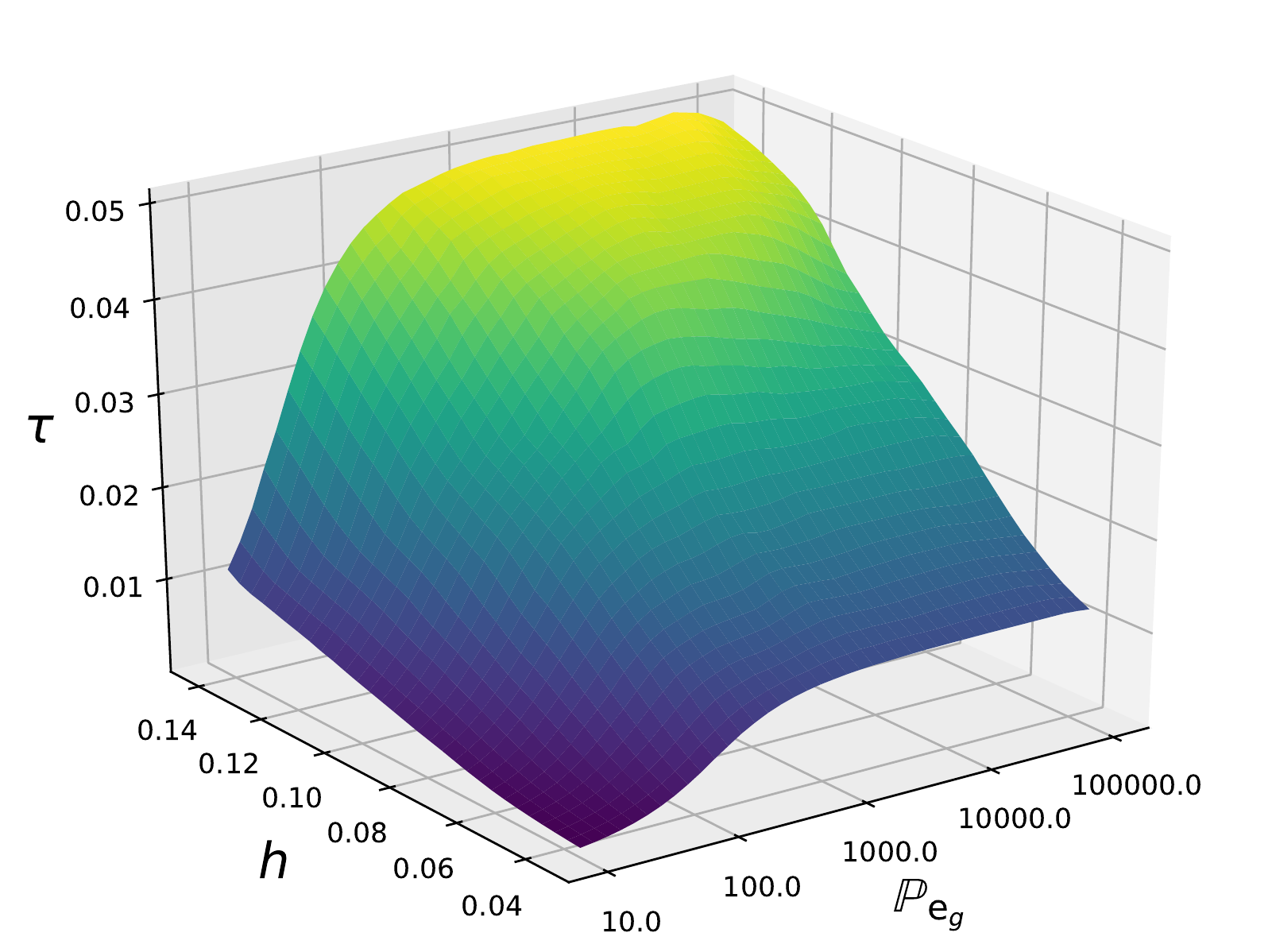}
         \caption{$\tau_{\mathrm{ANN}}$ for $r = 1$}
    \end{subfigure}
    \begin{subfigure}{0.49\textwidth}
         \includegraphics[width=\textwidth]{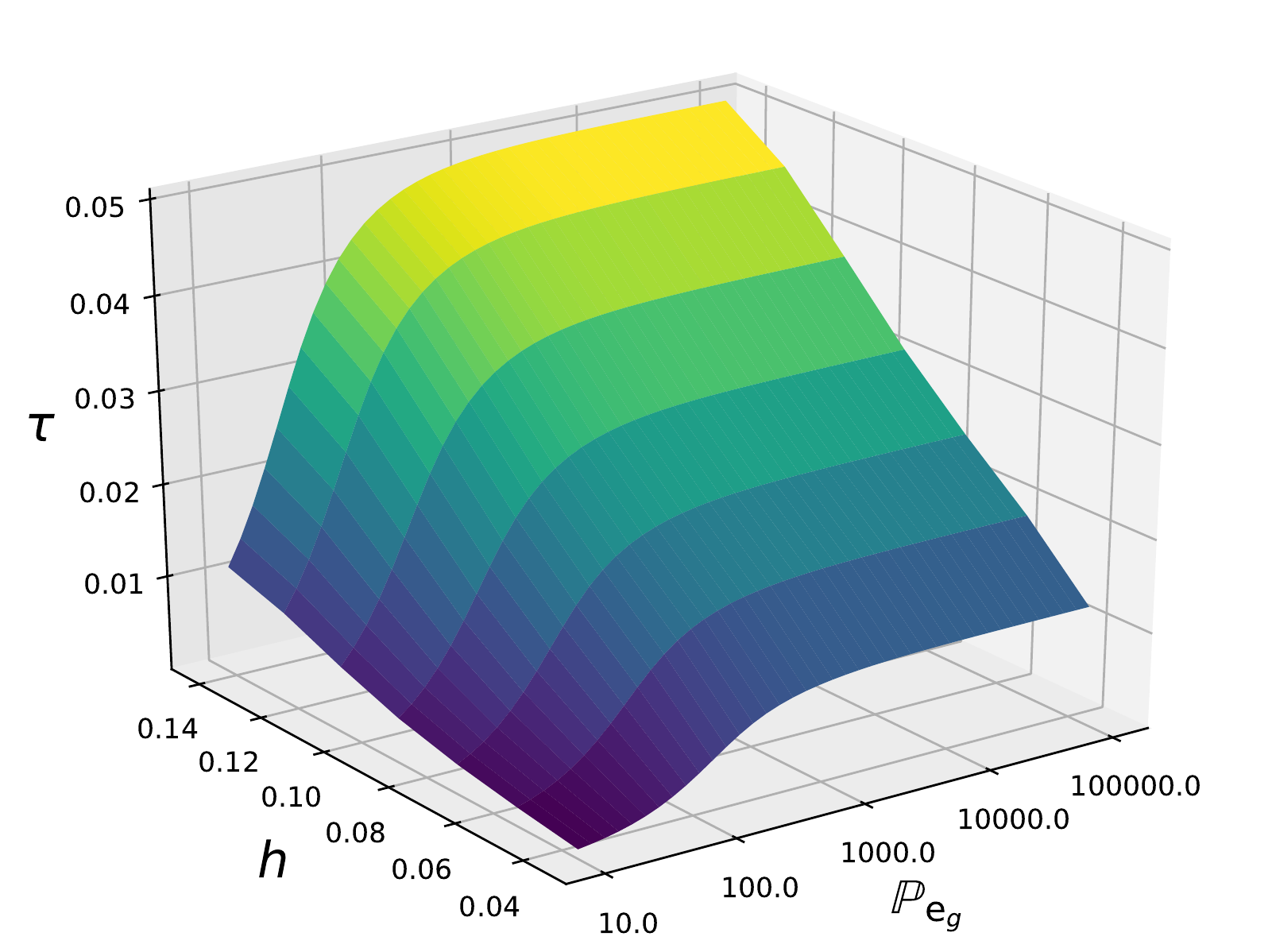}
         \caption{$\widetilde{\tau}_r$ for $r = 1$}
    \end{subfigure}
    \begin{subfigure}{0.49\textwidth}
         \includegraphics[width=\textwidth]{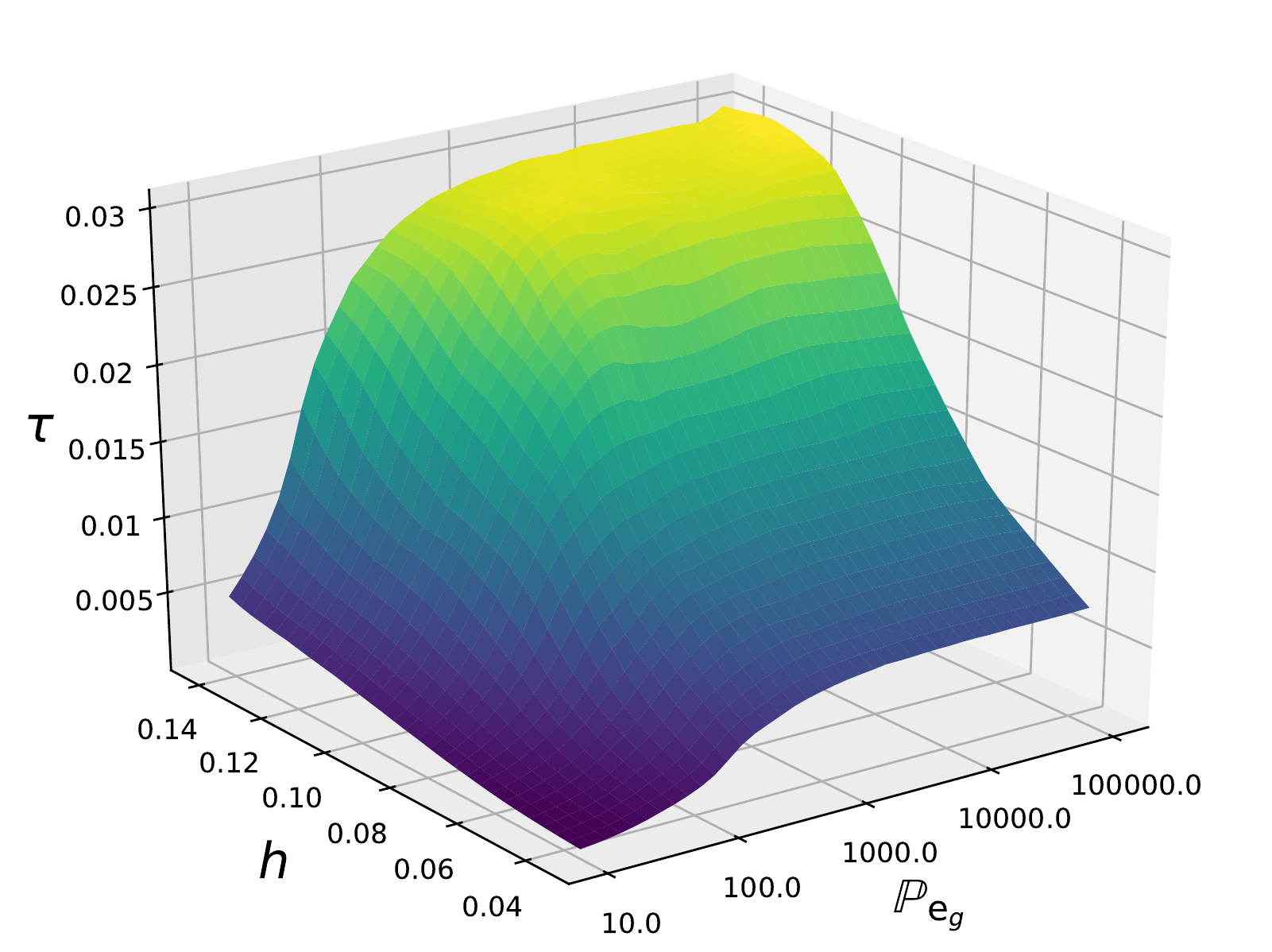}
         \caption{$\tau_{\mathrm{ANN}}$ for $r = 2$}
    \end{subfigure}
    \begin{subfigure}{0.49\textwidth}
         \includegraphics[width=\textwidth]{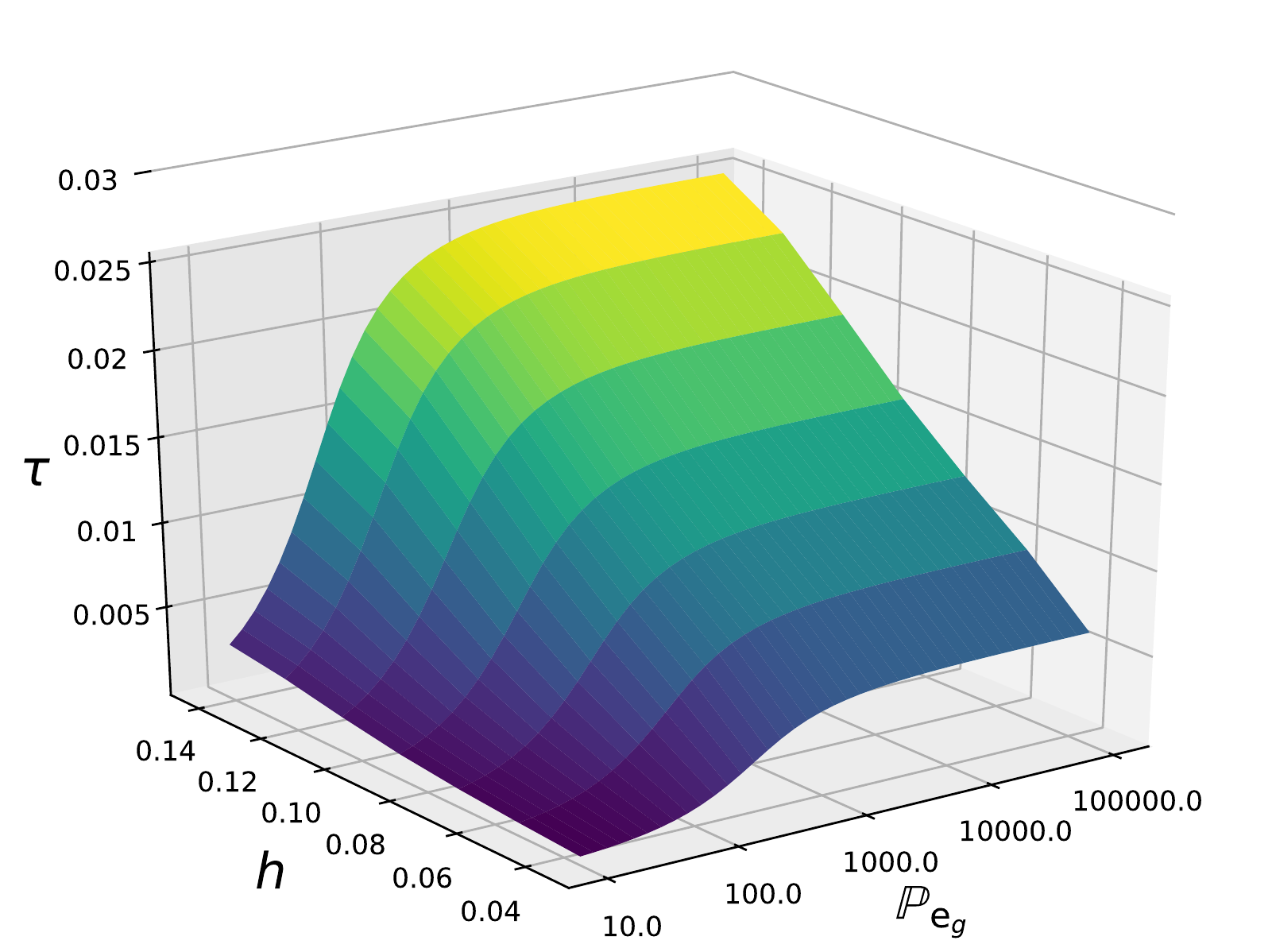}
         \caption{$\widetilde{\tau}_r$ for $r = 2$}
    \end{subfigure}
    \begin{subfigure}{0.49\textwidth}
         \includegraphics[width=\textwidth]{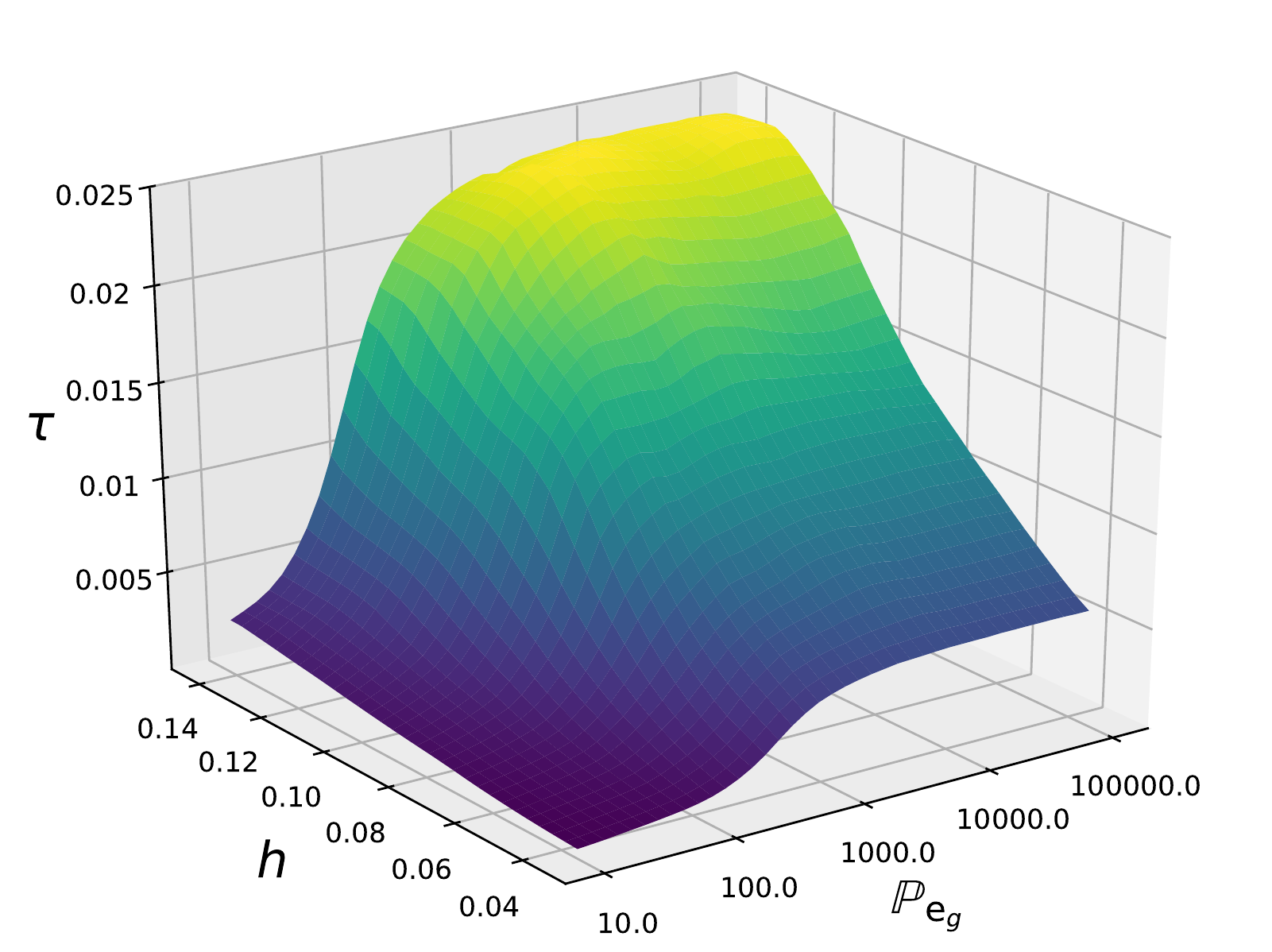}
         \caption{$\tau_{\mathrm{ANN}}$ for $r = 3$}
    \end{subfigure}
    \begin{subfigure}{0.49\textwidth}
         \includegraphics[width=\textwidth]{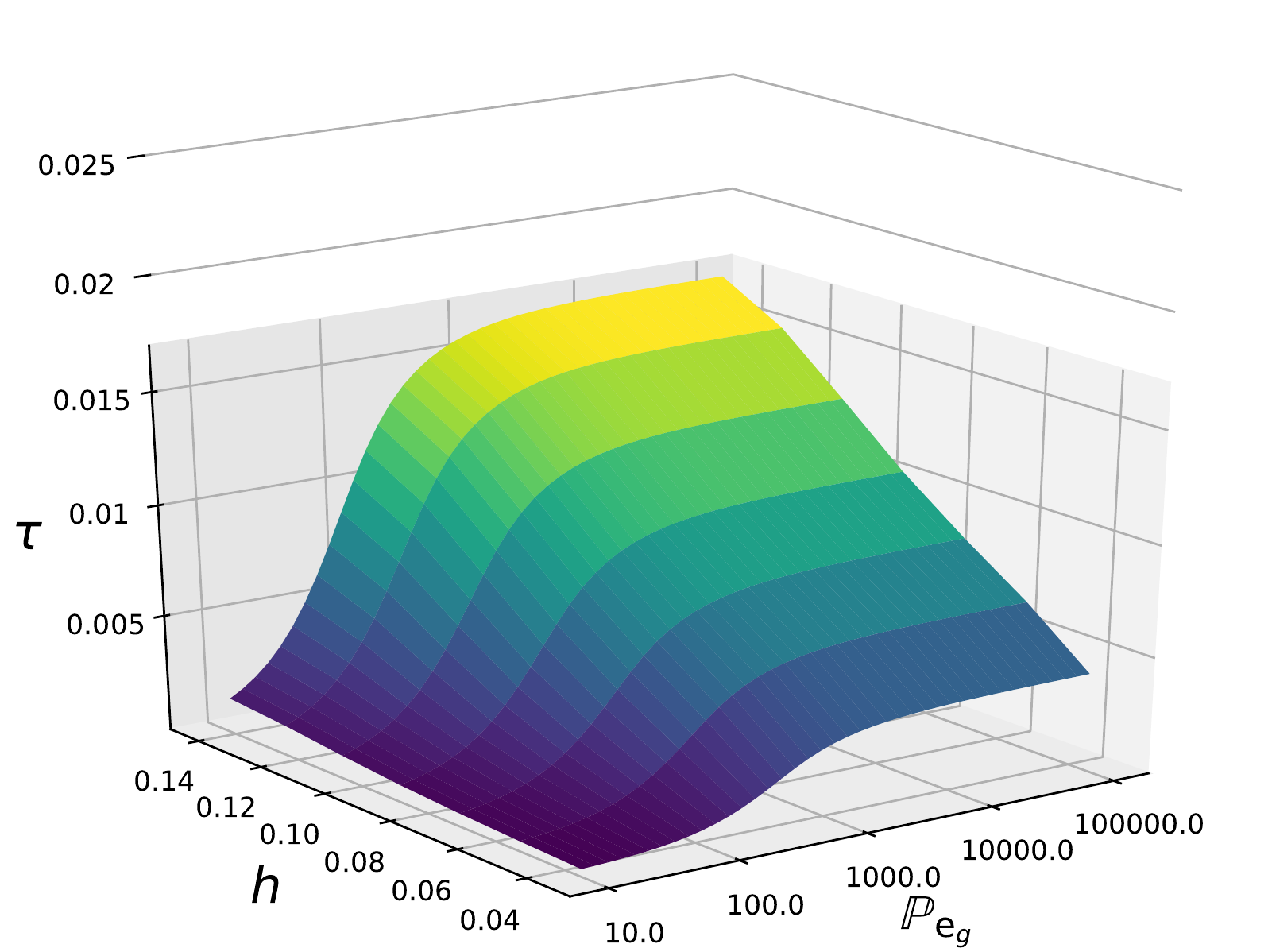}
         \caption{$\widetilde{\tau}_r$ for $r = 3$}
    \end{subfigure}
    \caption{Stabilization parameter $\tau_{\mathrm{ANN}}$ predicted by the ANN and theoretical one $\widetilde{\tau}_r$ for varying $\Peclet_g$ and $h$ at different FE degrees $r=1,2$, an $3$ for the 2D advection-diffusion problem of Section~\ref{sec:training_2D}.}
    \label{fig:td_2d_predictions_comparison}
\end{figure}


\begin{figure}[H]
    \centering
    \begin{subfigure}{0.49\textwidth}
         \includegraphics[width=\textwidth]{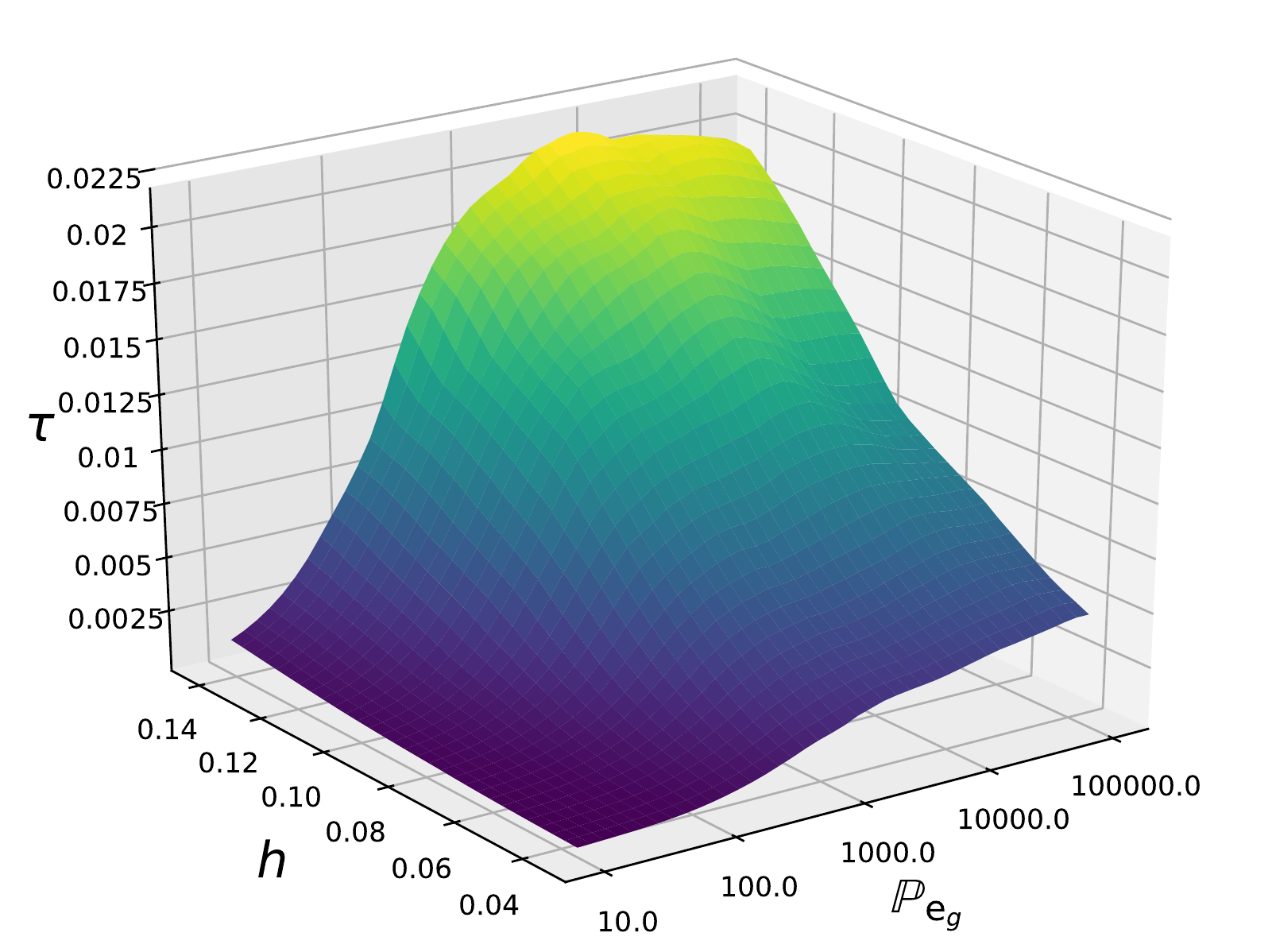}
         \caption{ANN's predicted $\tau$ at $r = 4$}
    \end{subfigure}
    \begin{subfigure}{0.49\textwidth}
         \includegraphics[width=\textwidth]{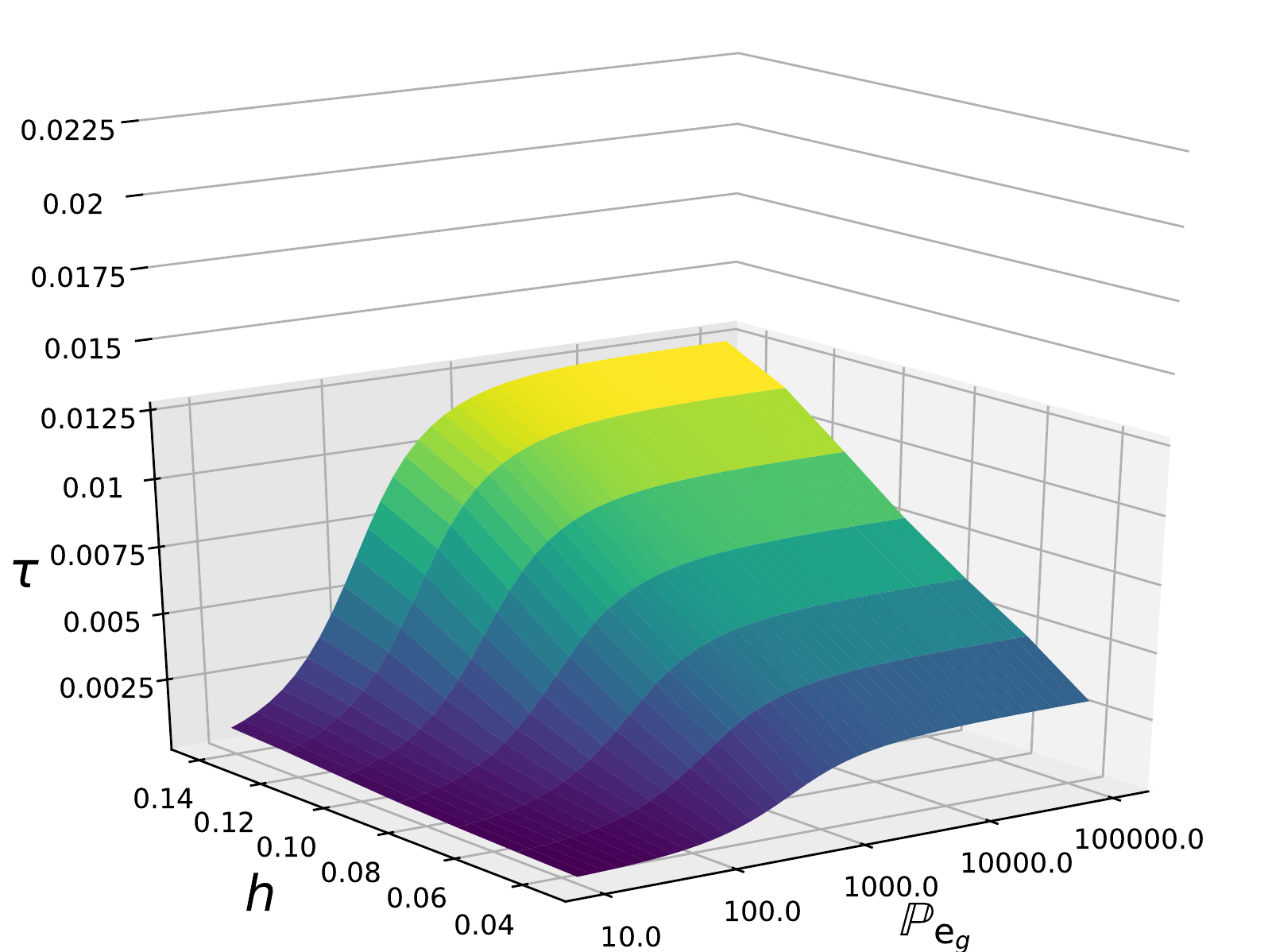}
         \caption{Theoretical $\tau$ at $r = 4$}
    \end{subfigure}
    \caption{Stabilization parameter $\tau_{\mathrm{ANN}}$ predicted by the ANN and theoretical one $\widetilde{\tau}_r$ for varying $\Peclet_g$ and $h$ with FE degree $r=4$ for the 2D advection-diffusion problem of Section~\ref{sec:training_2D}.}
    \label{fig:td_2d_predictions_comparison_r4}
\end{figure}

\subsubsection{Test 1:  predictions for the problem used in the ANN's training}
We compare the numerical solutions $u^\mathrm{ANN}_h$ and $\widetilde u_h$ obtained by means of the SUPG stabilization method with the parameter $\tau_\mathrm{ANN}$ predicted by the ANN and the theoretical one $\widetilde{\tau}_r$, respectively; the comparison also involves the exact solution $u$~(\ref{eq:td_2d_exact}) of the 2D advection-diffusion problem used for the training of the ANN. In particular, we display the comparison of the former 2D solutions in Figure~\ref{fig:td_2d_solutions_extracted} along the line $(1-h,y)$ for any $y \in [0,1]$ with: (a) $\Peclet_h = 2$, $r=1$, and $h=\sqrt{2}/10$ (Figure~\ref{fig:td_2d_solutions_extracted1}); (b) $\Peclet_h = 500$, $r=3$, and and $h=\sqrt{2}/20$ (Figure~\ref{fig:td_2d_solutions_extracted2}). We notice that 
in both the cases, the numerical solution $u_h^\mathrm{ANN}$ involving the stabilization parameter $\tau_\mathrm{ANN}$ provides more accurate results than with the theoretical stabilization paramter $\widetilde{\tau}_r$. In particular, in the case (a), $\widetilde{u}_h$ involves a much smoother boundary layer and overshoots the exact solutions $u$, conversely to the nearly nodally exact numerical solution $u_h^\mathrm{ANN}$. In the case~(b),
$u_h^\mathrm{ANN}$ provides a much better representation of the bundary layer, without the undershooting of the solution $u$ exhibited by $\widetilde u_h$. Furthermore, we report in Table ~\ref{tab:errors} absolute errors between the numerical solutions ($u_h^\mathrm{ANN}$ and $\widetilde{u}_h$) and the exact one $u$. The errors, computed in $L^2(\Omega)$ and $H^1(\Omega)$ norms, show that the ANN-based SUPG stabilization method very often produces more accurate results than the ones obtained with theoretical stabilization parameter $\widetilde{\tau}_r$, especially for $r=3$.

Moreover, to assess the ability of the ANN to predict the stabilization parameter out of the training range, we compare in Figure \ref{fig:td_2d_solutions_extracted_r4} the numerical solutions for high Péclet in the case $r=4$. We observe that $u_h^\mathrm{ANN}$ is more accurate than $\widetilde{u}_h$ even with FE degree out of the training range. \newpage
\begin{figure}[H]
    \centering
    \begin{subfigure}{0.4\textwidth}
         \includegraphics[width=\textwidth]{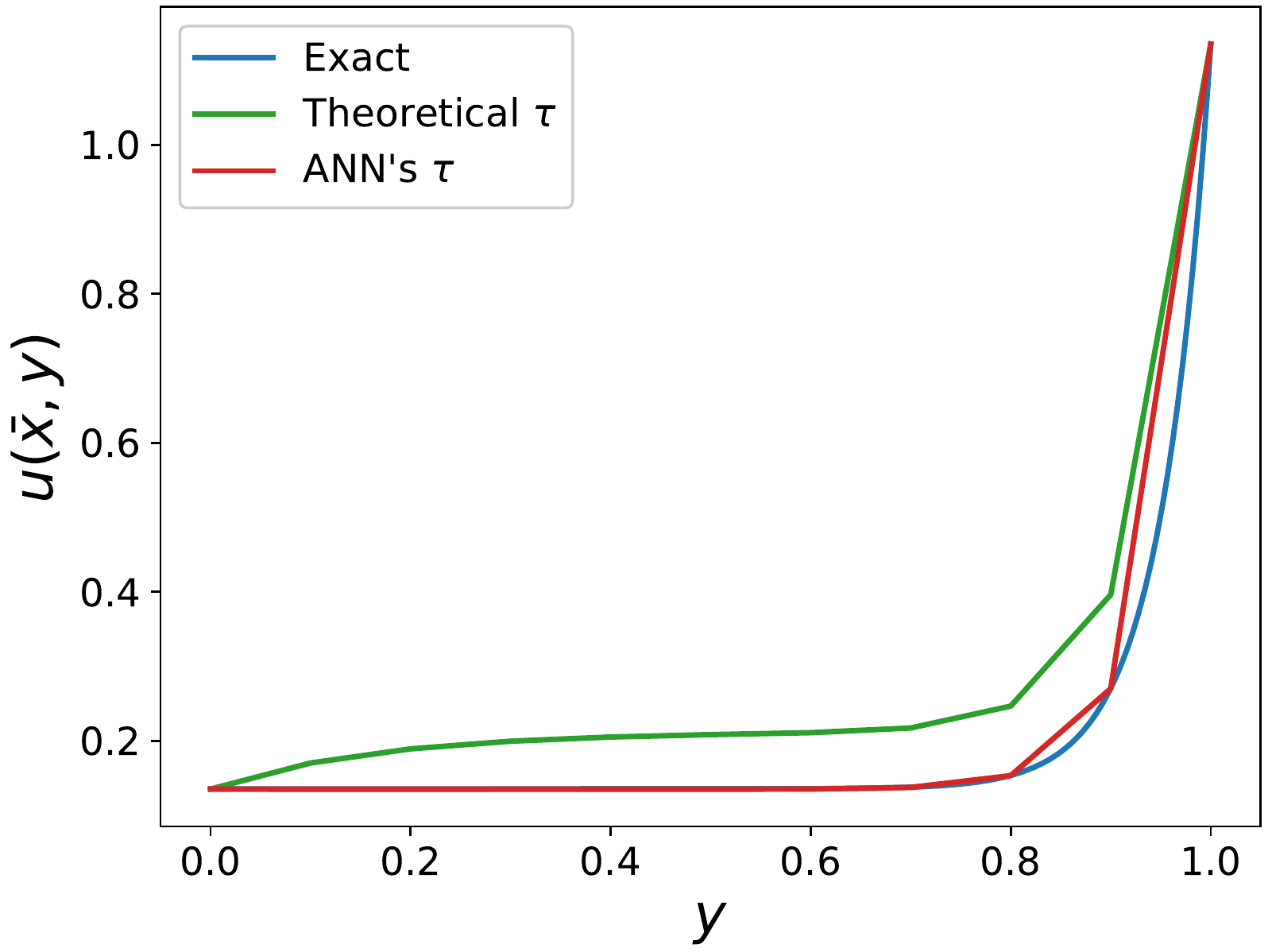}
         \caption{$\Peclet_h = 2$, $h = \sqrt{2}/10$, $r = 1$}
         \label{fig:td_2d_solutions_extracted1}
    \end{subfigure}
    \begin{subfigure}{0.4\textwidth}
         \includegraphics[width=\textwidth]{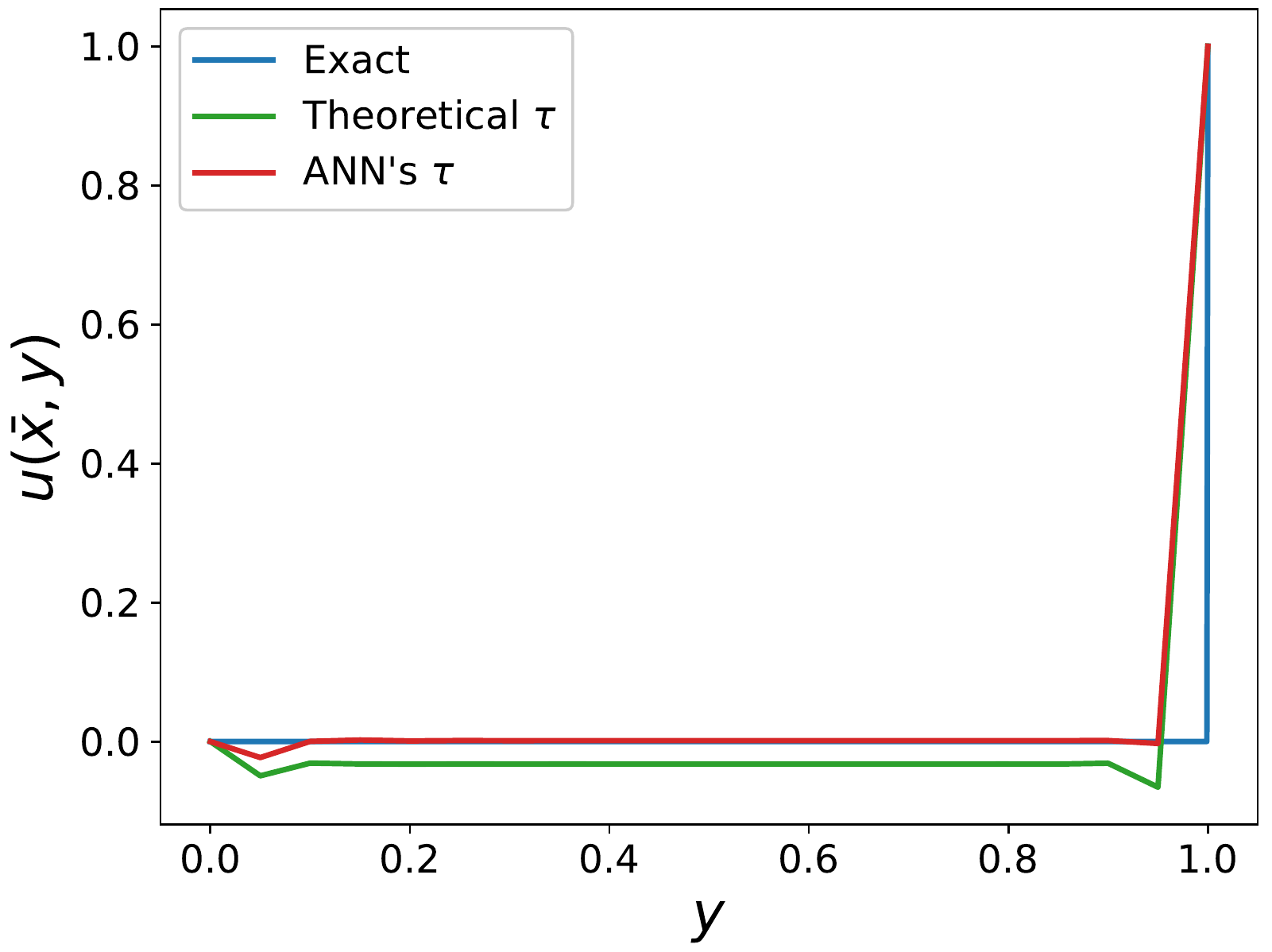}
         \caption{$\Peclet_h = 500$, $h = \sqrt{2}/20$, $r = 3$}
         \label{fig:td_2d_solutions_extracted2}
    \end{subfigure}
    \caption{Test 1:  comparison of the solutions $u$ \eqref{eq:td_2d_exact} (blue), $\widetilde{u}_h$ (green), and $u_h^\mathrm{ANN}$ (red) of the 2D advection-diffusion problem along the line $(1-h,y)$ for any $y \in [0,1]$.}
    \label{fig:td_2d_solutions_extracted}
\end{figure}

\begin{table}[H]
\centering
\caption{Test 1:  comparison of the errors in norms $L^2$ and $H^1$ between the numerical solutions ($u_h^\mathrm{ANN}$ and $\widetilde{u}_h$) and the exact one $u$ \eqref{eq:td_2d_exact} of the 2D advection-diffusion problem used for the training (Section~\ref{sec:training_2D}) for different values of the features.}
\label{tab:errors}
\begin{tabular}{lllllll}\hline
&  &  & \multicolumn{2}{c}{$\widetilde{\tau}_r$ ($e=\widetilde u_h - u$)}
& \multicolumn{2}{c}{${\tau}_\mathrm{ANN}$ ($e=u^\mathrm{ANN}_h - u$)}
\\
$r$ & $h$ & $\Peclet_h$ & $||e||_{L^2(\Omega)}$ & $||e||_{H^1(\Omega)}$ & $||e||_{L^2(\Omega)}$ & $||e||_{H^1(\Omega)}$ 
\\
\hline
1 & $\sqrt{2}/10$ & 2 & $9.56 \cdot 10^{-2}$ & $7.93 \cdot 10^{-1}$ & $6.49 \cdot 10^{-2}$ & $5.31 \cdot 10^{-1}$\\
2 & $\sqrt{2}/10$ & 2 & $6.48 \cdot 10^{-2}$ & $5.31 \cdot 10^{-1}$ & $6.50 \cdot 10^{-2}$ & $5.33 \cdot 10^{-1}$\\
3 & $\sqrt{2}/10$ & 2 & $6.50 \cdot 10^{-2}$ & $5.33 \cdot 10^{-1}$ & $6.49 \cdot 10^{-2}$ & $5.32 \cdot 10^{-1}$\\
\hline 
1 & $\sqrt{2}/20$ & 500 & $2.99 \cdot 10^{-4}$ & $5.30 \cdot 10^{-3}$ & $3.11 \cdot 10^{-5}$ & $5.50 \cdot 10^{-4}$\\
2 & $\sqrt{2}/20$ & 500 & $2.71 \cdot 10^{-2}$ & $4.86 \cdot 10^{-1}$ & $2.50 \cdot 10^{-2}$ & $4.48 \cdot 10^{-1}$\\
3 & $\sqrt{2}/20$ & 500 & $1.05 \cdot 10^{-2}$ & $1.90 \cdot 10^{-1}$ & $1.68 \cdot 10^{-3}$ & $3.46 \cdot 10^{-2}$\\\hline
\end{tabular}
\end{table}


\begin{figure}[H]
    \centering
         \includegraphics[width=0.55\textwidth]{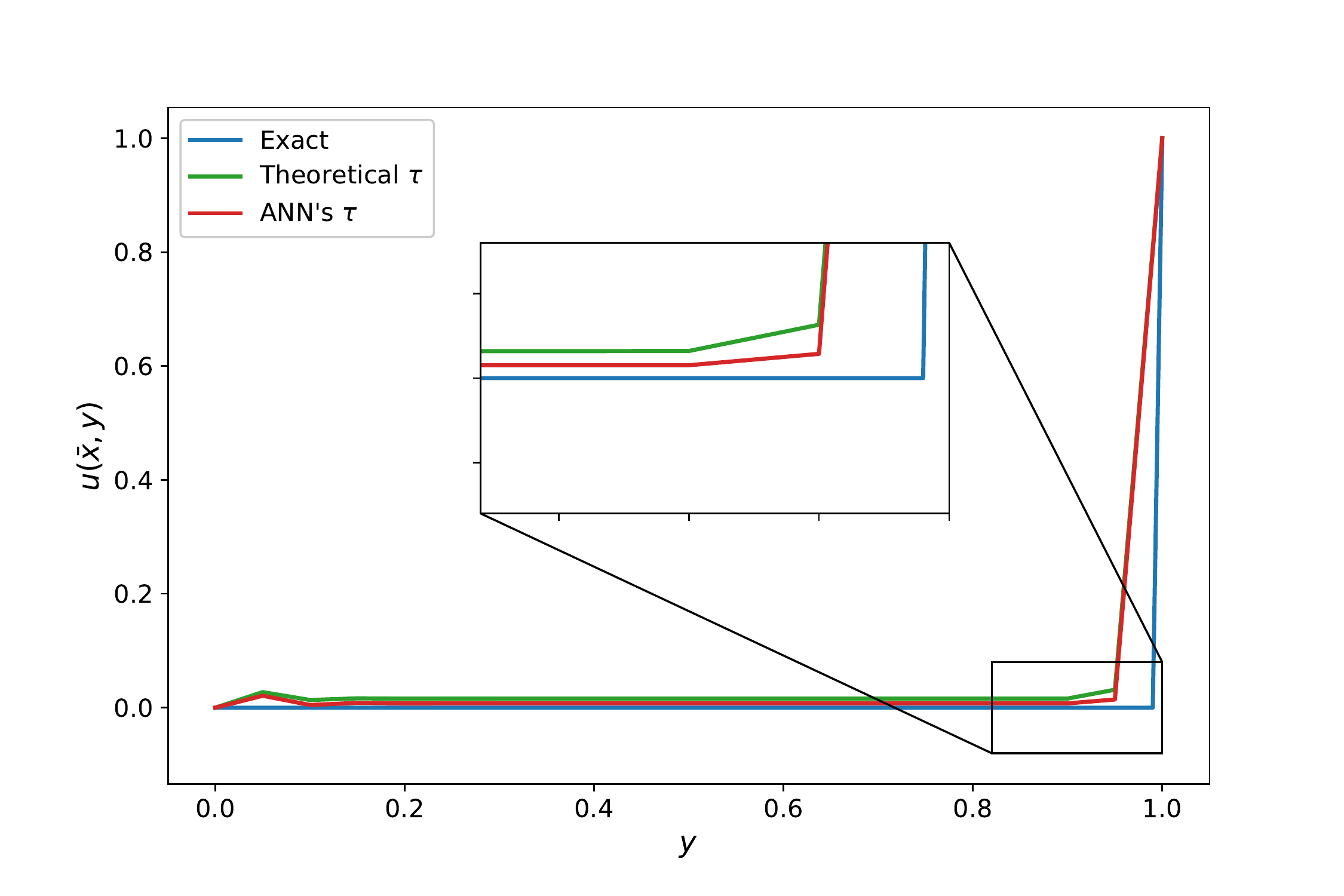}
    \caption{Test 1: comparison of the solutions $u$ \eqref{eq:td_2d_exact} (blue), $\widetilde{u}_h$ (green), and $u_h^\mathrm{ANN}$ (red) of the 2D advection-diffusion problem with $\Peclet_g = 7'071$ and $h = \sqrt{2}/20$,  along the line $(1-h,y)$ for any $y \in [0,1]$ for $r=4$ (outside of the training range).}
    \label{fig:td_2d_solutions_extracted_r4}
\end{figure}

\subsubsection{Test 2:  predictions for an unseen problem with constant forcing term}
\label{sec:newadvection-diffusion}
We check the model generalization of the ANN trained for the 2D advection-diffusion problem with exact solution~$u$ of Eq~\eqref{eq:td_2d_exact} by predicting $\tau_\mathrm{ANN}$ for an unseen 2D advection-diffusion problem. In particular, we consider the advection-diffusion problem of Eq~\eqref{eq:transport_diffusion} in $\Omega=(0,1)^2$, with $f=1$ and $\boldsymbol{\beta} = (1, 1)$.  We prescribe the following exact solution on the whole boundary $\partial \Omega$:
\begin{equation}
    \label{eq:td_2d_exact_new}
    u(x, y) = \frac{1}{2}(x+y) + \frac{1 - \frac{1}{2} (e^{(x / \mu)} + e^{(y / \mu)})}{e^{(1 / \mu)} - 1}.
\end{equation}

We compare in Figure~\ref{fig:td_2d_solutions_extracted_new} the numerical solutions $u^\mathrm{ANN}_h$ and $\widetilde u_h$ with the novel exact solution $u$  for different Péclet and mesh sizes. As for the previous numerical tests, the stabilization parameter $\tau_\mathrm{ANN}$ predicted by the network provides more accurate numerical solutions $u_h^\mathrm{ANN}$ than with the theoretical stabilization parameter $\widetilde{\tau}_r$. In particular, the boundary layers and overall solution behaviours are better represented.
\begin{figure}[H]
    \centering
    \begin{subfigure}{0.49\textwidth}
         \includegraphics[width=\textwidth]{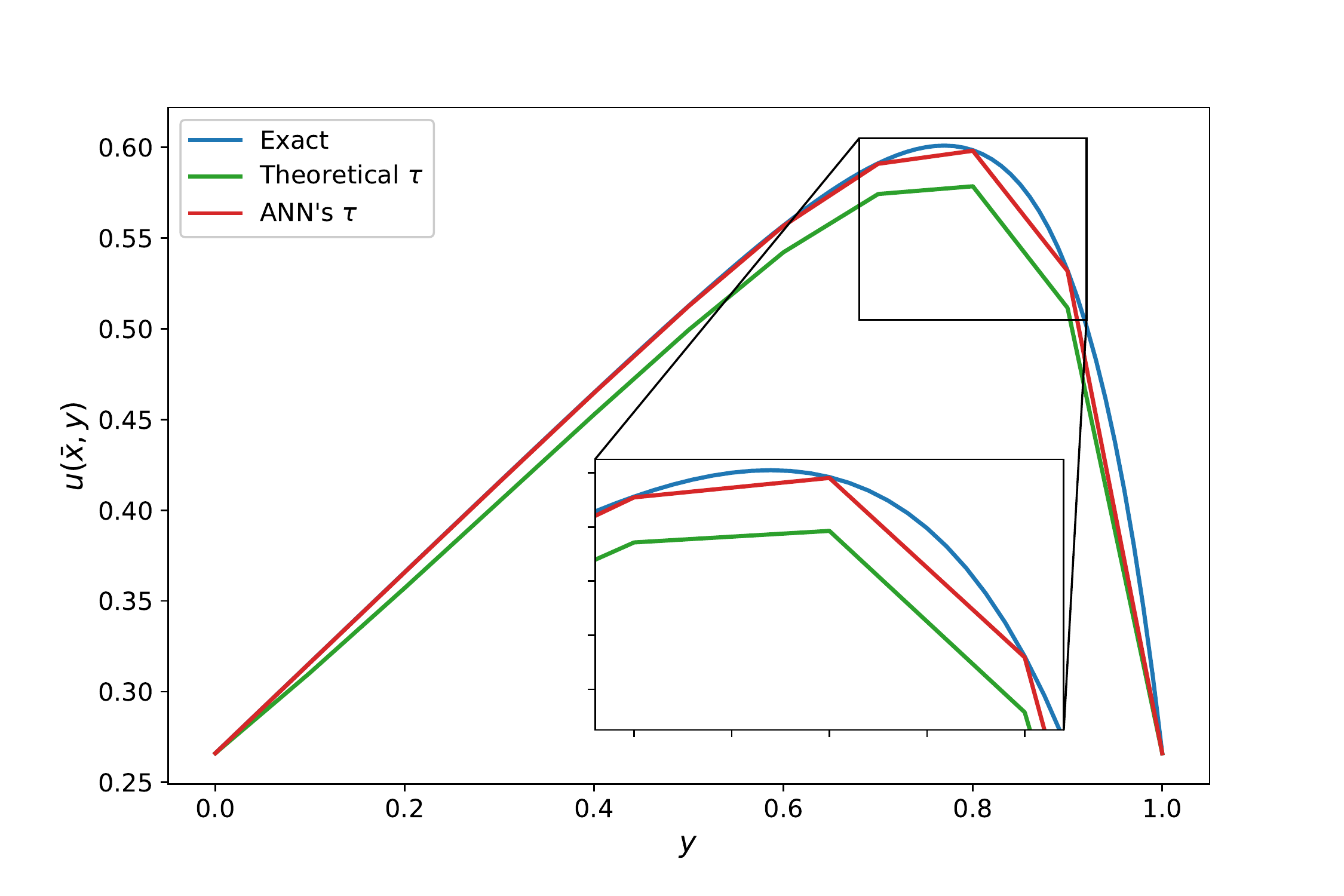}
         \caption{$\Peclet_g = 7$, $h = \sqrt{2}/10$, $r = 1$}
         \label{fig:td_2d_solutions_extracted_new_1}
    \end{subfigure}
    \begin{subfigure}{0.49\textwidth}
         \includegraphics[width=\textwidth]{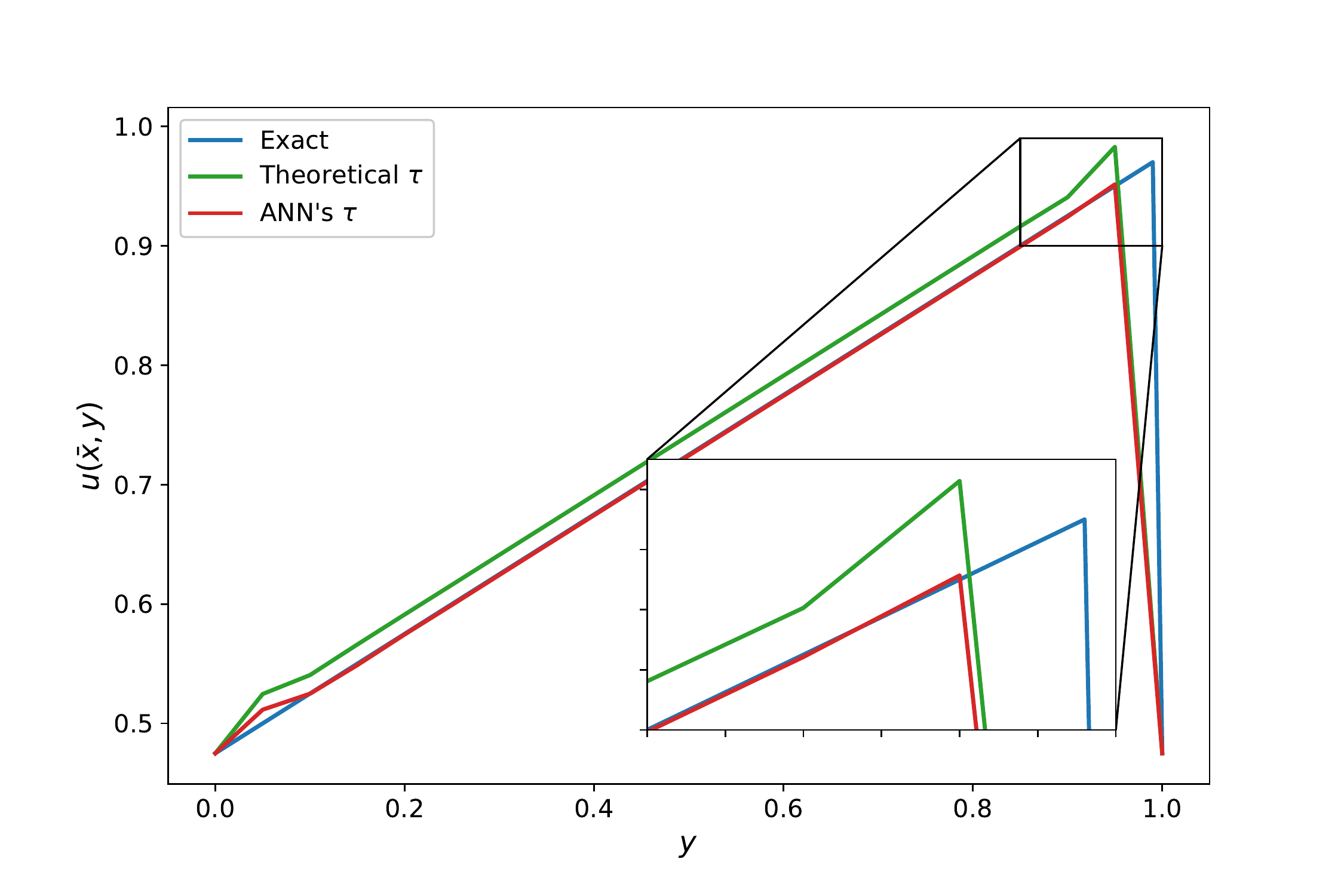}
         \caption{$\Peclet_g = 7'071$, $h = \sqrt{2}/20$, $r = 3$ }
         \label{fig:td_2d_solutions_extracted_new_2}
    \end{subfigure}
        \begin{subfigure}{0.49\textwidth}
         \includegraphics[width=\textwidth]{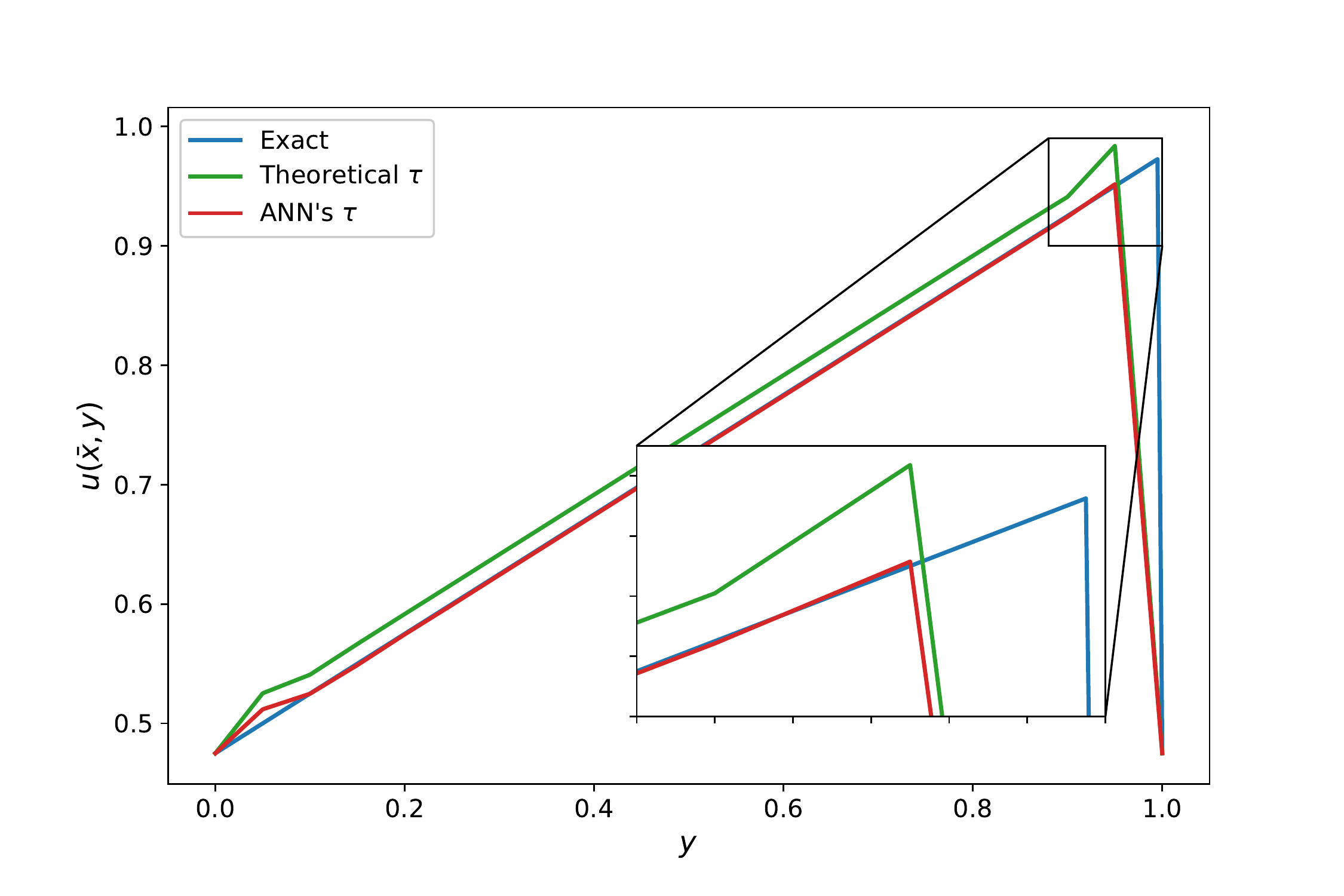}
         \caption{$\Peclet_g = 14'142$, $h = \sqrt{2}/10$, $r = 3$}
         \label{fig:td_2d_solutions_extracted_new_3}
    \end{subfigure}
    \begin{subfigure}{0.49\textwidth}
         \includegraphics[width=\textwidth]{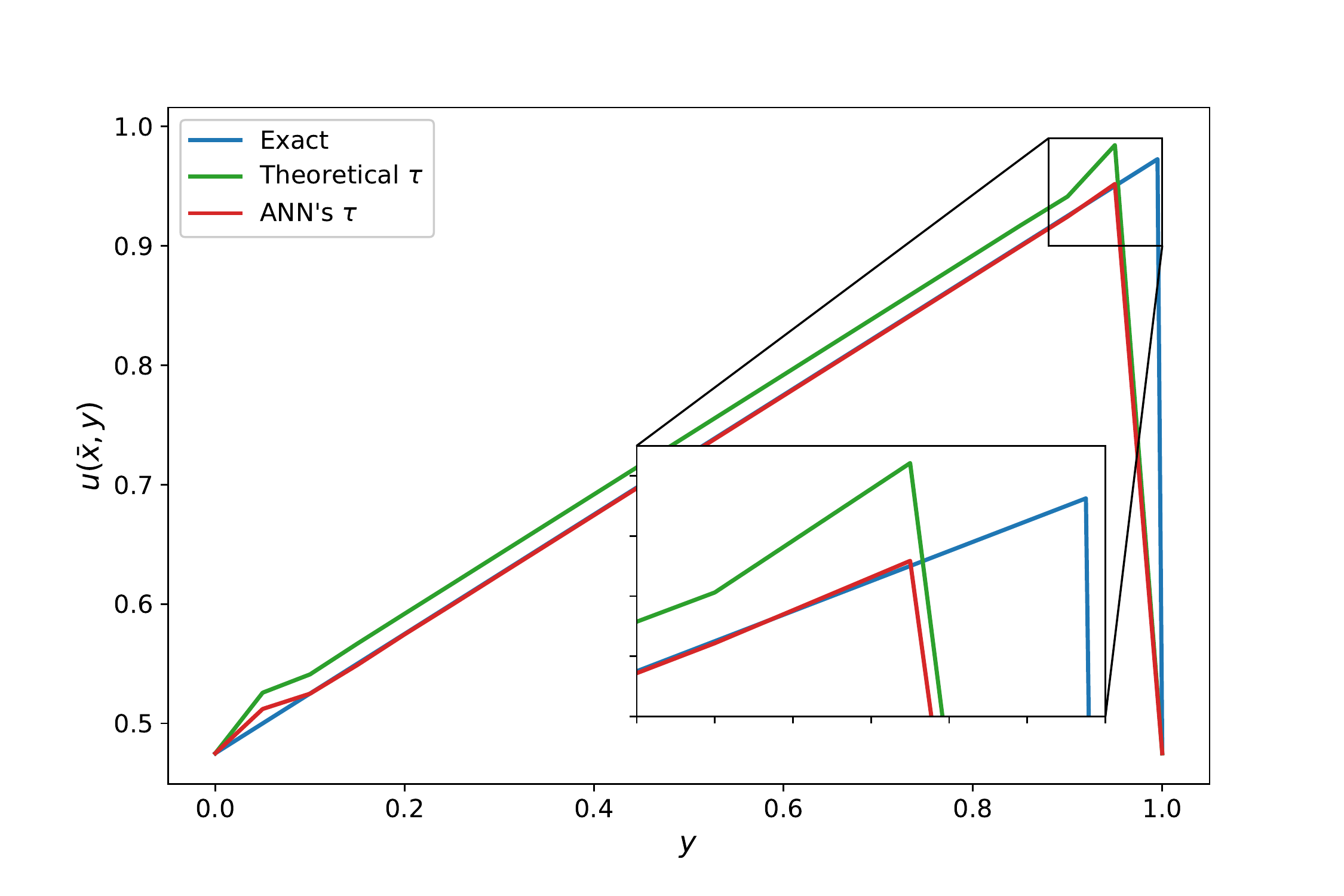}
         \caption{$\Peclet_g = 70'710$, $h = \sqrt{2}/20$, $r = 3$ }
         \label{fig:td_2d_solutions_extracted_new_4}
    \end{subfigure}
    \caption{Test 2:  Comparison of the solutions $u$~(\ref{eq:td_2d_exact_new}) (blue), $\widetilde{u}_h$ (green), and $u_h^\mathrm{ANN}$ (red) of the \emph{unseen} 2D advection-diffusion problem along the line $(1-h,y)$ for any $y \in [0,1]$.}
    \label{fig:td_2d_solutions_extracted_new}
\end{figure}

We better assess the former qualitative consideration, by computing the error $E(\tau)$ associated with the numerical solutions with SUPG stabilization for the parameters $\tau_\mathrm{ANN}$ and $\widetilde{\tau}_r$; different values of $\Peclet_g$ and $r$ are considered. In particular, Figure~\ref{fig:td_error_comparison_new} compares $E(\tau)$ against $\Peclet_g$ (in logarithmic scale) for different values of $r$: the error obtained by using $\tau_\mathrm{ANN}$ is always lower than the one achieved with $\widetilde{\tau}_r$. These results indicate that $\tau_\mathrm{ANN}$, although trained on a specific advection-diffusion problem, can be used to make predictions of the optimal $\tau$ for an unseen advection-diffusion problem in place of the theoretical stabilization parameter. 

\begin{figure}[H]
    \centering
    \begin{subfigure}{0.45\textwidth}
         \includegraphics[width=\textwidth]{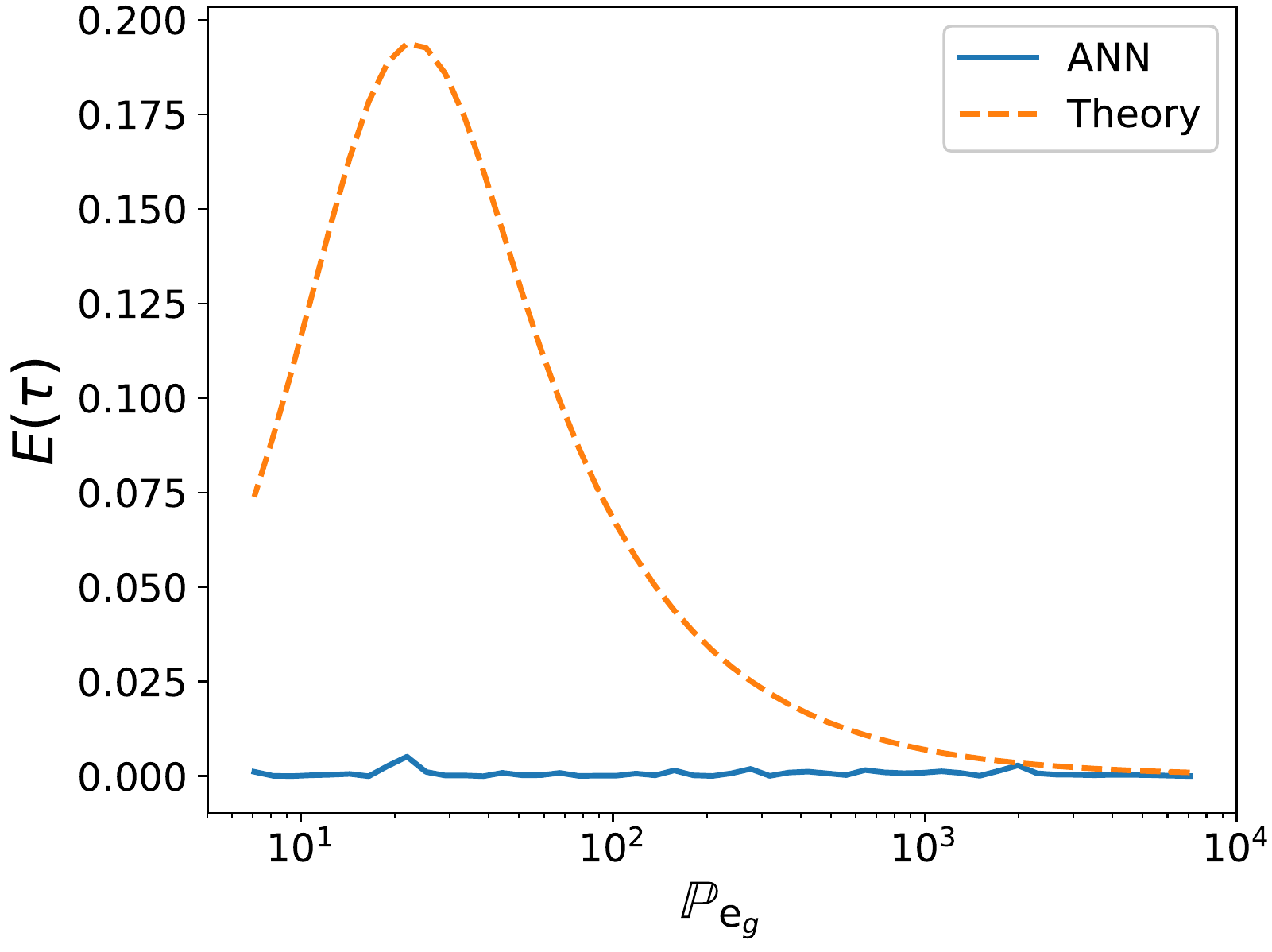}
         \caption{$r = 1$}
    \end{subfigure}
    \begin{subfigure}{0.45\textwidth}
         \includegraphics[width=\textwidth]{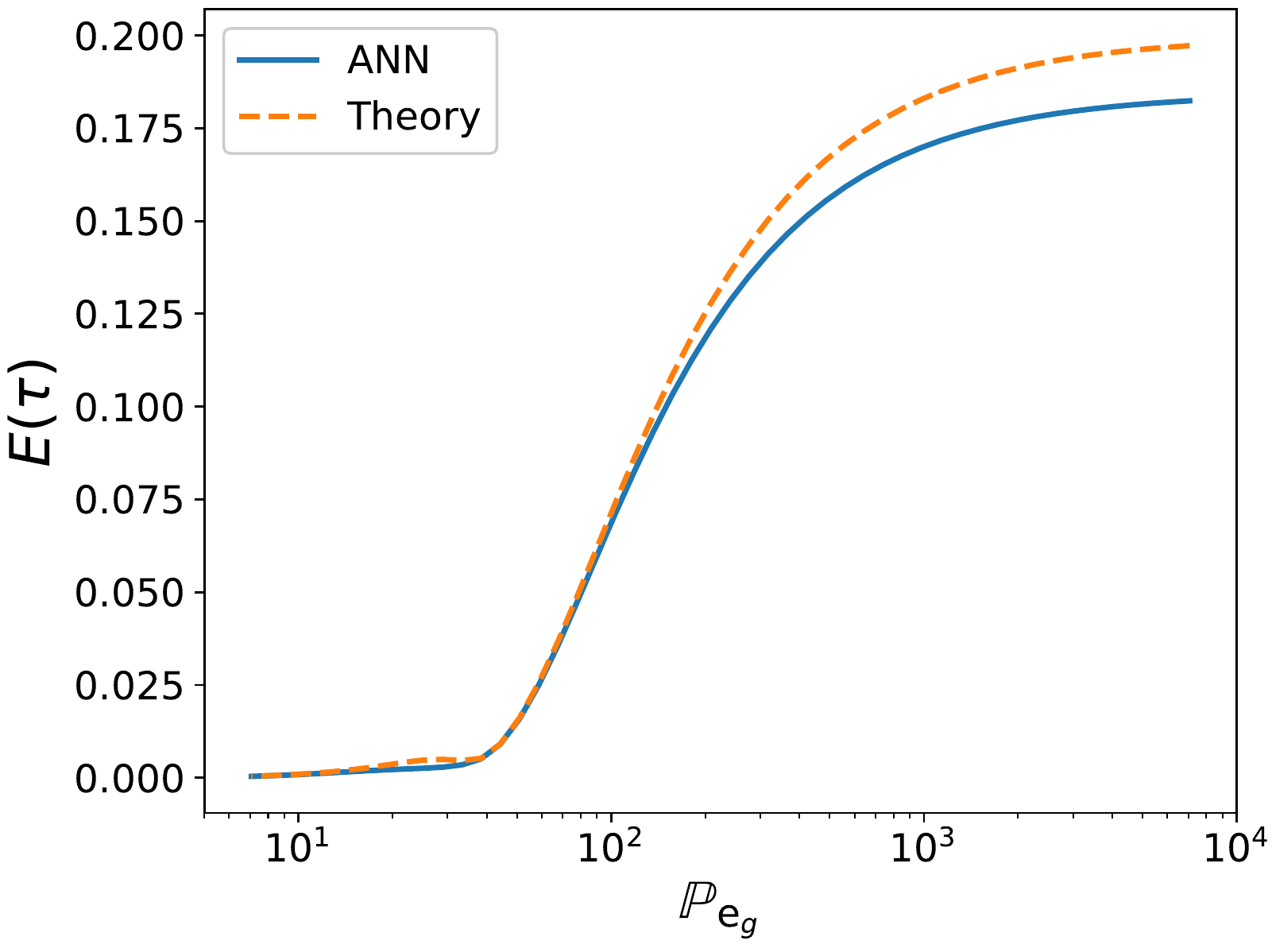}
         \caption{$r = 2$}
    \end{subfigure}
    \begin{subfigure}{0.45\textwidth}
         \includegraphics[width=\textwidth]{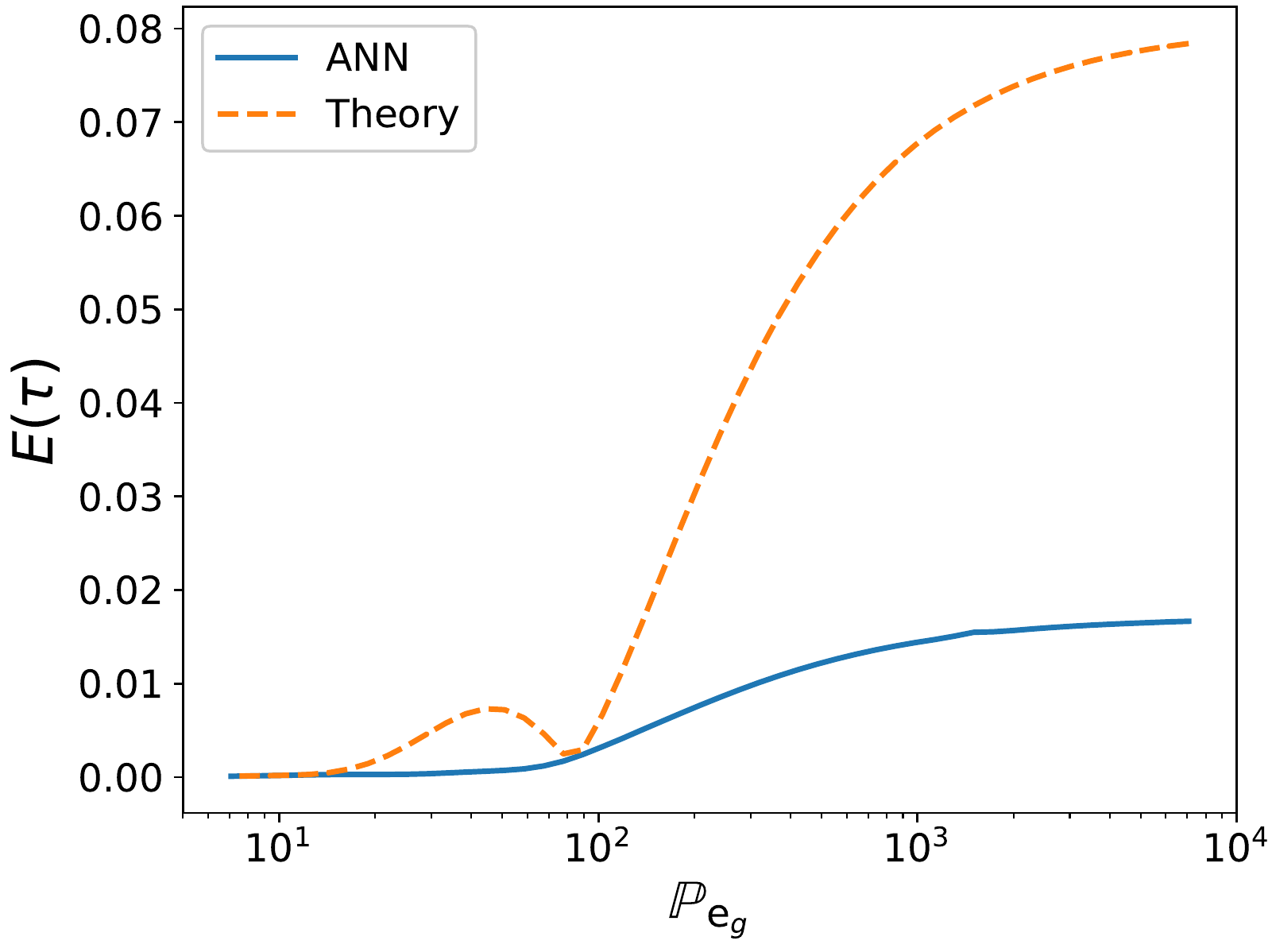}
         \caption{$r = 3$}
    \end{subfigure}
    \caption{Test 2:  Comparison of errors $E(\tau)$ in Eq~(\ref{eq:td_error_definition}) obtained with $\tau_\mathrm{ANN}$ (blue, full line) and $\widetilde{\tau}_r$ (orange, ashed line) applied to the unseen 2D advection-diffusion problem with the exact solution $u$~(\ref{eq:td_2d_exact_new}); $E(\tau)$ against $\Peclet_g$ with $h = \sqrt{2}/10$ and FE degrees $r=1,2$, and $3$.}
    \label{fig:td_error_comparison_new}
\end{figure}

\subsubsection{Test 3: predictions for and unseen problem with constant forcing term without exact solution}

In this section, we consider an advection-diffusion problem with constant forcing term $f=1$, advection coefficient $\boldsymbol{\beta} = (1, 1)$, and boundary conditions $u = 0$ on $\partial \Omega$. Although there is not an exact solution to the given problem, we consider as  ``ground truth'' (reference) solution a numerical solution obtained on a much finer grid ($h=\sqrt 2/ 400$) without stabilization.
In Figure \ref{fig:td_2d_solutions_extracted_new2}, we compare the SUPG numerical solutions $\widetilde{u}_h$, $u_h^\mathrm{ANN}$ against our ``ground truth'' solution for two different Péclet numbers. The numerical solution with the stabilization parameter $\tau_\mathrm{ANN}$ provides more accurate results with respect to the the solution obtained with the theoretical parameter. In particular, we highlight the largest discrepancies observed in the corners of the extracted solution. 
\begin{figure}[H]
    \centering
    \begin{subfigure}{0.49\textwidth}
         \includegraphics[width=\textwidth]{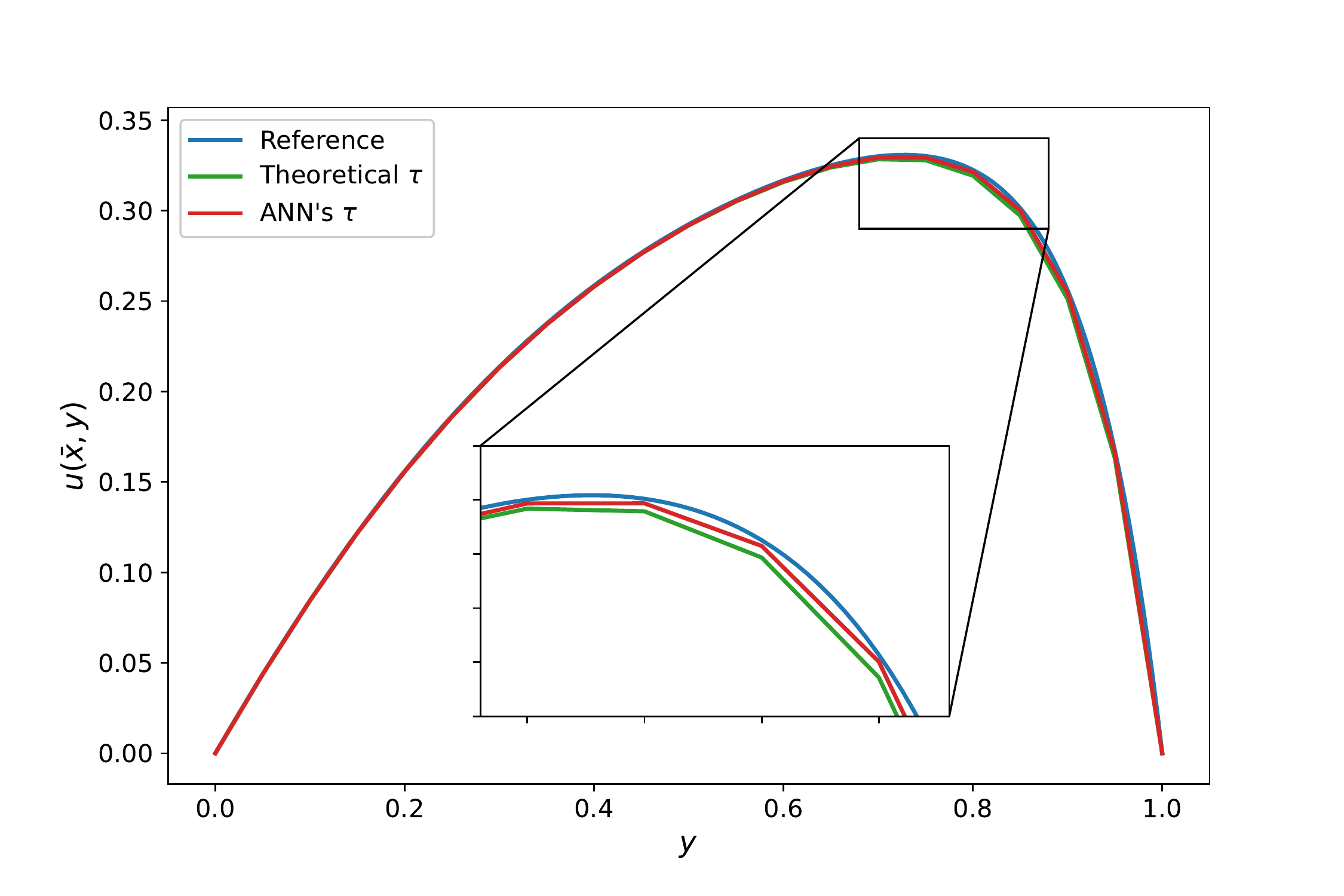}
         \caption{$\Peclet_g = 7$, $h = \sqrt{2}/20$, $r = 1$}
    \end{subfigure}
    \begin{subfigure}{0.49\textwidth}
         \includegraphics[width=\textwidth]{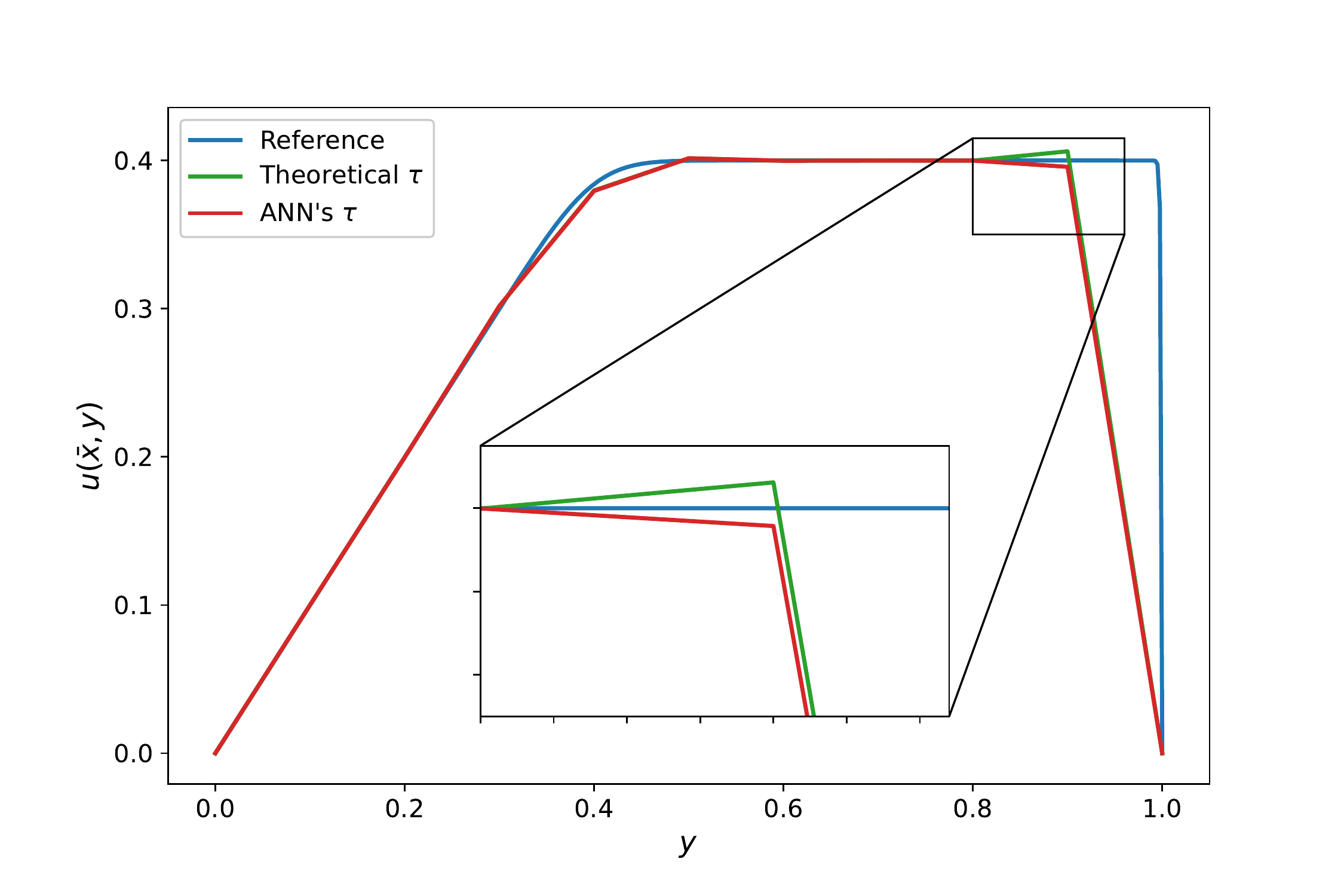}
         \caption{$\Peclet_g = 707$, $h = \sqrt{2}/10$, $r = 3$}
    \end{subfigure}
    \caption{Test 3: comparison of the solutions $u$ (blue), $\widetilde{u}_h$ (green), and $u_h^\mathrm{ANN}$ (red) of an \emph{unseen} 2D advection-diffusion problem along the line $(0.5,y)$ for any $y \in [0,1]$.}
    \label{fig:td_2d_solutions_extracted_new2}
\end{figure}

\FloatBarrier
\subsubsection{Test 4: prediction for an unseen problem with a non-constant forcing term}

In this section, we cope the case of a non costant forcing term by considering  advection-diffusion problem of Eq~\eqref{eq:transport_diffusion} in $\Omega=(0,1)^2$ with $\boldsymbol \beta = (1, 1)$, and we prescribe the following exact solution on the whole boundary $\partial \Omega$:
\begin{equation}
    \label{eq:exact_sol_noncostantf}
    u(x, y) = - \frac{\mathrm{atan}{\left( (x-1/2)^2 + (y-1/2)^2 - 1/16 \right)}}{\sqrt{\mu}}.
\end{equation}
We display the exact solution of the considered problem in Figure~\ref{fig:exact_sol_noncostantf}.
\begin{figure}[H]
    \centering
    \includegraphics[width=0.65\textwidth]{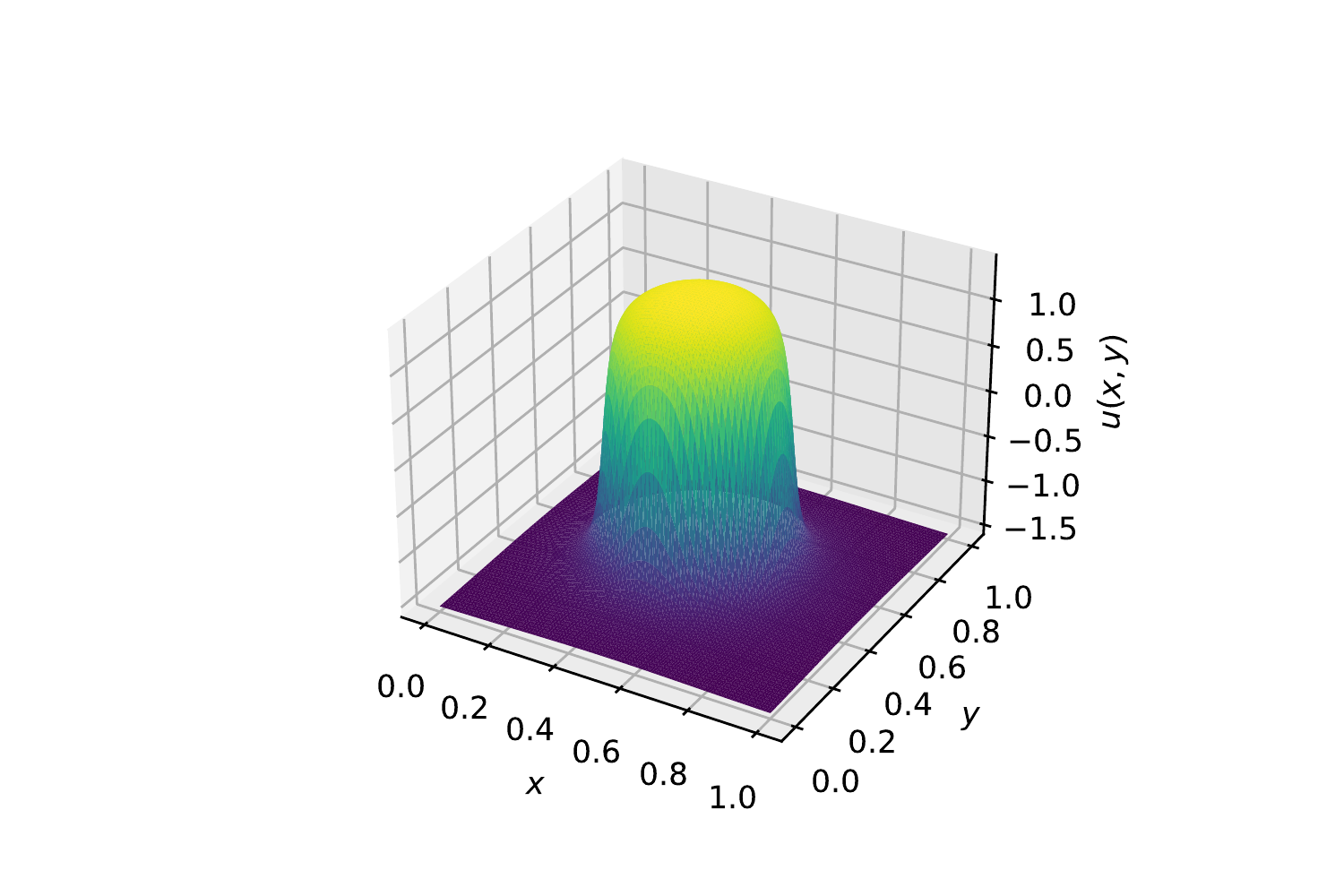}
    \caption{Test 4: exact solution in Eq~\eqref{eq:exact_sol_noncostantf} of the relative advection-diffusion problem.}
    \label{fig:exact_sol_noncostantf}
\end{figure}

In particular, Figures~\ref{fig:panettone_comparison_n10} and~\ref{fig:panettone_comparison_n20} show $E(\tau)$ against $\Peclet_g$ for different values of $r$ with two different meshes with $h = \sqrt{2}/10$ and $h = \sqrt{2}/20$, respectively. The error obtained by using $\tau_\mathrm{ANN}$ is always comparable to the one achieved with $\widetilde{\tau}_r$, for both mesh levels and for all the FE degrees considered. These results suggest that $\tau_\mathrm{ANN}$, although trained on a specific advection-diffusion problem with $f=0$, is robust with respect to unseen parameters and data in the advection-diffusion problem, including non constant $f$. Nevertheless, further studies can be conducted to better assess the role of a non-constant forcing term into the ANN's training phase, 
possibly encompassing the accuracy of the theoretical stabilization parameter in this scenario too.

\begin{figure}[H]
    \centering
    \begin{subfigure}{0.49\textwidth}
         \includegraphics[width=\textwidth]{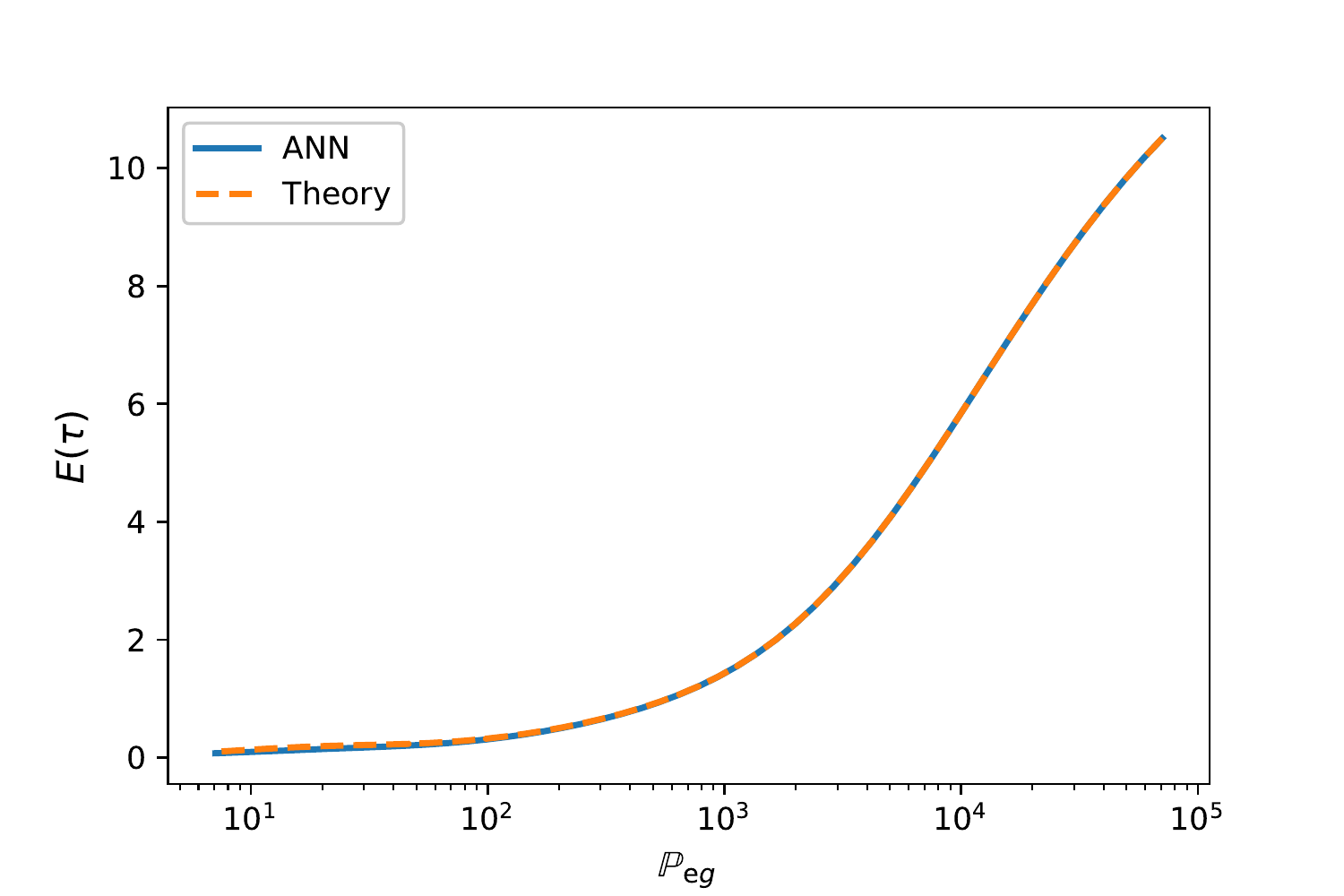}
         \caption{$r = 1$}
    \end{subfigure}
    \begin{subfigure}{0.49\textwidth}
         \includegraphics[width=\textwidth]{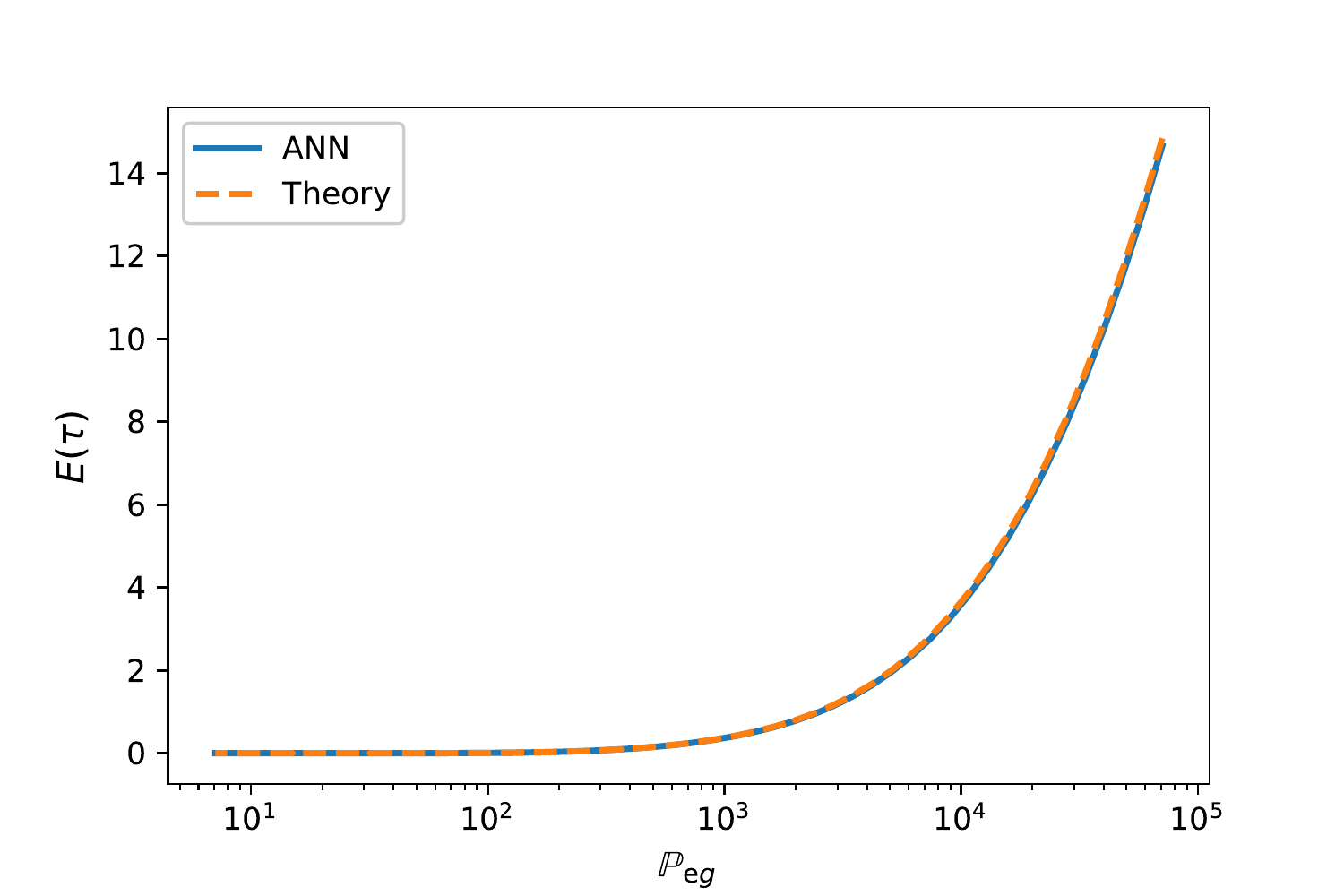}
         \caption{$r = 2$}
    \end{subfigure}
    \\
    \begin{subfigure}{0.49\textwidth}
         \includegraphics[width=\textwidth]{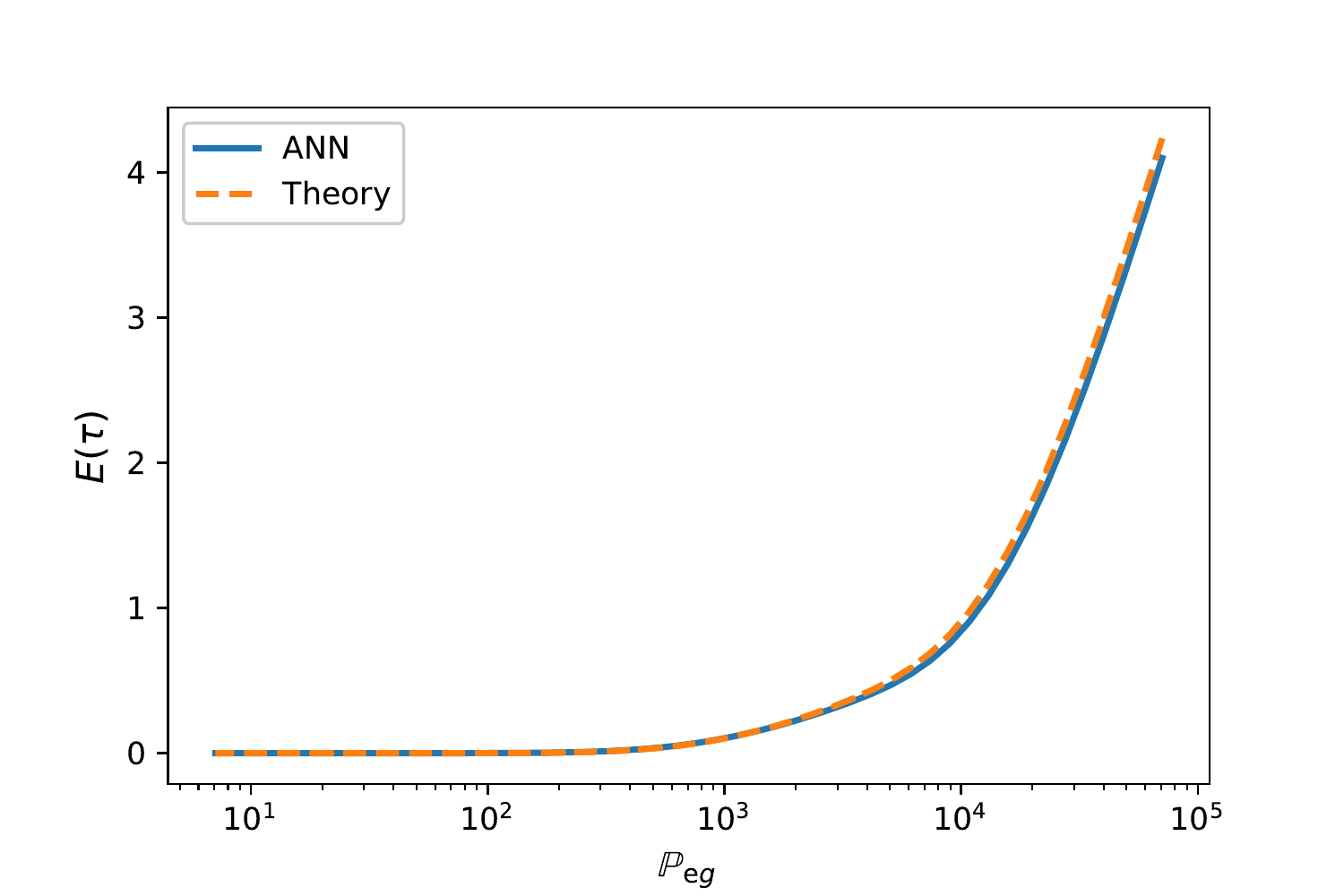}
         \caption{$r = 3$}
    \end{subfigure}
    \caption{Test 4: comparison of errors $E(\tau)$ in Eq~(\ref{eq:td_error_definition}) obtained with $\tau_\mathrm{ANN}$ (blue, full line) and $\widetilde{\tau}_r$ (orange, ashed line) applied to the unseen 2D advection-diffusion problem with the exact solution $u$ (\ref{eq:exact_sol_noncostantf}); $E(\tau)$ against $\Peclet_g$ with $h = \sqrt{2}/10$ and FE degrees $r=1,2$, and $3$. }
    \label{fig:panettone_comparison_n10}
\end{figure}
\newpage
\begin{figure}[H]
    \centering
    \begin{subfigure}{0.49\textwidth}
         \includegraphics[width=\textwidth]{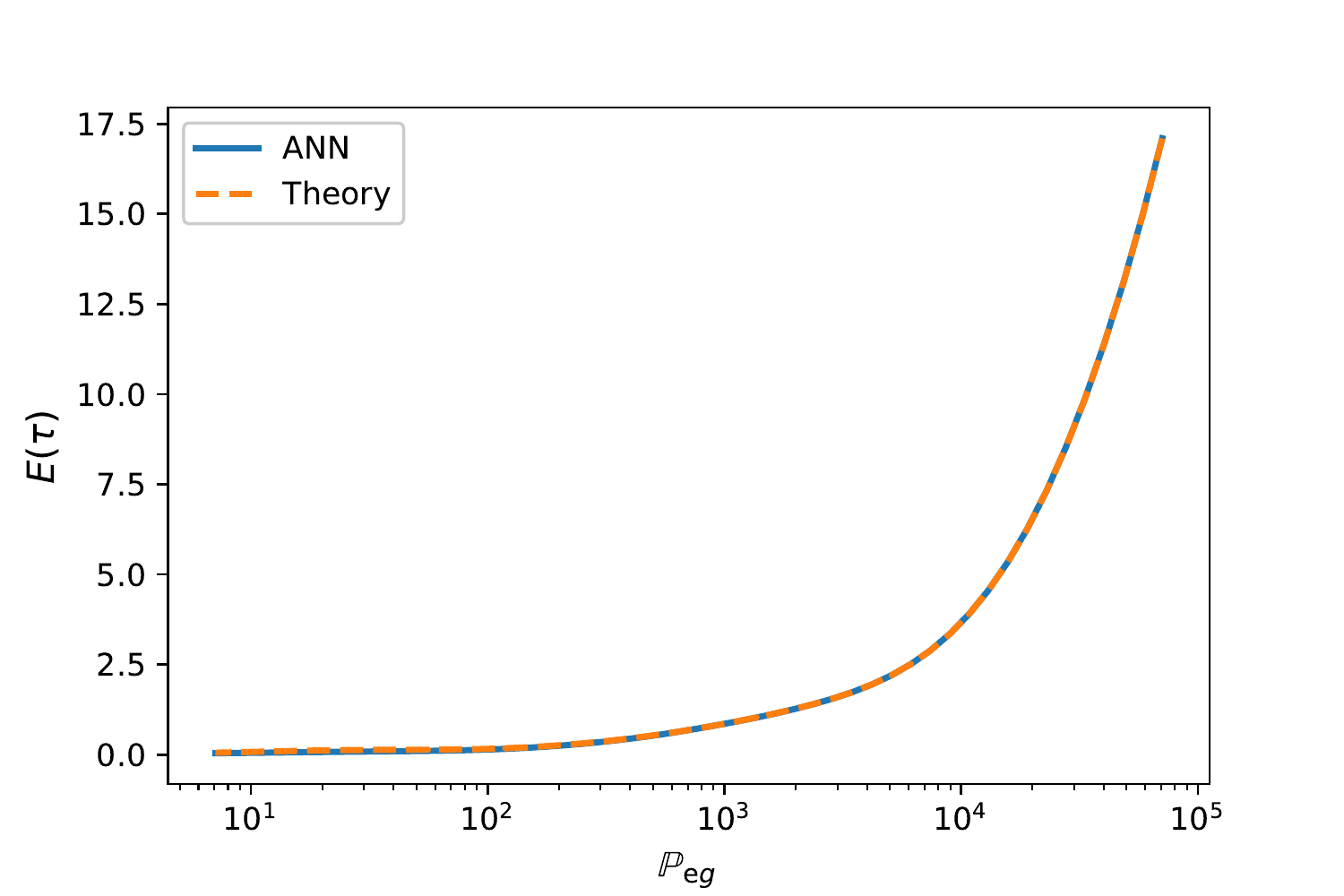}
         \caption{$r = 1$}
    \end{subfigure}
    \begin{subfigure}{0.49\textwidth}
         \includegraphics[width=\textwidth]{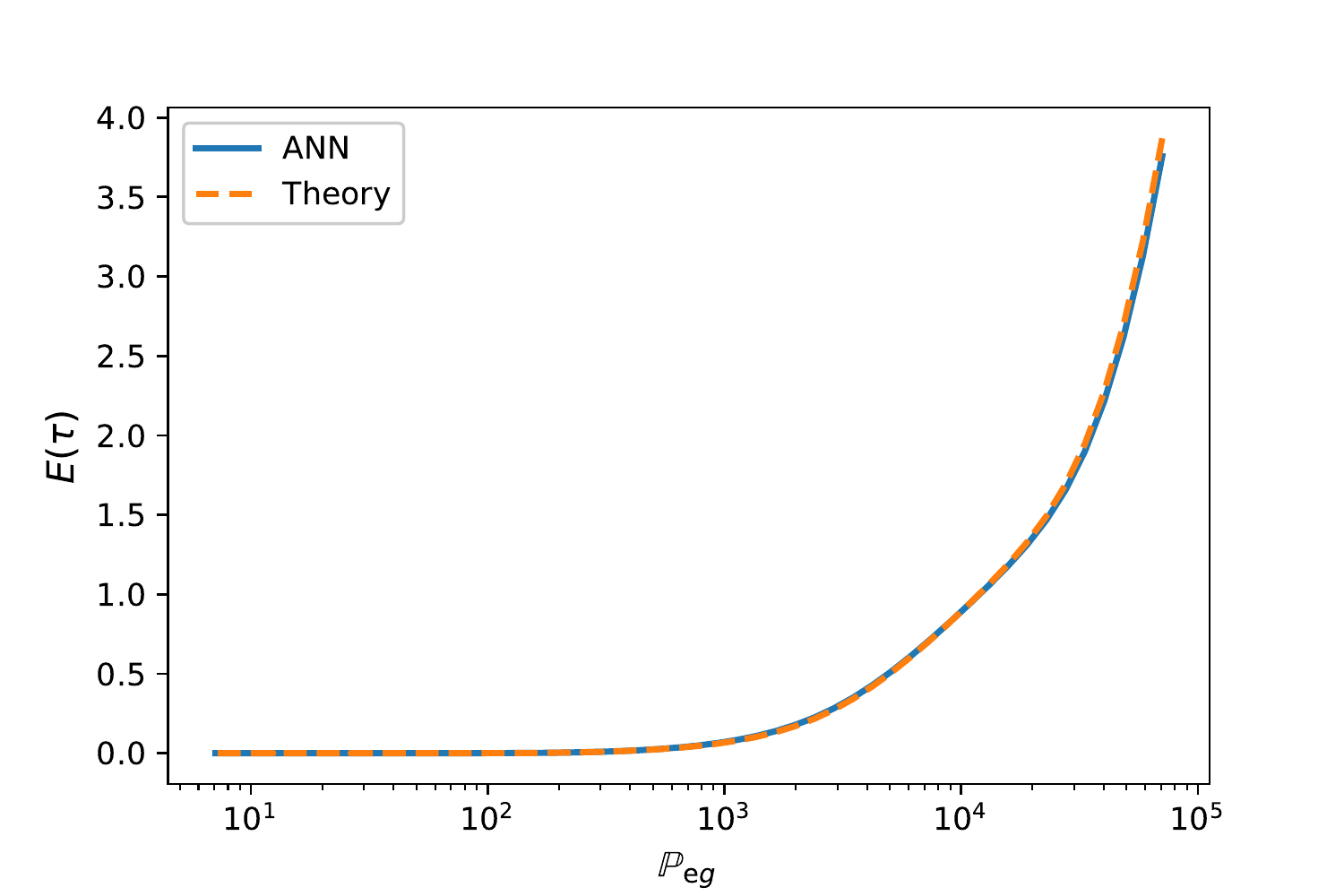}
         \caption{$r = 2$}
    \end{subfigure}
       \begin{subfigure}{0.49\textwidth}
         \includegraphics[width=\textwidth]{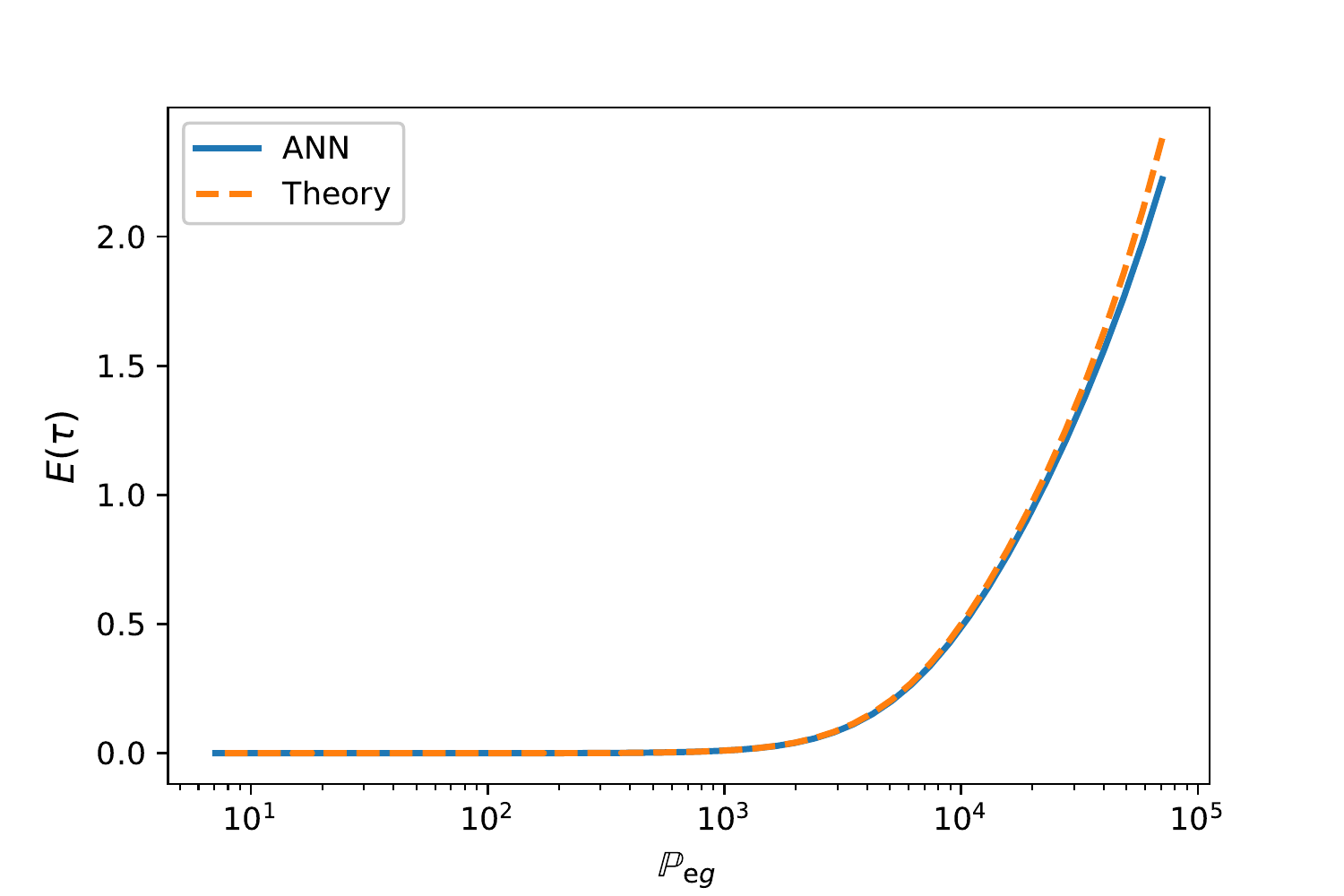}
         \caption{$r = 3$}
    \end{subfigure}
    \caption{ Test 4: comparison of errors $E(\tau)$ in Eq~(\ref{eq:td_error_definition}) obtained with $\tau_\mathrm{ANN}$ (blue, full line) and $\widetilde{\tau}_r$ (orange, ashed line) applied to the unseen 2D advection-diffusion problem with the exact solution $u$ (\ref{eq:exact_sol_noncostantf}); $E(\tau)$ against $\Peclet_g$ with $h = \sqrt{2}/20$ and FE degrees $r=1,2$, and $3$. }
    \label{fig:panettone_comparison_n20}
\end{figure}

\section{Direction of the advection field and stabilization parameters}

 We investigate the effect of the advection direction on the value of the stabilization parameter. We define $\theta$ as the angle between the advection velocity and the $x$-axis: $\boldsymbol{\beta} = (|\boldsymbol{\beta}| \cos \theta, \;|\boldsymbol{\beta}| \sin \theta )^T$. We carry out the optimization strategy introduced in Section~\ref{sec:strategy} including also $\theta$ in the input parameter $\mathbf x^{(i)}$. We apply the optimization to the advection-diffusion problem with exact solution in Eq~\eqref{eq:td_2d_exact_new} varying $\theta$ in $[\pi / 12, \pi/2]$ and we report a comparison between the optimal, theoretical and ANN's stabilization parameter in Figure~\ref{fig:beta_direction_sum2} (left). 
 We set $|\boldsymbol \beta | = \sqrt 2$, $\Peclet_g = 7'071$, $h = \sqrt{2}/20$, $r = 3$.
 We recall that, differently from the optimization strategy, in this plot we are still employing the ANN trained without accounting for the advection direction. We can observe that the optimal $\tau$ is affected by $\theta$, a feature that is not accounted by the theoretical stabilization parameter. Furthermore, in Figure~\ref{fig:beta_direction_sum2} (right), we compare the mean errors obtained with the optimal, theoretical and ANN's stabilization parameter. It can be observed that the ANN still provides an advantage with respect to the theoretical stabilization. However, we believe that by including the advection direction as an additional feature in the training phase of the ANN, a higher accuracy can be obtained that the one achieved in this paper, where a single direction ($\theta = \frac{\pi}{4}$) has been considered.

\begin{figure}[H]
    \centering
    \begin{subfigure}{0.49\textwidth}
         \includegraphics[width=\textwidth]{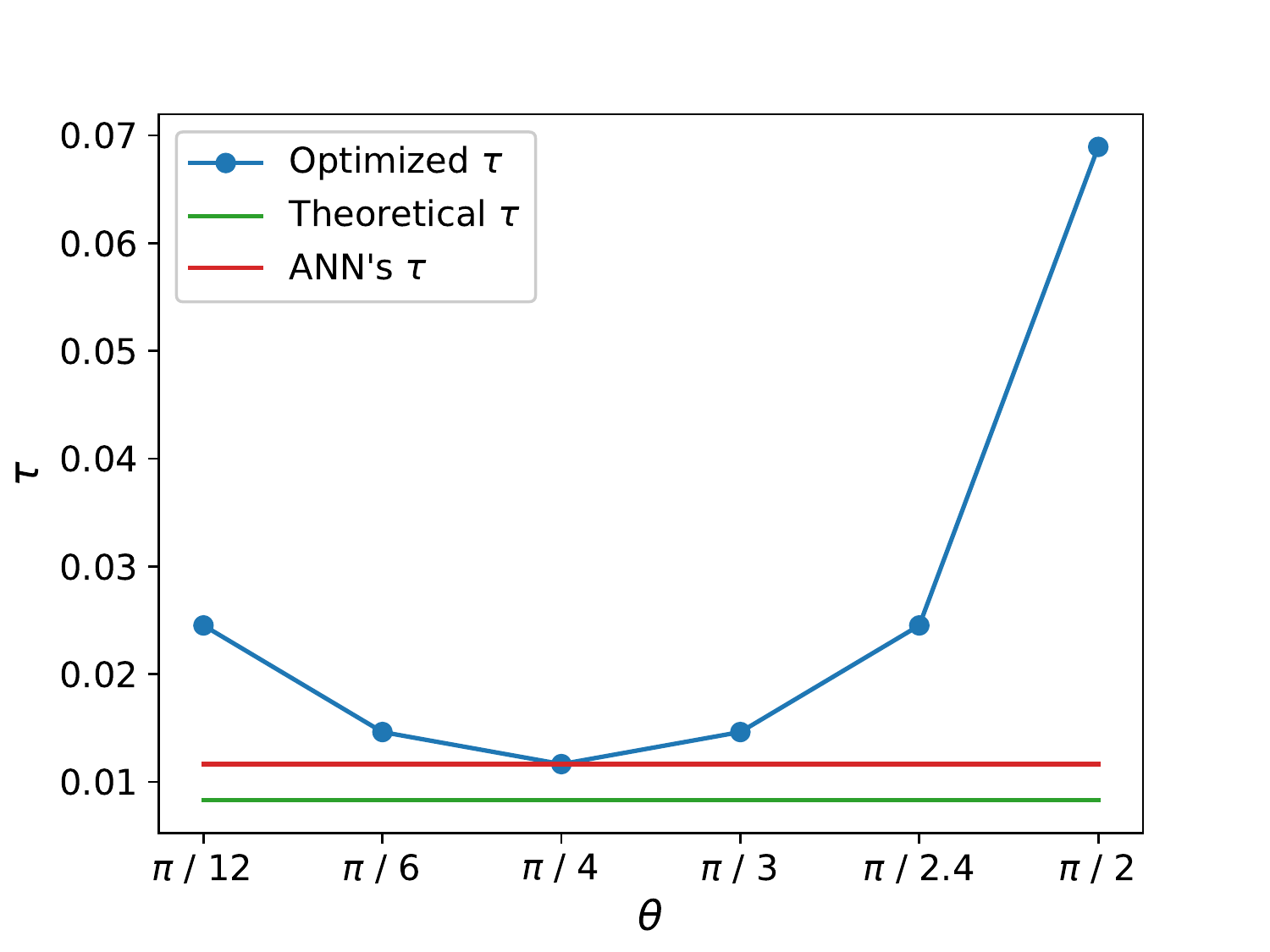}
    \end{subfigure}
    \begin{subfigure}{0.45\textwidth}
         \includegraphics[width=\textwidth]{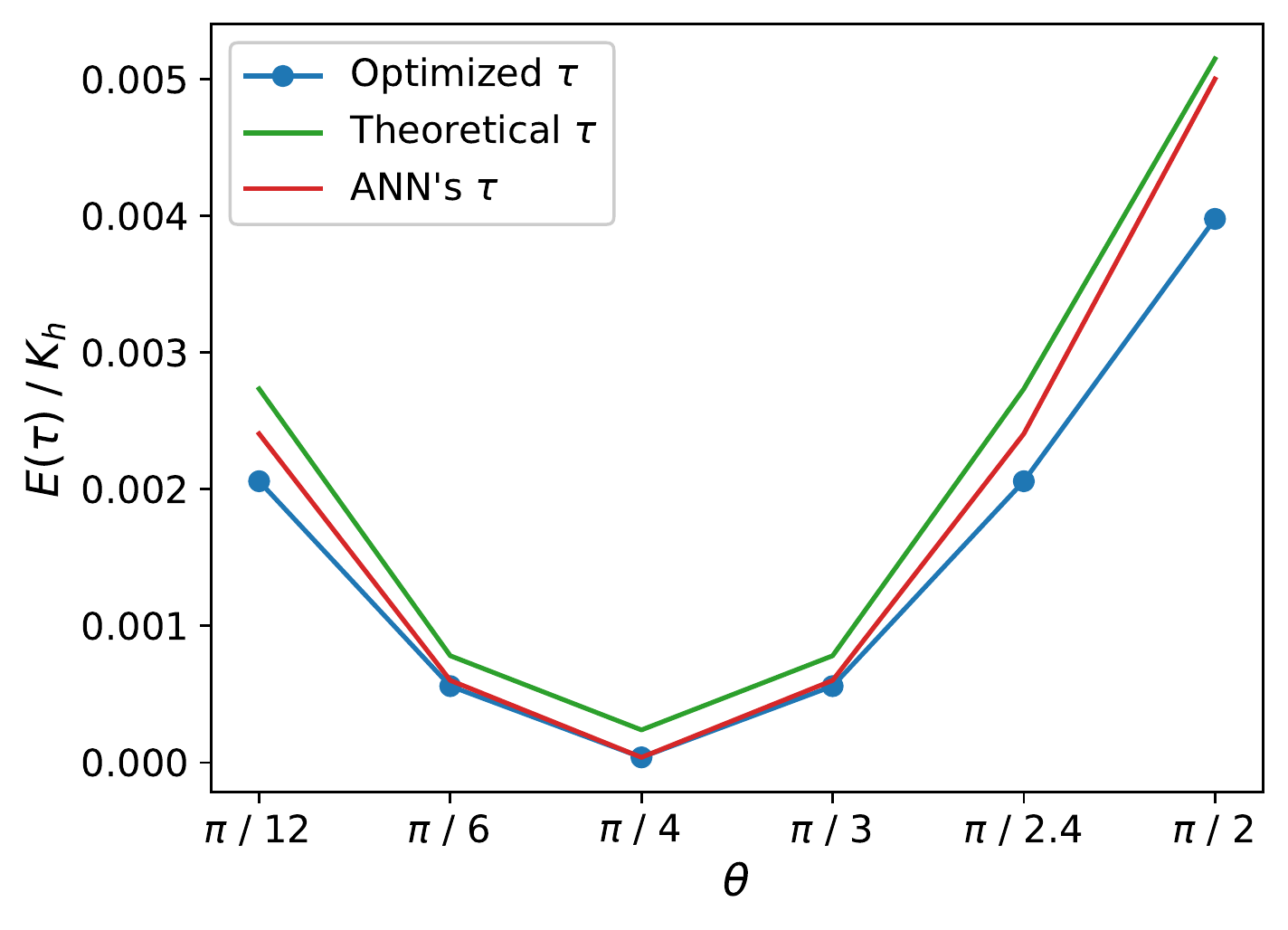}
    \end{subfigure}
    \caption{Comparison of optimal, theoretical and ANN's predicted stabilization parameters and their relative errors obtained for varying values of the advection angle $\theta$ for advection-diffusion problem with exact solution in Eq~\eqref{eq:td_2d_exact_new} with $\Peclet_g = 7'071$, $h = \sqrt{2}/20$, $r = 3$}
    \label{fig:beta_direction_sum2}
\end{figure}

\section{Conclusions}
\label{sec:conclusions}

In this work, we presented an approach based on machine learning and ANN to compute the optimal stabilization parameter to be used in the SUPG FE approximation of advection-diffusion problems. Indeed, albeit the expression of the stabilization parameter is available for 1D problems and FE of degree $r=1$, its extension to more general advection-diffusion problems (2D and 3D) and FE degrees $r>1$ is still lacking.

We validated our approach against the 1D case, for which the ANN-stabilization parameter matches the already optimal, theoretical one for $r=1$, leading to nodally exact numerical solutions Instead, for higher polynomials degrees $r>1$, remarkable differences are observed among the theoretical and optimal stabilization parameter: we observed better accuracy and stabilization properties of the numerical solution with the optimal stabilization parameter with respect to the theoretical one. 

Then, we generated the dataset on a 2D advection-diffusion problem, and we used it to train a fully-connected feed-forward ANN. We applied the predictions of the network on the original advection-diffusion problem for unseen input values and we also apply it to an unseen advection-diffusion problem to check for model generalization.
Our numerical results showed that the proposed ANN-based approach provides more accurate numerical solutions than using the theoretical stabilization parameter for the SUPG method.

This work represents a step towards the enhancements of stabilization methods for the FE approximation of advection-dominated differential problems. In particular, this work is limited to few features for the ANN-based stabilization parameter and provides as output a single optimal stabilization parameter meant for the whole mesh FE $\mathcal T_h$. Natural extensions of this work will therefore involve local stabilization parameters, i.e., element by element over the mesh $\mathcal T_h$, to better account for non-uniform meshes, differential problems with varying coefficients, and capturing the local behaviour of the solution. 
 In addition, to enhance the robustness of the ANN-based approach, possible additional inputs of the network are the angle of the advection velocity and the nodal value of the forcing term. A future development consists in the extension of the proposed approach to the 3D case for which, as for the 2D case, a universal definition of the stabilization parameter does not exist. The extension to the 3D case surely would present a larger cost in the training phase; however, once trained, the ANN could be used online real-time. 

Finally, we envision extending the proposed machine learning approach to SUPG and variational multiscale stabilization methods for the FE approximation of the Navier-Stokes equations. 




\section*{Acknowledgments}
The authors acknowledge Prof. C. Canuto, Politecnico di Torino, for the fruitful discussions and useful suggestions.  L. D. and A. Z. are members of the INdAM Research group GNCS. L. D. has been partially funded by the research project PRIN 2020, n.20204LN5N5, by MIUR.  The authors ackwnowledge the anonymous reviewers for their insightful comments and suggestions .

\section*{Conflict of interest}
The authors  declare no conflict of interest.


\begin{thebibliography}{999}
\bibitem{tensorflow}
M.~Abadi, A.~Agarwal, P.~Barham, E. Brevdo, Z. Chen, C. Citro. et al.,
 {TensorFlow}: Large-scale machine learning on heterogeneous systems, 	arXiv:1603.04467. 

\bibitem{fenics}
M. S. Aln{\ae}s, J.~Blechta, J.~Hake, A. Johansson, B. Kehlet, A. Logg, et al.,
 The fenics project version 1.5,
{\em Archive of Numerical Software}, {\bf3} (2015), 9--23. \url{//doi.org/10.11588/ans.2015.100.20553}

\bibitem{bazilevs2007variational}
Y.~Bazilevs, V. M.~Calo, J. A.~Cottrell, T. J. R. Hughes, A. Reali, G. Scovazzi,
 Variational multiscale residual-based turbulence modeling for large
  eddy simulation of incompressible flows,
 {\em Comput. Method. Appl. Mech. Eng.},
  \textbf{197} (2007), 173--201. \url{//doi.org/10.1016/j.cma.2007.07.016}

\bibitem{annurev_control}
N.~B{\'e}nard, J.~Pons-Prats, J.~P{\'e}riaux, G. Bugeda, P. Braud, J. P. Bonnet, et al.,
 Turbulent separated shear flow control by surface plasma actuator:
  experimental optimization by genetic algorithm approach,
 {\em Exp. Fluids}, \textbf{57} (2016), 22. \url{//doi.org/10.1007/s00348-015-2107-3}

\bibitem{Bochev}
P.~B. Bochev, C.~R. Dohrmann, M.~D. Gunzburger,
 Stabilization of low-order mixed finite elements for the {S}tokes
  equations,
{\em SIAM J. Numer. Anal.}, \textbf{44} (2016), 82--101. \url{//doi.org/10.1137/S0036142905444482}

\bibitem{Brooks_1982}
A.~N. Brooks, T.~J.~R. Hughes,
 Streamline upwind/{P}etrov-{G}alerkin formulations for convection
  dominated flows with particular emphasis on the incompressible
  {N}avier-{S}tokes equations,
 {\em Comput. Method. Appl. Mech. Eng.},
 \textbf{32} (1982), 199--259.
 \url{//doi.org/10.1016/0045-7825(82)90071-8}
  

\bibitem{spectralbook}
C.~Canuto, M. Y. Hussaini, A.~Quarteroni, T. A. Zang,
  {\em Spectral methods. {F}undamentals in single domains}, Berlin, Heidelberg:
  Springer, 2006.

\bibitem{CODINA200061}
R. Codina,
  On stabilized finite element methods for linear systems of
  convection-diffusion-reaction equations,
  {\em Comput. Method. Appl. Mech. Eng.},
  \textbf{188} (2000), 61--82. \url{//doi.org/10.1016/S0045-7825(00)00177-8}

\bibitem{CODINA2008264}
R. Codina,
  Analysis of a stabilized finite element approximation of the {O}seen
  equations using orthogonal subscales,
  {\em Appl. Numer. Math.}, \textbf{58} (2008), 264--283. \url{//doi.org/10.1016/j.apnum.2006.11.011}

\bibitem{CODINA20072413}
R. Codina, J. Principe, O. Guasch, S. Badia,
  Time dependent subscales in the stabilized finite element
  approximation of incompressible flow problems,
  {\em Comput. Method. Appl. Mech. Eng.},
  \textbf{196} (2007), 2413--2430. \url{//doi.org/10.1016/j.cma.2007.01.002}

\bibitem{annurev_wake}
B.~Colvert, M.~Alsalman, E.~Kanso,
  Classifying vortex wakes using neural networks,
  {\em Bioinspir. Biomim.}, \textbf{13} (2018), 025003. \url{//doi.org/10.1088/1748-3190/aaa787}

\bibitem{cottrell2009isogeometric}
J. A. Cottrell, T. J. R. Hughes, Y. Bazilevs,
  {\em Isogeometric analysis: {T}oward integration of {CAD} and
  {FEA}},   John Wiley \& Sons, 2009. \url{//doi.org/10.1002/9780470749081}

\bibitem{discacciati2020controlling}
N. Discacciati, J. S. Hesthaven,  D. Ray,
  Controlling oscillations in high-order discontinuous {G}alerkin
  schemes using artificial viscosity tuned by neural networks,
  {\em J. Comput. Phys.}, \textbf{409} (2020), 109304. \url{//doi.org/10.1016/j.jcp.2020.109304}

\bibitem{annurev_rans}
K. Duraisamy, G. Iaccarino,  H. Xiao,
  Turbulence modeling in the age of data,
  {\em Annu. Rev. Fluid Mech.}, \textbf{51} (2019), 357--377. \url{//doi.org/10.1146/annurev-fluid-010518-040547}

\bibitem{forti2015semi}
D. Forti, L. Dede',
  Semi-implicit {B}{D}{F} time discretization of the {N}avier--{S}tokes
  equations with {V}{M}{S}-{L}{E}{S} modeling in a high performance computing
  framework,
  {\em Comput. Fluids}, \textbf{117} (2015), 168--182. \url{//doi.org/10.1016/j.compfluid.2015.05.011}

\bibitem{franca1992stabilized}
L. P. Franca, S. L. Frey,  T. J. R. Hughes,
  Stabilized finite element methods: {I}. application to the
  advective-diffusive model,
  {\em Comput. Method. Appl. Mech. Eng.},  \textbf{95} (1992), 253--276. \url{//doi.org/10.1016/0045-7825(92)90143-8}

\bibitem{fresca2021comprehensive}
S. Fresca, L. Dede', A. Manzoni,
  A comprehensive deep learning-based approach to reduced order
  modeling of nonlinear time-dependent parametrized {PDEs},
  {\em J. Sci. Comput.}, \textbf{87} (2021), 61. \url{//doi.org/10.1007/s10915-021-01462-7}

\bibitem{galeao2004finite}
A. C. Galeao, R. C. Almeida, S. M. C. Malta, A. F. D. Loula,
  Finite element analysis of convection dominated reaction--diffusion
  problems,
  {\em Appl. Numer. Math.}, \textbf{48} (2004), 205--222. \url{//doi.org/10.1016/j.apnum.2003.10.002}

\bibitem{goodfellow2016deep}
I.~Goodfellow, Y.~Bengio, A.~Courville,
  {\em Deep learning},  MIT Press, 2016.

\bibitem{guo2019data}
M. Guo, J. S. Hesthaven,
  Data-driven reduced order modeling for time-dependent problems,
  {\em Comput. Method. Appl. Mech. Eng.},
  \textbf{345} (2019), 75--99. \url{//doi.org/10.1016/j.cma.2018.10.029}

\bibitem{hesthaven2018non}
J. S. Hesthaven, S. Ubbiali,
  Non-intrusive reduced order modeling of nonlinear problems using
  neural networks,
  {\em J. Comput. Phys.}, \textbf{363} (2018), 55--78. \url{//doi.org/10.1016/j.jcp.2018.02.037}

\bibitem{hughes2012finite}
T. J. R. Hughes,
  {\em The finite element method: linear static and dynamic finite element analysis},
  Courier Corporation, 2012. 

\bibitem{janssensadvancing}
M. Janssens, S. J. Hulshoff,
  Advancing artificial neural network parameterisation for atmospheric
  turbulence using a variational multiscale model,
  {\em J. Adv. Model. Earth Syst.}, \textbf{14} (2022),
  e2021MS002490. \url{//doi.org/10.1029/2021MS002490}

\bibitem{janssenthesis}
M. Janssens,  Machine learning of atmospheric turbulence in a variational multiscale model, 2019. Available from: \url{http://resolver.tudelft.nl/uuid:bd090309-305e-4c04-93b7-64f1b79df8d4}

\bibitem{john2007spurious}
V. John, P. Knobloch,
  On spurious oscillations at layers diminishing (sold) methods for
  convection--diffusion equations: Part I--A review,
  {\em Comput. Method. Appl. Mech. Eng.},
  \textbf{196} (2007), 2197--2215. \url{//doi.org/10.1016/j.cma.2006.11.013}

\bibitem{keras}
Keras.   Available from:   \url{//keras.io}.

\bibitem{Kutyniok_2021}
G. Kutyniok, P. Petersen, M. Raslan, R. Schneider,
  A theoretical analysis of deep neural networks and parametric {PDEs},
  {\em Constr. Approx.}, \textbf{55} (2021), 73--125. \url{//doi.org/10.1007/s00365-021-09551-4}


\bibitem{annurev_turbulent_wall}
M. Milano, P. Koumoutsakos,
  Neural network modeling for near wall turbulent flow,
  {\em J. Comput. Phys.}, \textbf{182} (2002), 1--26. \url{//doi.org/10.1006/jcph.2002.7146}

\bibitem{Mishra}
S. Mishra,
A machine learning framework for data driven acceleration of computations of differential equations,
{\em Mathematics in Engineering}, \textbf{1} (2019), 118--146. \url{//doi.org/10.3934/Mine.2018.1.118}

\bibitem{mitchell1997machine}
T. M. Mitchell,
  {\em Machine learning}, New York:  McGraw-hill, 1997.

\bibitem{Neittaanmaki202227}
P. Neittaanmaki, S. Repin,
  Artificial intelligence and computational science, In:
  {\em Computational sciences and artificial intelligence in industry}, \textbf{76} (2022), 27--35. \url{//doi.org/10.1007/978-3-030-70787-3\_3}

\bibitem{annurev_gliding}
G. Novati, L. Mahadevan, P. Koumoutsakos,
  Controlled gliding and perching through deep-reinforcement-learning,
  {\em Phys. Rev. Fluids}, \textbf{4} (2019), 093902. \url{//doi.org/10.1103/PhysRevFluids.4.093902}

\bibitem{QV94}
A. Quarteroni, A. Valli,
  {\em Numerical approximation of partial differential equations},
  Berlin, Heidelberg: Springer, 1994. \url{//doi.org/10.1007/978-3-540-85268-1}

\bibitem{Quarteroni_2017}
A. Quarteroni,
  {\em Numerical models for differential problems}, 3 Eds., Cham:
  Springer, 2017. \url{//doi.org/10.1007/978-3-319-49316-9}

\bibitem{raissi2018hidden}
M.~Raissi, G.~E. Karniadakis,
  Hidden physics models: Machine learning of nonlinear partial
  differential equations,
  {\em J. Comput. Phys.}, \textbf{357} (2018), 125--141. \url{//doi.org/10.1016/j.jcp.2017.11.039}

\bibitem{raissi2017machine}
M.~Raissi, P.~Perdikaris, G.~E. Karniadakis,
  Machine learning of linear differential equations using {G}aussian
  processes,
  {\em J. Comput. Phys.}, \textbf{348} (2017), 683--693. \url{//doi.org/10.1016/j.jcp.2017.07.050}

\bibitem{raissi2019physics}
M. Raissi, P. Perdikaris, G.~E. Karniadakis,
  Physics-informed neural networks: A deep learning framework for
  solving forward and inverse problems involving nonlinear partial differential
  equations,
  {\em J. Comput. Phys.}, \textbf{378} (2019), 686--707. \url{//doi.org/10.1016/j.jcp.2018.10.045}

\bibitem{REBOLLO2015406}
T. C. Rebollo, B. M. Dia,
  A variational multi-scale method with spectral approximation of the sub-scales: Application to the 1D advection–diffusion equations,
  {\em Comput. Method. Appl. Mech. Eng.},
  \textbf{285} (2015), 406--426. \url{//doi.org/10.1016/j.cma.2014.11.025}

\bibitem{REGAZZONI2020113268}
F.~Regazzoni, L.~Dede',  A.~Quarteroni,
  Machine learning of multiscale active force generation models for the
  efficient simulation of cardiac electromechanics,
  {\em Comput. Method. Appl. Mech. Eng.},
  \textbf{370} (2020),  113268. \url{//doi.org/10.1016/j.cma.2020.113268}

\bibitem{regazzoni2019machine}
F. Regazzoni, L. Dede',  A. Quarteroni,
  Machine learning for fast and reliable solution of time-dependent
  differential equations,
  {\em J. Comput. Phys.}, \textbf{397} (2019), 108852. \url{//doi.org/10.1016/j.jcp.2019.07.050}

\bibitem{Roos_1996}
H. G. Roos, M. Stynes, L. Tobiska,
  {\em Numerical methods for singularly perturbed differential
  equations}, Berlin, Heidelberg:  Springer,  1996. \url{//doi.org/10.1007/978-3-662-03206-0}

\bibitem{Schwab}
C. Schwab,
  {\em {p}- and hp-Finite element methods: {T}heory and application to solid and fluid
  mechanics},  Oxford University Press, 1998.

\bibitem{SCOVAZZI2007923}
G.~Scovazzi, M. A. Christon, T. J. R. Hughes,  J. N. Shadid,
  Stabilized shock hydrodynamics: {I}. {A} {L}agrangian method,
  {\em Comput. Method. Appl. Mech. Eng.},
  \textbf{196} (2007), 923--966. \url{//doi.org/10.1016/j.cma.2006.08.008}

\bibitem{ScovazziRossi}
G. Scovazzi, B. Carnes, X. Zeng, S. Rossi,
  A simple, stable, and accurate linear tetrahedral finite element for
  transient, nearly, and fully incompressible solid dynamics: a dynamic
  variational multiscale approach,
  {\em Int. J. Numer. Meth. Eng.},
  \textbf{106} (2016), 799--839. { https://doi.org/10.1002/nme.5138}

\bibitem{tezduyar2000finite}
T. E. Tezduyar, Y. Osawa,
  Finite element stabilization parameters computed from element
  matrices and vectors,
  {\em Comput. Method. Appl. Mech. Eng.},
  \textbf{190} (2000), 411--430. \url{//doi.org/10.1016/S0045-7825(00)00211-5}

\bibitem{university1988continuous}
University of~Illinois at Urbana-Champaign. Center~for Supercomputing~Research,
Development, and G~Cybenko, Continuous valued neural networks with two hidden layers are  sufficient,   1988. Available from: \url\url{//searchworks.stanford.edu/view/4620277}.

\bibitem{scipy}
P.~Virtanen, R.~Gommers, T. E. Oliphant, M. Haberland, T. Reddy, D. Cournapeau, et al.,
  {{SciPy} 1.0: fundamental algorithms for scientific computing in
  Python},
  {\em Nat. Methods}, \textbf{17} (2020), 261--272. \url{//doi.org/10.1038/s41592-019-0686-2}

\bibitem{xie2019artificial}
C. Xie, J. Wang, K. Li, C. Ma, 
  Artificial neural network approach to large-eddy simulation of
  compressible isotropic turbulence,
  {\em Phys. Rev. E}, \textbf{99} (2019), 053113. \url{//doi.org/10.1103/PhysRevE.99.053113}

\bibitem{xie2020modeling}
C. Xie, J. Wang, W. E,
  Modeling subgrid-scale forces by spatial artificial neural networks
  in large eddy simulation of turbulence,
  {\em Phys. Rev. Fluids}, \textbf{5} (2020), 054606. \url{//doi.org/10.1103/PhysRevFluids.5.054606}


\bibitem{Yegnanarayana}
B. Yegnanarayana,  {\em Artificial neural networks},
  PHI Learning Pvt. Ltd., 2009.

\bibitem{zancanaro2021hybrid}
M. Zancanaro, M. Mrosek, G. Stabile, C. Othmer,   G. Rozza,
  Hybrid neural network reduced order modelling for turbulent flows
  with geometric parameters,
    {\em Fluids}, \textbf{6} (2021), 296.
\url{//doi.org/10.3390/fluids6080296}

\bibitem{zhou2019subgrid}
Z. Zhou, G. He, S. Wang, G. Jin,
  Subgrid-scale model for large-eddy simulation of isotropic turbulent
  flows using an artificial neural network,
  {\em Comput. Fluids}, \textbf{195} (2019), 104319. \url{//doi.org/10.1016/j.compfluid.2019.104319}

\bibitem{ann_repo}
A. Zingaro, ANN-SUPG, Project ID: 30854063, GitLab repository. Available from:
\url{//gitlab.com/albertozingaro/ann-supg}.

\bibitem{zingaro2021hemodynamics}
A. Zingaro, F. Menghini, L. Dede', A. Quarteroni,
  Hemodynamics of the heart’s left atrium based on a Variational Multiscale-LES numerical method,
  {\em Eur. J. Mech. B/Fluids}, \textbf{89} (2021), 380--400. \url{//doi.org/10.1016/j.euromechflu.2021.06.014}


\end{thebibliography}
\end{document}